\documentclass[11pt,usenames,dvipsnames]{amsart}

\usepackage{amsmath, amssymb, amsthm, array, bm, calrsfs, comment, enumitem, graphicx, hyperref, latexsym, mathpazo, tikz-cd}

\hypersetup{
colorlinks, 
linkcolor=darkred,
urlcolor=darkred,
citecolor=cadmiumgreen,
pdfauthor={Omid Amini, Daniel Corey, Leonid Monin}, 
pdftitle={tropical-Abel-Jacobi-theory},
pdfstartview ={FitV},
    }
\definecolor{cadmiumgreen}{rgb}{0.0, 0.42, 0.24}
\definecolor{darkred}{rgb}{.85,0,0}

\usepackage{manfnt}
\usepackage{scalerel}
\usepackage{xparse}
\usepackage[margin=1.3in]{geometry}
\usepackage[cmtip, all]{xy}
\usepackage{mathtools}
\usepackage{stmaryrd}

\RequirePackage[capitalize, noabbrev]{cleveref}
\newtheorem{theorem}{Theorem}[section]
\newtheorem{conjecture}[theorem]{Conjecture}
\newtheorem{lemma}[theorem]{Lemma}
\newtheorem{proposition}[theorem]{Proposition}

\newtheorem*{theorem*}{Theorem}

\newtheorem{question}[theorem]{Question}

\newtheorem*{question*}{Question}

\theoremstyle{definition}
\newtheorem{defii}[theorem]{Definition}
\newenvironment{definition}
  {\pushQED{\qed}\defii}
  {\popQED\enddefii}
\newtheorem{remm}[theorem]{Remark}
\newenvironment{remark}
  {\pushQED{\qed}\remm}
  {\popQED\endremm}

\numberwithin{equation}{section}
\numberwithin{table}{section}
\numberwithin{figure}{section}

\usepackage{todonotes} %% use disable option to disable

\newcommand\shorttitle[1]{\renewcommand\@shorttitle{#1}}

\newcommand{\Alb}{\mathrm{Alb}}
\newcommand{\bul}{\bullet}
\newcommand{\cc}{\mathbf{v}}
\newcommand{\CCt}{\C(\!(t)\!)}
\newcommand{\coker}{\mathrm{coker}}
\renewcommand{\div}{\mathrm{div}}
\newcommand{\Div}{\mathrm{Div}}

\newcommand{\Gr}{\mathrm{Gr}}
\newcommand{\image}{\mathrm{im}}
\newcommand{\Jac}{\mathrm{Jac}}
\newcommand{\JH}{\mathrm{JH}}
\newcommand{\Mod}{\mathrm{Mod}}
\newcommand{\Pic}{\mathrm{Pic}}
\newcommand{\Prin}{\mathrm{Prin}}

\newcommand{\reg}{\mathrm{reg}}
\newcommand{\sgn}{\ssub{\mathrm{sgn}}}
\newcommand{\sgnbar}{\ssub{\overline{\mathrm{sgn}}}}

\DeclareMathOperator{\Span}{span}
\newcommand{\ve}{\varepsilon}

\newcommand{\aab}[3]{(\mathsf{a}_{#1} \wedge \mathsf{a}_{#2}) \otimes \mathsf{b}_{#3}}
\newcommand{\abb}[3]{\mathsf{a}_{#1} \otimes (\mathsf{b}_{#2} \wedge \mathsf{b}_{#3})}
\newcommand{\bk}[2]{\langle #1, #2 \rangle}
\newcommand{\ev}[1]{\mathbf{I}_{#1}}
\newcommand{\rquot}[2]{#1\big/#2}

\newcommand{\src}[1]{{s}_{#1}}
\newcommand{\dst}[1]{{t}_{#1}}

\newcommand{\cycb}[1]{\sa_{#1}}
\newcommand{\unt}[1]{\sfb_{#1}}

\newcommand{\C}{\mathbb{C}}

\newcommand{\R}{\ssub{\mathbb{R}}}
\newcommand{\Z}{\ssub{\mathbb{Z}}}
\newcommand{\Q}{\mathbb{Q}}

\newcommand{\cC}{\mathcal{C}}
\newcommand{\cI}{\mathcal{I}}

\newcommand{\cM}{\mathcal{M}}
\newcommand{\cZ}{\mathcal{Z}}

\newcommand{\rmA}{\mathrm{A}}
\newcommand{\rmB}{\mathrm{B}}
\newcommand{\rmC}{\mathrm{C}}
\newcommand{\rmH}{\mathrm{H}}
\newcommand{\rmL}{\mathrm{L}}
\newcommand{\rmK}{\mathrm{K}}
\newcommand{\rmN}{\mathrm{N}}
\newcommand{\rmQ}{\mathrm{Q}}
\newcommand{\rmT}{\ssub{\mathrm{T}}}
\newcommand{\rmZ}{\mathrm{Z}}

\newcommand{\sa}{{\mathsf{a}}}
\newcommand{\sfb}{{\mathsf{b}}}

\newcommand{\fv}{\mathfrak{v}}

\makeatletter
\let\oldsum\sum
\renewcommand{\sum}{\@ifnextchar_\@mysum\oldsum}
\def\@mysum_#1{\oldsum_{\substack{#1}}}
\let\oldbigoplus\bigoplus
\renewcommand{\bigoplus}{\@ifnextchar_\@mybigoplus\oldbigoplus}
\def\@mybigoplus_#1{\oldbigoplus_{\substack{#1}}}
\let\oldprod\prod
\renewcommand{\prod}{\@ifnextchar_\@myprod\oldprod}
\def\@myprod_#1{\oldprod_{\substack{#1}}}
\makeatother

\let\oldbigwedge\bigwedge
\renewcommand{\bigwedge}{{\textstyle\oldbigwedge\!}}

\newcommand{\st}{\bigm|} % such that in sets
\newcommand{\Hom}{\ssub{\mathrm{Hom}}} % Hom
\newcommand{\rest}[1]{\raisebox{-1pt}{$\vert$}_{#1}}

\newcommand{\dual}{\star}
\newcommand{\conezero}{{\underline0}}

\newcommand{\suppaux}[2]{\scalebox{1}[1.2]{$#1\lvert$}#2\scalebox{1}[1.2]{$#1\rvert$}}
  \newcommand{\supp}[1]{\mathpalette\suppaux{#1}}
\newcommand{\dims}[1]{\ss{d}_{#1}}
\newcommand{\subface}{\prec}
\newcommand{\subfaceq}{\preceq}
\newcommand{\ssubface}{\mathbin{\mathchoice
  {\subface\!\!\!\cdot}%
  {\subface\!\!\!\cdot}%
  {\subface\!\cdot}%
  {\subface\!\cdot}%
}} % subface of codimension one codimension one
\newcommand{\supface}{\succ}
\newcommand{\supfaceq}{\succeq}
\newcommand{\ssupface}{\mathbin{\mathchoice
  {\cdot\!\!\!\supface}%
  {\cdot\!\!\!\supface}%
  {\cdot\!\supface}%
  {\cdot\!\supface}%
}}

\DeclareMathOperator{\sed}{sed} % sedentarity
\DeclareMathOperator{\SedT}{Sed} % sedentarity
\newcommand{\TT}{\mathrm{T}} % tangent space
\newcommand{\comp}[1]{\overline{#1}} % compactification

\newcommand{\RpMod}{\mathcal M} % category of \R_+ modules

\newcommand{\nvect}{\ssub{\mathfrak n}} % normal unit vector to gamma in delta

\newcommand{\eR}{\ssub{\mathbb T}} % tropical semiring
\newcommand{\bbU}{\ssub[-1pt]{\mathbb U}}
\newcommand{\bbO}{\ssub[-1pt]{\mathbb O}}

\newcommand{\TP}{\ssub{\T\P}} % for tropical toric variety
\newcommand{\T}{\mathbb T} % tropical semiring
\renewcommand{\P}{\mathbb P} % projective space
\newcommand{\e}{\ss{\mathfrak e}} % basis before quotienting
\newcommand{\trop}{{\scaleto{\mathrm{trop}}{5pt}}} % trop

\newcommand{\SF}{\textrm{\bf F}} % for multi-(co)tangent space
\newcommand{\AJ}{\ssub{\mathrm{AJ}}}
\newcommand{\ChowTriv}[2]{\rmA_{#1}^{\circ}(#2)}

\newcommand{\ie}{i.e.}
\newcommand{\eg}{e.g.}

\newcommand{\loccit}{\emph{loc. cit.}\ }
\newcommand{\ibid}{\emph{ibid.}\ }

\NewDocumentCommand{\ssub}{O{0pt} O{0pt} O{.8} m e{_^}}{
  #4%
  \IfValueT{#5}{
    \sb{\hspace{#1}\scaleobj{#3}{#5}}
  }
  \IfValueT{#6}{
    \sp{\hspace{#2}#6}
  }
}

\NewDocumentCommand{\ssup}{O{0pt} O{0pt} O{.8} m e{_^}}{
  #4%
  \IfValueT{#5}{
    \sb{\hspace{#1}#5}
  }
  \IfValueT{#6}{
    \sp{\hspace{#2}\scaleobj{#3}{#6}}
  }
}
\RenewDocumentCommand{\ss}{O{0pt} O{0pt} O{.9} m e{_^}}{
  #4%
  \IfValueT{#5}{
    \sb{\hspace{#1}\scaleobj{#3}{#5}}
  }
  \IfValueT{#6}{
    \sp{\hspace{#2}\scaleobj{#3}{#6}}
  }
}

\ExplSyntaxOn
\NewDocumentCommand{\tossub}{o o m}{
  \expandafter\let\csname old\cs_to_str:N #3\endcsname#3
  \renewcommand#3%
  {\ss[#1][#2]{\csname old\cs_to_str:N #3\endcsname}}
}
\ExplSyntaxOff

\tossub\omega

\tossub\varpi
\tossub\nu
\tossub\mu
\tossub\deg
\tossub{\iota}
\tossub{\pi}
\tossub{\Delta}
\tossub{\delta}
\tossub\phi

\newcommand{\sssigma}{\ssub[-2pt]{\sigma}}
\newcommand{\sstau}{\ssub[-2pt]{\tau}}
\newcommand{\ssinfty}{\ssub{\infty}}
\newcommand{\ssM}{\ss{M}}
\newcommand{\ssN}{\ss{N}}
\newcommand{\ssSigma}{\ssub{\Sigma}}
\newcommand{\compSigma}{\ssub{\comp\Sigma}}

\newcommand{\ssf}[1]{\ss f_{#1}}

\newcommand{\ssH}{\ss{H}}
\newcommand{\ssW}{\ss{W}}
\newcommand{\ssU}{\ss{U}}
\newcommand{\ssX}{\ss{X}}
\newcommand{\ssY}{\ss{Y}}

\newcommand{\ssdelta}{\ssub \delta}
\newcommand{\ssQ}{\ssub Q}
\newcommand{\ssF}{\ssub{F}}

\newcommand{\rflat}{\reflectbox{$\flat$}}
\newcommand{\rflatind}{\scaleto{\rflat}{6.6pt}}

\newcommand{\op}[1]{o_{\scaleto{#1}{5.5pt}}}
\newcommand{\sso}[1]{\ss o_{#1}}

\newcommand{\pibkk}[1]{\llbracket #1\rrbracket}
\newcommand{\bkk}[1]{\llbracket #1\rrbracket}

\newcommand{\cl}{\ss{\mathrm{cl}}}  %%%cycle class map

\newcommand{\cube}{\Pi}

\title{Tropical Abel--Jacobi theory}

\author[Omid Amini]{Omid Amini}
\address{Centre de math\'ematiques Laurent Schwartz \\
\'Ecole Polytechnique \\
91128 Palaiseau, France }
\email{\href{mailto:omid.amini@polytechnique.edu}{omid.amini@polytechnique.edu}}
\author[Daniel Corey]{Daniel Corey}
\address{University of Nevada, Las Vegas, 4505 S. Maryland Pkwy, Las Vegas NV 89154}
\email{\href{mailto:daniel.corey@unlv.edu}{daniel.corey@unlv.edu}}
\author[Leonid Monin]{Leonid Monin}
\address{\'Ecole Polytechnique F\'ed\'erale de Lausanne, 1015 Lausanne, Switzerland}
\email{\href{mailto:leonid.monin@mis.mpg.de}{leonid.monin@epfl.ch}}
\date{April 2025}


\begin{document}

\maketitle

\begin{abstract}

To a compact tropical variety of arbitrary dimension, we associate a collection of intermediate Jacobians defined in terms of tropical homology and tropical monodromy. 
We then develop an Abel--Jacobi theory in the tropical setting by defining functorial Abel--Jacobi maps. 
We introduce, in particular, tropical Albanese varieties and formulate obstructions to algebraic equivalence of tropical cycles. 
In dimension 1, we show that this recovers the existing Abel--Jacobi theory for tropical curves. 

As an application, we consider the Ceresa class of a tropical curve which is defined as the image of the Ceresa cycle in an appropriate intermediate Jacobian under the Abel--Jacobi map. 
We give an explicit formula for this class entirely in terms of the combinatorics of the tropical curve. 
\end{abstract}

\setcounter{tocdepth}{1}
\tableofcontents
\section{Introduction}

Associated to a smooth and projective variety $X$ of dimension $d$ is a collection of intermediate Jacobians $\JH_{2k+1}(X)$, for $0\leq k\leq d-1$, which were introduced by Griffiths in \cite{Griffiths-68}.  
These are complex tori that interpolate between the well-known abelian varieties associated with $X$: the Picard variety $\Pic(X) = \JH_{2d-1}(X)$ and the Albanese $\Alb(X) = \JH_{1}(X)$. 
For algebraic curves these two abelian varieties are naturally isomorphic, and this common variety is called the Jacobian and denoted by $\Jac(X)$. 
Furthermore, the Abel--Jacobi map from the theory of algebraic curves generalizes to higher-dimensional varieties as a map
\begin{equation*}
	\AJ\colon \rmA_{k}^{\circ}(X) \to \JH_{2k+1}(X)
\end{equation*}
where $\rmA_{k}^{\circ}(X)$ is the Chow group of homologically trivial $k$--dimensional cycles in $X$ modulo rational equivalence. 
Therefore, intermediate Jacobians have been used to distinguish between homological, algebraic, and rational equivalence of cycles \cite{Ceresa, Clemens, Griffiths-69, Voisin-00}.

The theory of tropical curves has seen intense development in recent years out of a desire to understand the asymptotic geometry of a degenerating family of algebraic curves.   
Much of the theory described above for algebraic curves exists in an analogous way for tropical curves. 
Notably, Mikhalkin and Zharkov in \cite{MikhalkinZharkov:TropicalCurves} associate to a tropical curve $\Gamma$ its Jacobian $\Jac(\Gamma)$ and an Abel--Jacobi map
\begin{equation*}
	\AJ\colon \Gamma \to \Jac(\Gamma).
\end{equation*}
A version of this construction for finite graphs appears in \cite{BacherHarpeNagnibeda, BakerNorine}. 

A degenerating family of algebraic curves may be described succinctly as a (smooth and complete) curve $X$ over a field $K$ equipped with a nontrivial nonarchimedian norm. 
To such a curve, we can associate a tropical curve $\Gamma$, by taking a \emph{skeleton} of the Berkovich analytification of $X$, see \cite{Berkovich, BakerPayneRabinoff}. Furthermore, the Jacobian and Abel--Jacobi map of $\Gamma$ are related in a similar way to the Jacobian and Abel--Jacobi map of $X$ \cite{BakerRabinoff, AminiNicolussi25}. 
Given this connection, tropical geometry incorporates combinatorial and computational techniques to cast new light on old theorems and to derive new results in algebraic geometry. We refer to the survey papers \cite{BJ16, JP21} for an account of recent development in this direction.

\medskip

The aim of this paper is to initiate the development of an Abel--Jacobi theory for tropical varieties of arbitrary dimension. 

Let $X$ be a $d$-dimensional tropical variety. For $A = \Z$ or $\R$, denote by $\rmH_{p,q}(X,A)$ the $(p,q)$--tropical homology of $X$ for $0\leq p,q \leq d$, introduced in~\cite{ItenbergKatzarkovMikhalkinZharkov}.  
The tropical monodromy map is a linear map $\rmN \colon \rmH_{q,p}(X,A) \to \rmH_{q+1,p-1}(X, \R)$. 
Let
\begin{equation*}
	\rmL_{p,q}(X) \coloneqq \image(\rmN^{p-q} \colon \rmH_{q,p}(X,\Z) \to \rmH_{p,q}(X,\R)).
\end{equation*}
For $p\geq q$, we define the \emph{$(p,q)$--th (tropical) intermediate Jacobian of $X$} by
\begin{equation*}
	\JH_{p,q}(X) \coloneqq \frac{\rmH_{p,q}(X,\R)}{\rmL_{p,q}(X)}.
\end{equation*}

In general, $\JH_{p,q}(X)$ is noncompact, but if $X$ satisfies the \emph{weight-monodromy property} as in \eqref{eq:WMP}, \eg, if $X$ is a projective tropical manifold or more generally if $X$ is a K\"ahler tropical variety, then $\JH_{p,q}(X)$ is a real torus; see \S~\ref{sec:KahlerWMP} for details. 
In \cite{MikhalkinZharkov:Eigenwave}, Mikhalkin and Zharkov define a tropical intermediate Jacobian for projective tropical manifolds when $p+q=d$, see Remark \ref{rem:MZIntermediateJacobian} for a comparison to our definition in this setting.

\medskip

Denote by $\rmA_p(X)$ the Chow group of $p$-dimensional tropical cycles of $X$ modulo rational equivalence, see \S~\ref{sec:rationalFunctionsDivisors}. There is a cycle class map:
\begin{equation*}
	\cl_{X} \colon \rmA_{p}(X) \to \rmH_{p,p}(X,\Z).
\end{equation*}
Define $\ChowTriv{p}{X}$ to be the kernel of the cycle class map, \ie, the group of \emph{homologically trivial $p$-cycles modulo rational equivalence}. 
Our first result is the following theorem.

\begin{theorem}
\label{thm:IntroAJ}
There is a well-defined (tropical) Abel--Jacobi map
\begin{equation*}
	\AJ\colon \ChowTriv{p}{X} \to \JH_{p+1,p}(X).
\end{equation*}
It is functorial in the following sense. 
For any morphism $\varphi \colon X \to X'$ of tropical varieties, we have the following commutative diagram:
\begin{equation*}
	\begin{tikzcd}
	\ChowTriv{p}{X} \arrow[r, "\AJ"]  \arrow[d, "\varphi_{*}"] & \JH_{p+1,p}(X) \arrow[d, "\varphi_{*}"]  \\
	\ChowTriv{p}{X'} \arrow[r, "\AJ"]   & \JH_{p+1,p}(X')\mathrlap{.}
	\end{tikzcd}
\end{equation*}
\end{theorem}
We briefly describe our construction of the Abel--Jacobi map. 
Given a homologically trivial cycle in $\ChowTriv{p}{X}$, we choose a $(p,p+1)$--chain $\gamma$ that bounds this cycle. 
We apply the monodromy map to $\gamma$ and show that, by pairing with $(p+1,p)$--cocycles, it yields a well defined element in $\JH_{p+1,p}(X)$. 
Furthermore, we show that the Abel--Jacobi map vanishes on cycles rationally equivalent to 0. 
See \S~\ref{sec:JH_and_AJ}. 

Our first application of this construction is to find an obstruction to algebraic equivalence between tropical cycles. 
We define the following subgroup
\begin{equation*}
	\rmK_{p,q}(X) \coloneqq \image(\rmN^{p-q-1}\colon \rmH_{q+1,p-1}(X,\Z) \to \rmH_{p,q}(X,\R)).
\end{equation*}
Set 
\begin{equation*}
	\rmQ_{p,q}(X) \coloneqq \frac{\rmH_{p,q}(X,\R)}{\rmK_{p,q}(X)}.
\end{equation*}
If the monodromy map $\rmN$ is integral, then $\rmK_{p,q}(X)$ contains $\rmL_{p,q}(X)$, but in general the two spaces are unrelated. 
In any case, $\rmN \colon \rmH_{p,q}(X,\R)\to \rmH_{p+1,q-1}(X,\R)$ induces a map
\begin{equation*}
	\rmN \colon \JH_{p,q}(X) \to \rmQ_{p+1,q-1}(X).
\end{equation*}

\newtheorem*{thm:obstructionAlgebraicEquivalence}{\cref{thm:obstructionAlgebraicEquivalence}}
\begin{thm:obstructionAlgebraicEquivalence}
	\emph{If a $p$-cycle $\alpha$ is algebraically equivalent to 0, then 
	\begin{equation*}
		\rmN(\AJ(\alpha)) = 0 \quad \text{in} \quad \rmQ_{p+2,p-1}(X).
	\end{equation*}
	}
\end{thm:obstructionAlgebraicEquivalence}

\noindent In the case of tropical curves, Zharkov \cite{Zharkov} and Ritter \cite{Ritter} describe similar obstructions for tropical Ceresa cycles to be algebraically equivalent to 0, see \S~\ref{sec:CZ-class}. 

\medskip

Next, we define the Albanese of a K\"ahler tropical variety $X$ as $\Alb(X) = \JH_{1,0}(X)$. 
This is a  tropical abelian variety when the monodromy of $X$ is rational. 
We prove that $\Alb(X)$ satisfies a universal property similar to that in the classical setting, see Proposition \ref{prop:albaneseUniversalProperty}. 

In the case of tropical abelian varieties, or more generally tropical compact tori, we provide an explicit description of the monodromy map, yielding a concrete  description of the intermediate Jacobians.

\medskip

We now discuss the case of tropical curves and their Jacobians. 
A \emph{tropical curve} is a tropical variety of dimension 1, which roughly speaking means that it is a compact metric graph endowed with an integral affine structure. 
Recall that a metric graph is a metric space associated to a pair $(G,\ell)$ where $G = (V,E)$ is a finite graph  and $\ell \colon E\to \R_{>0}$ is an edge-length function. 
Precisely, we plug in an interval of length $\ell(e)$ between the two endpoints of each edge $e$, and the metric is the path metric. 
Different pairs might yield the same tropical curve, a fixed pair $(G,\ell)$ is called a \emph{model} of $\Gamma$.

Any tropical curve is K\"ahler. 
We show that $\JH_{1,0}(\Gamma)$ is naturally isomorphic to the  Jacobian $\Jac(\Gamma)$ of $\Gamma$ and via this isomorphism, the above $\AJ \colon \ChowTriv{0}{\Gamma} \to \JH_{1,0}(\Gamma)$ is the Abel--Jacobi map from \cite{MikhalkinZharkov:TropicalCurves}. See Theorem \ref{thm:tropicalJHCurves}. 

Classically, a natural and important instance of a homologically trivial cycle is the Ceresa cycle of an algebraic curve. Let $X$ be an algebraic curve. 
To a point  $\flat \in X$ corresponds two maps:
\begin{equation}
\label{eq:XToDiv0}
	X \to \Div^{0}(X), \quad x\mapsto [x] - [\flat] \quad \text{ and } \quad x\mapsto  [\flat] -[x].
\end{equation}
Composing with the Abel--Jacobi map defines two embedding
\begin{equation*}
	X \hookrightarrow \Jac(X).
\end{equation*}
The images of $X$ determine two dimension-$1$ cycles $[X_{\flat}]$ and $[X_{\flat}^{-}]$ which have the same homology class. 
The \emph{Ceresa cycle} is the algebraic cycle $[X_{\flat}] - [X_{\flat}^{-}]$ in $\ChowTriv{1}{\Jac(X)}$. 
Ceresa's theorem \cite{Ceresa} asserts that the cycle $[X_{\flat}] - [X_{\flat}^{-}]$ is not algebraically equivalent to 0 for a very general curve of genus at least 3 in the corresponding moduli space.

In the tropical setting, the Ceresa cycle of a connected tropical curve $\Gamma$ is defined in an analogous way. 
First, the Jacobian $\Jac(\Gamma)$ of $\Gamma$ is a K\"ahler tropical variety of dimension $g$, where $g$ is the first Betti number of $\Gamma$; this is more commonly referred to as the \emph{genus} of $\Gamma$.
Using the tropical analog of the maps in \eqref{eq:XToDiv0}, a point $\flat \in \Gamma$ yields two homologous $1$-cycles $[\Gamma_{\flat}]$ and $[\Gamma_{\flat}^{-}]$. 
The Ceresa cycle relative to the basepoint $\flat$ is the $1$-cycle $[\Gamma_{\flat}] - [\Gamma_{\flat}^{-}]$. 
This is homologically trivial and it gives rise, via  the Abel--Jacobi map, to the Ceresa class $\cc_{\flat}(\Gamma)$ of $\Gamma$ in $\JH_{2,1}(\Jac(\Gamma))$. Our next theorem provides an explicit description of $\cc_{\flat}(\Gamma)$.

As in the classical setting, $\rmH_{p,q}(\Jac(\Gamma),\R)$ can be expressed in terms of the homology of $\Gamma$:
\begin{equation}
\label{eq:tropHomologyJac}
	\rmH_{p,q}(\Jac(\Gamma), \R) \cong \wedge^{q} \rmH_{0,1}(\Gamma,\R) \otimes \wedge^{p} \rmH_{1,0}(\Gamma,\R).
\end{equation}
 Explicitly, $\rmH_{0,1}(\Gamma, \R)$ is the first singular homology group $\rmH_{1}(\Gamma,\R)$ of $\Gamma$ and $\rmH_{1,0}(\Gamma,\R)$ is the dual of $\rmH_{0,1}(\Gamma,\R)$ with respect to the Poincar\'e duality pairing
 \begin{equation}
 \label{eq:tropicalCurvePoincarePairing}
 	\cap \colon \rmH_{0,1}(\Gamma,\R) \times \rmH_{1,0}(\Gamma,\R) \to \rmH_{0,0}(\Gamma,\R) = \R.
 \end{equation}
 Let $(G, \ell)$ be a model of $\Gamma$ with $G = (V,E)$.  Fix an orientation of each edge; denote by $\src{e}$ the source of $e$ and by $\dst{e}$ the target of $e$.  
 Choosing a spanning tree $T = (V,F)$, we obtain a basis $\{ \cycb{\ve} \st \ve \in \ssF^c \}$ of $\rmH_{0,1}(\Gamma,\R)$ where $\ssF^c = E\setminus F$. 
 Let $\{\unt{\ve} \st \ve \in F^{c}\}$ be the dual basis in $\rmH_{1,0}(\Gamma,\R)$ with respect to the pairing in \eqref{eq:tropicalCurvePoincarePairing}. 
 For each edge $e \in F$, we define $\unt{e} \in \rmH_{1,0}(\Gamma, \R)$. 
 The removal of $e$ from the spanning tree separates $T$ into two connected components, one connected component, denoted by $S_1$, contains $\src{e}$ and the other, denoted by $S_2$, contains $\dst{e}$.  
 Then,
 \begin{equation}
 \label{eq:generalbe}
 	\unt{e} = \sum_{\src{\ve} \in S_2 \\\dst{\ve} \in S_1} \unt{\ve}
 	-
 	\sum_{\src{\ve} \in S_1 \\\dst{\ve} \in S_2 } \unt{\ve}.
 \end{equation}
\noindent 
Geometrically, $\unt{e}$ is the unit tangent vector to the oriented edge $e$ of $\Gamma$ inside $\Jac(\Gamma)$. 
We also have a map $\sgn_{T}^{\flat}\colon F \times \ssF^c \to \{0, \pm 1, \pm 2\}$ which, roughly speaking, records the positions of each pair of edges $e\in F$ and $\ve\in \ssF^c$ relative to $\flat$, see \S~\ref{sec:Ceresa-class-explicit}.

\newtheorem*{thm:ceresa-explicit}{\cref{thm:ceresa-explicit}}
\begin{thm:ceresa-explicit}
\emph{Given a tropical curve $\Gamma$ with model $(G = (V,E),\ell)$, a point $\flat\in V$, and spanning tree $T = (V, F)$, the tropical Ceresa class $\cc_{\flat}(\Gamma)$ is given by
\begin{equation*}
	\cc_{\flat}(\Gamma) = \sum_{e\in F \\ \ve \in \ssF^c} \sgn_{T}^{\flat}(e,\ve) \, \ell(e) \, \cycb{\ve} \otimes (\unt{\ve} \wedge \unt{e}). 
\end{equation*}
}
\end{thm:ceresa-explicit}

Next we consider the dependence on the base point $\flat$. 
Define the homology class
\begin{equation*}
	\omega = \sum_{\ve \in \ssF^c} \cycb{\ve} \otimes  \unt{\ve}
	\quad \text{ in } \quad
	\rmH_{0,1}(\Gamma,\R) \otimes \rmH_{1,0}(\Gamma,\R) \cong \rmH_{1,1}(\Jac(\Gamma),\R).
\end{equation*}
Given two points $\flat, \rflat \in \Gamma$, we prove in Theorem \ref{thm:dependenceOnBasePoint} that their Ceresa classes are related by the simple formula
\begin{equation}
\label{eq:dependenceOnBasePoint}
	\cc_{\flat}(\Gamma) - \cc_{\rflatind}(\Gamma) = -2 \AJ(\flat - \rflat) \wedge \omega \quad \text{ in } \quad \JH_{2,1}(\Jac(\Gamma)).
\end{equation}
This is a tropical analog of a formula provided by Hain and Reed in \cite[p.204]{HainReed} see also \cite{Pulte}. 
We define
\begin{equation*}
	\overline{\JH}_{p+1,p}(\Jac(\Gamma)) \coloneqq \frac{\rmH_{p+1,p}(\Jac(\Gamma),\R)}{\omega \wedge \rmH_{p,p-1}(\Jac(\Gamma),\R) + \rmL_{p+1,p}(\Jac(\Gamma))}.
\end{equation*}
The image of $\cc_{\flat}(\Gamma)$ under the natural quotient map $\JH_{2,1}(\Jac(\Gamma)) \to \overline{\JH}_{2,1}(\Jac(\Gamma))$ is independent of $\flat$ by Equation \eqref{eq:dependenceOnBasePoint}. 
We call this class the \emph{unpointed Ceresa class of $\Gamma$}, and denote it by $\overline{\cc}(\Gamma)$. 
We have a formula for $\overline{\cc}(\Gamma)$ similar to that of $\cc_{\flat}(\Gamma)$ appearing in Theorem \ref{thm:ceresa-explicit}.

\begin{figure}	
\includegraphics[width=0.7\textwidth]{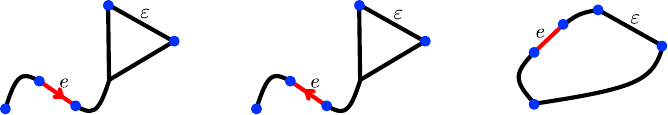}
	\caption{The values of $\sgnbar_{T}(e,\ve)$ are, left to right: $1$, $-1$, $0$}
	\label{fig:sgn_bar}
\end{figure}

We define a sign function $\sgnbar_{T} \colon F \times \ssF^c \to \{0,\pm 1\} $ in the following way, see Figure \ref{fig:sgn_bar}. 
As before, an edge $e\in F$ separates the spanning tree $T$ into two connected components: $S_1$, which contains $\src{e}$,  and $S_2$, which contains $\dst{e}$. 
Given another edge $\ve \in \ssF^c$, we set
\begin{equation*}
	\sgnbar_{T}(e,\ve) =
	\begin{cases}
		1 & \text{ if } \ve \in S_1, \\
		-1 & \text{ if } \ve \in S_2, \\
		0 & \text{ otherwise.}
	\end{cases}
\end{equation*}

\newtheorem*{thm:unpointedCeresaFormula}{\cref{thm:unpointedCeresaFormula}}
\begin{thm:unpointedCeresaFormula} \emph{Given a tropical curve $\Gamma$ with model $(G = (V,E), \ell)$ and spanning tree $T = (V,F)$, its unpointed tropical Ceresa class is 
\begin{equation*}
	\overline{\cc}(\Gamma) = \sum_{e\in F \\ \ve \in \ssF^c}  \sgnbar_{T}(e,\ve) \, \ell(e) \, \cycb{\ve} \otimes (\unt{\ve}\wedge \unt{e}) \quad \text{ in }\quad \overline{\JH}_{2,1}(\Jac(\Gamma)).
\end{equation*}
}
\end{thm:unpointedCeresaFormula}
\noindent When $\Gamma$ is hyperelliptic and $\flat$ is a Weierstrass point, the Ceresa cycle $[\Gamma_{\flat}] - [\Gamma_{\flat}^{-}]$ equals $0$, and so $\overline{\cc}(\Gamma) = 0$.

In the case where the edge lengths of $\Gamma$ are positive integers, the unpointed tropical Ceresa class relates to the tropical Ceresa class defined in \cite{CoreyEllenbergLi}, which in turn relates to the $\ell$-adic Ceresa of an algebraic curve defined over $\CCt$ as studied in \cite{HainMatsumoto}. 
In \cite{CoreyEllenbergLi}, there is a finite group $\overline{B}(\ssdelta_{\Gamma})$ and the tropical Ceresa class $n(\Gamma)$ is an element of $\overline{B}(\ssdelta_{\Gamma})$ which is defined using the Johnson homomorphism on the Torelli group of a topolocial surface (more accurately, the extension of this homomorphism to the mapping class group by Morita \cite{Morita:Extension-Johnson}).

In order to relate these two classes, we discovered a new conjectural formula for the Johnson homomorphism which we think could be of independent interest. 
We briefly describe the formula and refer to Appendix \ref{sec:MoritaClass} for details.

Let $\ssSigma_{g}$ be an orientiable, compact, topological surface of genus $g\geq 3$. 
For a simple closed curve $\gamma$, denote by $\rmT_{\gamma}$ the (left-hand) Dehn-twist about $\gamma$. 
Suppose there are nonseparating simple curves $\gamma, \gamma'$ and $\beta_1,\ldots,\beta_{g}$ such that 
\begin{enumerate}
	\item $\gamma$ and $\gamma'$ are homologous and disjoint from each $\beta_i$, and
	\item the curves $\beta_1,\ldots, \beta_g$  are pairwise disjoint and their homology classes form a basis for a Lagrangian subspace of $\rmH_1(\ssSigma_g, \Z)$.
\end{enumerate}

\begin{figure}
	\includegraphics[scale=.4]{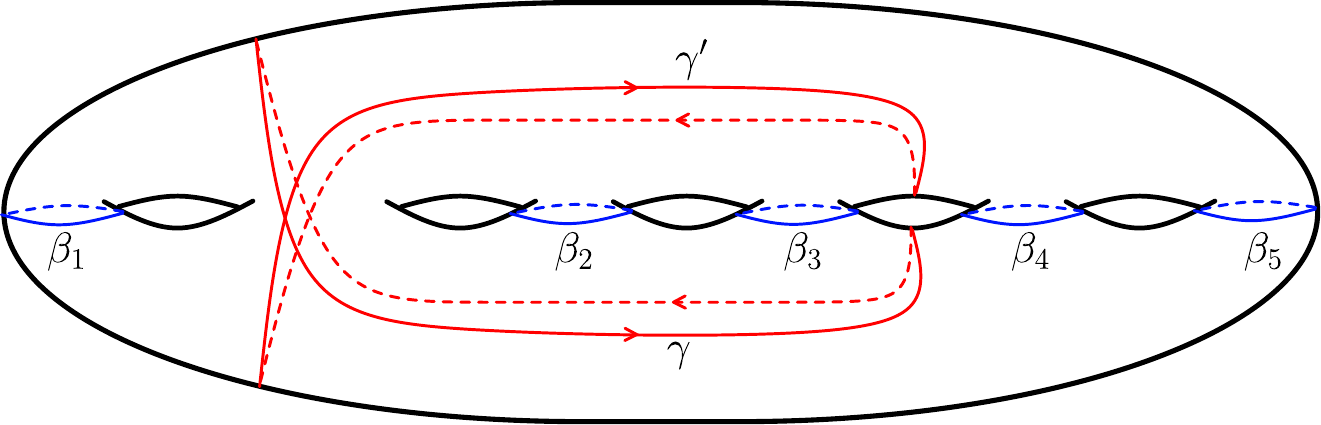}
	\caption{For this arrangement of curves, the values of $\sgn_{\gamma,\gamma'}$ on $\beta_1$, $\beta_2$, $\beta_3$, $\beta_4$, $\beta_5$ are $-1$, $1$, $1$, $0$, $0$, respectively.}
	\label{fig:sgn-surface}
\end{figure}

We define a sign function $\sgn\colon \{\beta_1,\ldots, \beta_g\} \to \{0,\pm 1\}$ that records the positions of the $\beta_i$'s relative to $\gamma$ and $\gamma'$, see Figure \ref{fig:sgn-surface} and \S~\ref{app:Comparison}. 
Finally, suppose $\alpha_1, \ldots, \alpha_g$ are curves such that the homology classes $[\alpha_1]$, $\ldots$, $[\alpha_g]$, $[\beta_1]$, $\ldots$, $[\beta_g]$ form a symplectic basis of $\rmH_{1}(\ssSigma_g, \Z)$. 

\newtheorem*{conj:Johnson}{\cref{conj:Johnson}}
\begin{conj:Johnson}
\emph{The image of the Johnson homomorphism on the mapping class $\rmT_{\gamma} \rmT_{\gamma'}^{-1}$ is
\begin{equation*}
	\sum_{i=1}^{g} \sgn_{\gamma,\gamma'}(\beta_i) \, [\alpha_{i}] \wedge [\beta_{i}] \wedge [\gamma].
\end{equation*}
}
\end{conj:Johnson}
 
\noindent When $\gamma$ and $\gamma'$ are disjoint, this recovers the well-known formula for the Johnson homomorphism of a bounding pair map as described in \cite[\S~6.6.2]{FarbMargalit}. 
If this conjecture holds, we obtain an explicit formula for the Johnson homomorphism of mapping classes of the  form $\rmT_{\gamma} \rmT_{\gamma'}^{-1}$ which are not necessarily bounding pair maps.

We prove in Appendix \ref{sec:MoritaClass} that for a tropical curve $\Gamma$ with integer edge-lengths, there is an embedding $\Phi_{\Gamma}$ of  $\overline{B}(\ssdelta_{\Gamma})$ into the torsion part of $\overline{\JH}_{2,1}(\Jac(\Gamma))$. 
Furthermore, we show the following comparison result.
\newtheorem*{thm:Morita-to-Ceresa}{\cref{thm:Morita-to-Ceresa}}
\begin{thm:Morita-to-Ceresa}
	\emph{Let $\Gamma$ be a tropical curve with integral edge lengths. 
	If Conjecture \ref{conj:Johnson} holds, then we have
	\begin{equation*}
		\Phi_{\Gamma}(n(\Gamma)) = \overline{\cc}(\Gamma).
	\end{equation*}
	}
\end{thm:Morita-to-Ceresa}

\subsection*{Acknowledgements} This project began in Michael Joswig's Diskrete Mathematik / Geometrie group at Technische Universit\"at Berlin. 
We thank them for providing a warm and collaborative environment.  O.A. thanks Math+, the Berlin Mathemaics Research Center for support.
We also thank the Centre de math\'ematiques Laurent Schwartz at \'Ecole polytechnique for hosting us in the late stages of this project. 
We thank Andrew Putman for providing us with examples that illustrate the difficulty in establishing Conjecture \ref{conj:Johnson}. 

O.A. was partially supported by ANR project AdAnAr (ANR-24-CE40-6184).
D.C. was partially supported by SFB-TRR project ``Symbolic Tools in Mathematics and their Application'' (project- ID 286237555). 
L.M. was partially supported by SNSF grant -- 200021E\_224099.

\section{Preliminaries on tropical varieties}
In this section, we gather basic definitions and results about tropical varieties. 
These constructions are known in the literature, and our presentation is close to~\cite{JellShawSmacka, JellRauShaw,AminiPiquerez:Hodge}, to which we refer for more details.

\subsection{Fans} Let $N$ be a free abelian group of finite rank and denote its dual by $M=\ssN^\dual = \Hom(N, \Z)$. Let $\ssN_{\Q}$, $\ssN_{\R}$, $\ssM_{\Q}$, $\ssM_{\R}$ be the corresponding rational and real vector spaces, we thus have $\ssM_\Q = \ssN_\Q^\dual$ and $\ssM_{\R} = \ssN_{\R}^\dual$.

 If $\sigma$ is a polyhedral cone in $\ssN_{\R}$, we denote by $\ssN_{\sigma, \R}$ the real vector subspace of $\ssN_{\R}$ generated by points of $\sigma$, and set $\ssN_{\R}^{\sigma} \coloneqq \rquot{\ssN_{\R}}{\ssN_{\sigma, \R}}$. 
 We denote by  $\dims\sigma$ the dimension of $\sigma$ which is equal to that of $\ssN_{\sigma, \R}$. 
 If $\sigma$ is rational, we also get natural full rank lattices $\ssN^{\sigma} \subset \ssN_{\R}^\sigma$ and $\ssN_{\sigma} \subset \ssN_{\sigma, \R}$.

A polyhedral cone in $\ssN_{\R}$ is called \emph{strongly convex} if it does not contain any line. 
A \emph{fan} $\Sigma$ of dimension $d$ in $\ssN_{\R}$ is a finite, non-empty, collection of strongly convex polyhedral cones in $\ssN_{\R}$ which verifies the following two properties:
\begin{itemize}
\item[(1)] for any cone $\sigma \in \Sigma$, any face $\tau$ of $\sigma$ belongs to $\Sigma$, and
\item[(2)] for any pair of cones $\sigma, \eta \in \Sigma$, the intersection $\sigma \cap \eta$ is a common face of $\sigma$ and $\eta$.
\end{itemize}
Given two cones $\sigma$ and $\tau$, we write $\tau \subfaceq \sigma$ if $\tau$ is a face of $\sigma$, and we write $\tau \ssubface \sigma$ if $\tau$ is a codimension one face of $\sigma$. 
This defines a partial order on the set of cones of $\Sigma$. 
The fan $\Sigma$ is \emph{rational} if each cone $\sigma \in \Sigma$ is rational. All the fans which appear in this paper are assumed to be rational.

We denote by $\ssSigma_k$ the set of $k$-dimensional cones of $\Sigma$. 
 We call elements of $\ssSigma_1$ rays of $\Sigma$. 
The cone $\{0\}$ is denoted by $\conezero$. 
Note that any $k$-dimensional cone $\sigma$ in $\Sigma$ is determined by its set of rays in $\ssSigma_1$.
A fan $\Sigma$ is called \emph{pure dimensional} if all its maximal cones have the same dimension.

We denote by $\supp \Sigma$ the \emph{support} of $\Sigma$, which by definition is the union of the cones of $\Sigma$. 
A \emph{fanfold} $X$ in $\ssN_{\R}$ is a closed subset which is the support of a rational fan in $\ssN_{\R}$.

\subsection{Normal vectors}

Let $\Sigma$ be a rational fan of pure dimension $d$. For a cone $\sigma$ in $\Sigma$ and a face $\tau$ of codimension one in $\sigma$, the space $\ssN_{\tau,\R}$ divides $\ssN_{\sigma,\R}$ into two closed half-spaces.  
One of these half spaces contains $\sigma$, and by a \emph{unit normal vector to $\tau$ in $\sigma$} we mean any vector $v\in \ssN_\sigma$ which lies in this half space, and which satisfies $\ssN_\tau + \Z v = \ssN_\sigma$. 
We denote a choice of such an element by $\nvect_{\sigma/\tau}$ in this paper. 
The corresponding vector in the quotient lattice $\ssN_{\sigma}^{\tau} \coloneqq \rquot{\ssN_{\sigma}}{\ssN_{\tau}}$  is well-defined and is denoted by $\e^\tau_\sigma$.

\subsection{Tropical fans and tropical fanfolds} A \emph{tropical fan} is a pair $(\Sigma,\omega)$ such that $\Sigma$ in $N_{\R}$ is a rational fan of pure dimension $d$, for a natural number $d\in \Z_{\geq0}$, and $\omega\colon \ssSigma_d \to \Z\setminus \{0\}$ is a \emph{tropical weight function}. 
This means that $\omega$ satisfies the \emph{balancing condition}: for any cone $\tau$ of dimension $d-1$ in $\Sigma$, we have
\begin{align}\label{eq:codimension_one_balancing} \sum_{\sigma \ssupface \tau} \omega(\sigma)\nvect_{\sigma/\tau} \in N_\tau, \qquad
\text{equivalently} \qquad
\sum_{\sigma \ssupface \tau} \omega(\sigma)\e^\tau_{\sigma} =0 \quad \textrm{in $N^\tau$}.
\end{align} 
The tropical fan is \emph{effective} if $\omega(\sigma)>0$ for all $\sigma\in \ssSigma_d$.
The support of a tropical fan is called a \emph{tropical fanfold}.

\subsection{Canonical compactifications of fans and fanfolds} A rational fan $\Sigma$ in $\ssN_{\R}$ gives rise to the tropical toric variety $\TP_\Sigma$ which is a partial compactification of $\ssN_{\R}$~\cite{BGJK, Kaj08, OR11, Pay09, Thu07}. 
This partial compactification coincides with the tropicalization of the toric variety $\P_\Sigma$ associated to $\Sigma$. 
The canonical compactification $\comp \Sigma$ of $\Sigma$ is defined as the closure of $\Sigma$ in $\TP_\Sigma$. Its support is the canonical compactification of the fanfold $X=\supp{\Sigma}$ with respect to $\Sigma$. 
We briefly discuss the construction and the induced stratification.

We denote by $\eR \coloneqq \R \cup \{\infty\}$ the extended real line endowed with its natural topology that extends the topology $\R$ by adding a basis of open neighborhoods of $\infty$ consisting of half open intervals $(a, \infty]$, $a\in \R$. 
The addition of real numbers extends to $\eR$. The resulting monoid $(\eR, +)$ is called the monoid of tropical numbers. Let $\eR_+ \coloneqq \R_+ \cup\{\infty\}$ be the submonoid of non-negative tropical numbers.
Denote by $\RpMod$ the category of $\R_+$-modules.

First we consider the tropical toric variety associated to a rational polyhedral cone $\sigma$ in $\ssN_{\R}$. The \emph{dual cone} $\sssigma^\vee$ and \emph{orthogonal plane} $\sssigma^\perp$ are defined by
\begin{equation*}
\sssigma^\vee \coloneqq 
\Bigl\{u\in \ssM_{\R} \st \langle u, v \rangle \geq 0 \ \forall\ v \in \sigma\Bigr\}, 
\quad \textrm{ and } \quad 
\sssigma^\perp \coloneqq \Bigl\{u \in \ssM_{\R} \st \langle u, v \rangle = 0 \ \forall\ v \in \sigma\Bigr\}.
\end{equation*}
Let $\bbU_\sigma \coloneqq \Hom_{\RpMod}(\sssigma^\ve e, \eR)$. 
The topology on $\eR$ induces a natural topology on $\bbU_\sigma$. 
Furthermore, the natural pairing $\langle u, v \rangle$ on $\ssM_{\R}\times \ssN_{\R}$ extends to a pairing $\ssM_{\R}\times \bbU_{\sigma} \to \eR$. 
Given a face $\tau$ of $\sigma$, the \emph{tropical torus orbit} of $\tau$ is
\begin{equation*}
	\mathbb{O}(\tau) = \{v\in \bbU_{\sigma} \st \langle u, v \rangle = \infty \;\; \text{ exactly for } u \in \sssigma^{\vee} \setminus \tau^{\perp}  \}.
\end{equation*}
The subspaces $\{\bbO(\tau) \st \tau \subfaceq \sigma \}$ form a  stratification of $\bbU_{\sigma}$.   
Moreover, $\bbO(\tau)$ may be identified with $\ssN_{\R}^{\tau}$ in the following way. 
Denote by $\infty_{\tau}$ the element of $\bbU_{\sigma}$ given by
\begin{equation*}
	\ss \infty_{\tau}(u) = 
	\begin{cases}
	0 & \text{ if } \; u\in \sstau^{\perp}, \\
	\infty & \text{ if }\; u\in \sssigma^\vee \setminus \sstau^{\perp}.
	\end{cases}
\end{equation*}
Under the natural identification 
$\ssN_{\R}^{\tau} = \Hom_{\RpMod}(\sstau^\perp, \R)$, we have
\begin{equation*}
	\bbO(\tau) = \ss \infty_{\tau} + \ssN_{\R}^{\tau} \subset \bbU_{\sigma}. 
\end{equation*}

Viewing $\tau$ as a cone by itself, we may also form $\bbU_{\tau}$ and we have an inclusion $\bbU_{\tau} \subset \bbU_{\sigma}$. 
The \emph{tropical toric variety} $\TP_\Sigma$ of a rational polyhedral fan $\Sigma$ in $\ssN_{\R}$ is obtained by gluing $\bbU_\sigma$, for  $\sigma\in \Sigma$, along these inclusions. 
If $\eta \subfaceq \tau \subfaceq\sigma$, then we may form $\bbO(\eta)$ in either $\bbU_{\tau}$ or $\bbU_{\sigma}$, but they are naturally identified under the inclusion $\bbU_{\tau} \subset \bbU_{\sigma}$. 
Therefore, we have a stratification into tropical torus orbits
\begin{equation*}
	\TP_{\Sigma} = \bigsqcup_{\sigma\in \Sigma} \bbO(\sigma).
\end{equation*}

The \emph{canonical compactification} of $\sigma$, denoted by $\comp{\sigma}$ is the closure of $\sigma$ in $\bbU_{\sigma}$. 
Under the natural identification $\sigma = \Hom_{\RpMod}(\sssigma^{\vee}, \R_+) \subset \bbU_\sigma$, we have
\begin{equation*}
	\comp{\sigma} = \Hom_\RpMod(\sssigma^\vee, \eR_+) \subset \bbU_\sigma.
\end{equation*}
Similarly, the \emph{canonical compactification} $\comp{\Sigma}$ of $\Sigma$ is defined as the closure of  $\Sigma$ in $\TP_{\Sigma}$. 
We can explicitly describe the polyhedral structure on $\comp{\Sigma}$ in the following way. Consider the cones $\tau \subface \sigma$. 
We define 
\begin{equation*}
\ss \sigma_{\infty}^{\tau}  \coloneqq \ss \infty_{\tau} + \sigma
	= \{v \in \comp{\sigma} \st \langle u, v \rangle = \infty \text{ exactly for } u \in \sstau^{\vee} \setminus \sstau^{\perp} \} \subset \bbO(\tau).
\end{equation*}
Under the identification $\bbO(\tau) \simeq \ssN_{\R}^{\tau}$, we have that $\ss \sigma_{\infty}^{\tau}$ is the image of $\sigma$ under the canonical projection $\ssN_{\R} \to  \ssN_{\R}^{\tau}$. 
The collection of cones $\ss \sigma_{\infty}^{\tau}$, for $\sigma \supfaceq \tau$ in $\Sigma$, form a fan $\ssSigma_{\infty}^{\tau}$ in $\bbO(\tau)$ with origin $\ssinfty_{\tau}$. 
Under the isomorphism $\ssN_{\R}^{\tau} \simeq \bbO(\tau)$, the fan $\ss \Sigma_{\infty}^{\tau}$ is identified with the star fan $\ssSigma^{\tau}$ of $\tau$ in $\Sigma$. 
For $\conezero \in \Sigma$, we have $\ss \Sigma_{\infty}^{\conezero} = \Sigma$. 

The \emph{sedentarity} of a point $x\in \TP_{\Sigma}$ denoted $\sed(x)$ is defined as the cone $\sigma \in \Sigma$ so that $x$ lies in the stratum $\bbO(\sigma)$. 
For an open subset $U$ in  $\comp\Sigma$, we define the \emph{sedentarity fan of $U$} as the subfan $\ssSigma_U$ of $\Sigma$ generated by the $\sed(x)$ for $x\in U$, and all of their faces. 
This is the smallest subfan in $\Sigma$ such that $U$ lives in $\TP_{\Sigma_{U}}$.

Given a rational polyhedral fan $\Xi$ and $y\in \Xi$, a \emph{basic open neighborhood} of $(\infty, y) \in \eR^{k} \times \Xi$ is an open set containing $(\infty, y)$ of the form 
\begin{equation}
\label{eq:prototype-basic-open}
	(t,\infty]^{k} \times V  \subset \eR^{k} \times \Xi
\end{equation} 
where $V\subset \Xi$ is a small open ball in $\Xi$ centered at $y$ which does not meet any cone $\tau$ of $\Xi$ with $y \not \in \tau$.

Suppose now that $\Sigma$ is \emph{simplicial}, i.e., each cone $\sigma$ of $\Sigma$ contains precisely $\dims{\sigma}$ rays of $\Sigma$. 
Fix a point $x\in \compSigma$, let $\delta = \sed(x)$, and let $x_{0}$ be the image of $x$ under the natural isomorphism $\Sigma_{\infty}^{\delta} \simeq \Sigma^{\delta}$.   
A \emph{basic open neighborhood} of $x$ is an open neighborhood  $\Z$-linearly isomorphic to a basic open neighborhood of $(\infty, x_{0}) \in \eR^{\dims\delta} \times \Sigma^{\delta}$ as in Equation \eqref{eq:prototype-basic-open}.

\subsection{Extended polyhedral spaces}
Let $X$ be a connected Hausdorff topological space. 
A \emph{polyhedral chart} is a pair $(W,\phi)$ where $W\subset X$ is a nonempty open subset and  $\phi\colon \ssW\to \ssU$ is a homeomorphism to a basic open neighborhood $U$ of a point in the canonical compactification of a simplicial fan. 
Suppose $(\ssW_i,\phi_i)$ and $(\ssW_j,\phi_j)$, with $\phi_i\colon \ssW_i \to \ssU_i$ and $\phi_j\colon \ssW_j \to \ssU_j$, are two polyhedral charts. 
Denote by $\ssN_i$ and $\ssN_j$ the finite-rank lattices such that $\ssN_{i,\R}$ and $\ssN_{j,\R}$ are the ambient vector spaces containing the simplicial fans used to define $\ssU_i$ and $\ssU_j$, respectively. 
These charts are \emph{compatible} if $\ssW_i\cap \ssW_j = \emptyset$ or $\ssW_i\cap \ssW_j \neq \emptyset$ and the transition map
\begin{equation*}
	\phi_i\circ \phi_j^{-1}\colon \phi_j(\ssW_i\cap \ssW_j) \to \phi_{i}(\ssW_i\cap \ssW_j) 
\end{equation*}
is induced by a global $\Z$-linear map $\psi_{ij}\colon \ssN_j\to \ssN_i$ that induces a $\Z$-linear isomorphism between the sedentarity fan of $\phi_j(\ssW_i\cap \ssW_j)$ and the sedentarity fan of $\phi_i(\ssW_i\cap \ssW_j)$. 

An \emph{extended polyhedral space} $X$ is a connected Hausdorff topological space $X$ with an atlas of pairwise compatible polyhedral charts $\{(\ssW_i, \phi_i)\}_{i\in I}$. 
A map $\phi\colon W \to U$ is \emph{compatible with the atlas} $\{(\ssW_{i}, \phi_{i})\}_{i\in I}$ if adding the pair $(W, \phi)$ to this atlas also produces an atlas.

\begin{proposition}
	If $\varphi:W\to U$ and $\varphi':W'\to U'$ are two polyhedral charts that are each compatible with the atlas $\{(\ssW_i, \phi_i)\}_{i\in I}$, then $(W,\varphi)$ and $(W',\varphi')$ are compatible with each other. 
\end{proposition}

\begin{proof}
Without loss of generality assume that $W \cap W' \neq \emptyset$. 
Let $\varphi_i\colon \ssW_i \to \ssU_i$ be a chart in the atlas such that $\ssW_{i} \cap W \cap W'\neq \emptyset$. 
Denote by $N$, $N'$, and $\ssN_{i}$ the lattices underlying  $U$, $U'$, and $\ssU_{i}$, respectively. 
By compatibility, we have $\Z$-linear maps $\psi_{i} \colon N \to \ssN_{i}$ and  $\psi_{i}' \colon \ssN_{i}\to N'$ that restrict to the transition maps $\phi_{i} \circ \varphi^{-1}$ and $\varphi' \circ \phi_{i}^{-1}$. 
These transition maps preserve sedentarity fans. 
Let $\theta_{i}\colon N\to N'$ be the composition $\psi_{i}' \circ \psi_{i}$. 
Thus, we have a collection of $\Z$-linear maps $\{\theta_i\colon N\to N' \st i\in I \text{ and } \ssW_{i} \cap W \cap W'\neq \emptyset \}$ that agree locally on intersections of the form $\phi(\ssW_i\cap \ssW_j\cap W \cap W')$. 
Denote by $\ssN_0$ and $\ssN_0'$ the lattices in the vector spaces generated by $\phi(W\cap W')$ and $\phi'(W\cap W')$, respectively. 
It follows from the local compatibility above that the maps $\theta_i$ induce an isomorphism $\theta \colon \ssN_0 \to \ssN_0'$, from which obtain the existence of a global $\Z$-linear map $N \to N'$ inducing the transition map $\phi(W\cap W') \to \phi'(W\cap W')$. 
Moreover, we claim that $\theta$ induces an isomorphism between the sedentarity fans. To see this, observe that each $\theta_i$ induces a bijection between sedentarity types, and by the compatibility of the $\theta_i$'s on the intersections, this implies that $\theta$ also induces a bijection between sedentarity types. 
From this we deduce that the image under $\theta$ of the relative interior of each cone in the sedentarity fan of $\phi(W\cap W')$ intersects the relative interior of a unique cone of the sedentarity fan of $\phi'(W\cap W')$. 
Since this also holds for $\theta^{-1}$, we get the desired claim. 
\end{proof}

\noindent Thus, we may form the \emph{completion} of $\{(\ssW_{i}, \phi_{i})\}_{i\in I}$ by adding to it all compatible polyhedral charts $(W, \phi)$. 
An extended polyhedral space is \emph{finite} if it has a finite atlas of charts. If $X$ is compact, then it is automatically finite.  

\subsection{Sedentarity types} 
\label{sec:sedentarity-types}
We say that two points $x , y \in \ssW_i$ for a chart $(\ssW_i, \phi_i)$ \emph{have the same sedentarity} if $\phi(x)$ and $\phi(y)$ have the same sedentarity in $\ssU_i$. 
The transitive closure of this relation defines a stratification of $X$ into sedentarity types. 
We denote by $\SedT(X)$ the set of sedentarity types of $X$. 
In this way, sedentarity types play the role of generic points in the sedentarity stratification of $X$ similar to the generic points of stratified algebraic varieties. 

The set $\SedT(X)$ comes with a natural partial order defined as follows. 
For two sedentarity type $s_1$ and $s_2$, we have $s_1 \prec s_2$ if one (equivalently, each) point of sedentarity type $s_2$ is the limit of points of  sedentarity type $s_1$. 
For each $s\in \SedT(X)$, we denote by $X_s$ the set of points of sedentarity type $s$. 

A point $x\in X$ is called a \emph{regular point} if $x$ has a neighborhood homeomorphic to an open subset of $\R^{d}$ for some $d$. 
The \emph{regular locus} of $X$ is the subset $\ssub X_{\reg}$ of regular points. 

\begin{proposition}\label{prop:minimum_sedenarity}
Let $X$ be a connected extended polyhedral space. Then, $\SedT(X)$ has a minimum.
\end{proposition}
\noindent By a slight abuse of notation, we denote by $\conezero$ the minimum element of $\SedT(X)$.
\begin{proof}[Proof of Proposition \ref{prop:minimum_sedenarity}] Denote by $d$ the dimension of $X$. It follows from the connectivity of $X$ that all the charts appearing in the defining atlas of $X$ are of dimension $d$. Let $s$ be a minimal element of $\SedT(X)$. The sedentarity stratum $X_{s}$ is of dimension $d$. We claim that the closure $\comp{X}_s$ of $X_s$ in $X$ is open. To see this, let $p$ be a point of $\comp{X}_s$. Since in any chart which contains $p$, there is a unique minimal sedentarity type, and $p$ is in the closure of $X_s$, this implies that the full chart is included in $\comp{X}_s$, proving the claim. 
	
By the connectivity of $X$, we infer that $X=\comp{X}_s$, which proves the proposition. 
\end{proof}

\subsection{Tropical varieties} Let $X$ be an extended polyhedral space. A \emph{tropical weight function} is a function $\omega\colon \ssub X_{\reg} \to \Z\setminus \{0\}$ that satisfies the balancing condition in the following sense. 
On each open chart $(\ssW_i, \phi_i)$, $\phi_{i}\colon \ssW_i\to \ssU_i$,  the restriction $\omega\rest{\ssW_{i,\reg}}$ is equal to the pullback of a tropical weight function on $\ssSigma_{i}$ (restricted to $\ssub U_{i,\reg}$). 
Note that the regular locus of $\compSigma_i$ does not meet the boundary at infinity. 

A \emph{tropical variety} is a pair $(X, \omega)$ consisting of an extended polyhedral space $X$ and a tropical weight function $\omega$. If $s$ is a sedentarity type of $X$, then $\comp{X}_{s}$ inherits the structure of a tropical variety from $X$.

We now discuss (rational) polyhedral structures on tropical varieties. 
An affine half-space is a subset of the form $\{ v\in \ssN_{\R} \st \langle u, v\rangle  \leq b \}$ where $u\in \ssM_{\R}$ and $b\in \R$.  
A \emph{polyhedron} in $\ssN_{\R}$ is a finite intersection of affine half-spaces.
Given a polyhedron $Q$, its tangent space $\ss \TT_{Q}$ is the linear space generated by differences of pairs of points in $Q$.  
A polyhedron $Q$ is \emph{rational} if the intersection $N \cap \ss \TT_{Q}$ is a full rank lattice in $\ss \TT_{Q}$. 
A  \emph{(rational) polytope} in $\compSigma$ is the closure of a (rational) polyhedron $Q$ in a face $\eta$ of $\compSigma$ such that the recession cone of $Q$ is a face of $\eta$.

Let $X$ be an extended polyhedral space. A (rational) polytope in $X$ is a closed subset $P\subset X$ such that, if $(W,\phi)$ is a chart with $P\subset W$, then $\phi(P)$ is a (rational) polytope. 
A \emph{(rational) polyhedral structure} on $X$ is a collection of (rational) polytopes $\cC$ of $X$ satisfying the following conditions. 
\begin{itemize}
	\item If $P\in \cC$ is contained in the chart $(W, \phi)$ and $F$ is a face of $\phi(P)$, then $\phi^{-1}(F) \in \cC$. 
	\item If $P_1, P_2 \in \cC$ and $(W_1, \phi_1)$, $(W_2, \phi_2)$ are charts containing $P_1$ and $P_2$, respectively, then $P_{12} = P_1\cap P_2$ is in $\cC$, the set $\phi_1(P_{12})$ is a face of $\phi_1(P_1)$, and  $\phi_2(P_{12})$ is a face of $\phi_2(P_2)$. 
\end{itemize}

\noindent Not all tropical varieties admit a rational polyhedral structure, however, they all admit a polyhedral structure. Furthermore, all full-dimensional polytopes are rational with our definition. 

Throughout the rest of this section, $X$ will be a tropical variety, endowed if needed with a (rational) polyhedral structure $\cC$, in which case we denote by $\subfaceq$ the partial order given by inclusion of faces.    We typically suppress $\omega$ from the notation.

\subsection{Tropical homology}
\label{sec:tropical-homology}
We briefly review the definition of tropical homology and cohomology groups of a tropical variety $X$, introduced by Itenberg--Katzarkov--Mikhalkin--Zharkov~\cite{ItenbergKatzarkovMikhalkinZharkov}; see~\cite{MikhalkinZharkov:Eigenwave, JellShawSmacka, JellRauShaw, GS-sheaf, AminiPiquerez:Hodge} for further results. See also related constructions for affine manifolds with singularities in \cite{GrossSiebert-2010, Ruddat}. 

We begin with the definition via (co)sheaves. 
First, consider the case $X = \compSigma$. 
Let $\eta\subset X$ be a face contained in some $\bbO(\tau)$. 
The \emph{tangent space} $\TT_{\eta}$ is the real vector subspace of $\ssN_{\R}^{\tau} \simeq \bbO(\tau)$ spanned by $\eta$. 
The tangent space has a natural full-rank lattice which we denote by $\ssN_{\eta}$.   

The cosheaf $\SF_p$ is defined in the following way. 
Let $x\in \overline{\Sigma}$. 
Then, for any basic open subset $U$ of $x$, define
\begin{equation*}
	\SF_p(U) = \sum_{\eta \ni x \\ \sed(\eta) = \sed(x)} \bigwedge^{p}\TT_{\eta} 
	\quad \text{ and } \quad
	\ssup\SF_p^{\Z}(U) = \sum_{\eta \ni x \\ \sed(\eta) = \sed(x) } \bigwedge^{p}\ssN_{\eta}.
\end{equation*}
If $U_i$ and  $U_j$ are basic open sets with $U_i\subset U_j$, then we have natural maps $\SF_p(U_i) \to \SF_p(U_j)$ and $\ssup\SF_p^{\Z}(U_i) \to \ssup\SF_p^{\Z}(U_j)$. 
As basic open sets form a basis for the topology of $\compSigma$, we may define $\SF_{p}(U)$ and $\ssup\SF_{p}(U)$ for an arbitrary open set $U$ to be
\begin{equation*}
	\SF_{p}(U) = \varinjlim_{V\subseteq U \\ \text{ basic}} \SF_{p}(V),
	\quad \text{ and } \quad
 \ssup\SF_{p}^{\Z}(U) = \varinjlim_{V\subseteq U \\ \text{ basic}} \ssup\SF_{p}^{\Z}(V).
\end{equation*} 
These define constructible cosheaves $\SF_p$ and $\ssup\SF_{p}^{\Z}$. 
We denote the restrictions to an open subset $U$ of $\comp{\Sigma}$ by $\SF_{p}\rest{U}$ and $\ssup\SF_{p}^{\Z}\rest{U}$, respectively. 
Also we define
\begin{equation*}
	\SF^{p}(U) = \SF_p(U)^{\dual} \hspace{20pt} \text{ and } \hspace{20pt}
	\ssub\SF^p_{\Z}(U) = (\ssup\SF_p^{\Z}(U))^{\dual}.
\end{equation*}
In a similar way, $\SF^{p}$ and $\ssub\SF_{\Z}^{p}$ are constructible sheaves on $X$.

Now consider the case of a general compact tropical variety $X$. 
Let $\phi\colon W\to U$ be a chart. 
Then, we define
\begin{equation*}
	\SF_{p}\rest{W} = \phi^*(\SF_p\rest{U}) 
	\quad \text{ and } \quad
	\ssup\SF_{p}^{\Z}\rest{W} = \phi^*(\ssup\SF_p^{\Z}\rest{U}).
\end{equation*}
These glue together to form cosheaves $\SF_p$ and $\ssup\SF_p^{\Z}$, respectively, on $X$. 
The \emph{(integral) tropical cosheaf homology} of $X$, denoted by $\rmH_{p,q}(X, \R)$ (respectively, $\rmH_{p,q}(X, \Z)$) is the $q$-cosheaf homology of $\SF_p$ (respectively, $\ssup\SF_{p}^{\Z}$)
\begin{equation*}
	\rmH_{p,q}(X, \R) = \rmH_{q}(X, \SF_{p}),
	\quad \text{ respectively, } \quad
	\rmH_{p,q}(X, \Z) = \rmH_{q}(X, \ssup\SF^{\Z}_{p}).
\end{equation*}

In a similar way, we glue the restrictions
\begin{equation*}
	\SF^{p}\rest{W} = \phi^*(\SF_p\rest{U}) 
	\hspace{20pt} \text{ and } \hspace{20pt}
	\ssub\SF_{\Z}^{p}\rest{W} = \phi^*(\ssup\SF^{\Z}_p\rest{U})
\end{equation*}
to form the sheaves $\SF^{p}$ and $\ssub\SF_{\Z}^{p}$ on $X$. 
The \emph{(integral) tropical sheaf cohomology} of $X$ is 
\begin{equation*}
	\rmH^{p,q}(X, \R) = \rmH^{q}(X, \SF^{p}),
	\quad \text{ respectively, } \quad
	\rmH^{p,q}(X, \Z) = \rmH^{q}(X, \ssub\SF_{\Z}^{p}).
\end{equation*}

When $X$ is endowed with a polyhedral structure, the above groups can be computed using cellular homology and cohomology.

Let $P$ be a $k$-dimensional polytope in $\compSigma$ for a fan $\Sigma$ in $\ssN_{\R}$, and denote by $\ss \TT_{P}$ the tangent space to $P$. When $P$ is rational, we have the full rank integral lattice $\ssN_{P} \subset \ss \TT_{P}$. 

An \emph{orientation} of $P$ is a choice of one of the two half-lines issued from the origin in the one dimensional real space $\wedge^{k}\ss \TT_{P}$. When $P$ is rational, we denote by $\ss \fv_P$ the primitive integral vector which generates this half-line. In this case, we frequently refer to this canonical multivector $\ss \fv_P$ as the orientation of $P$.  
 
If $P$ is a (rational) polytope of the extended polyhedral space $X$, then an orientation of $P$ is given by a choice of orientation of $\phi(P)$ for some chart $(W,\phi)$ with $P\subset W$. 
This induces an orientation on $\phi'(P)$ where $(W',\phi')$ is another compatible chart with $P\subset W'$.
Now suppose $X$ has a polyhedral structure $\cC$ and assign an orientation to each element $P$ of $\cC$. 
For each pair of cells $Q \ssubface P$ (that is, $Q$ is a codimension one face of $P$), we define $\sgn(P,Q)$ by $\sgn(P,Q) = 1$ if the induced orientation by $P$ on $Q$ coincides with the fixed orientation on $Q$, and $-1$ otherwise. When $P$ and $Q$ are rational and live in the same sedentarity stratum, we get the equation $\ss \fv_{P} =\sgn(P,Q) \nvect_{P/Q} \ss \wedge \ss{\mathfrak{v}}_{Q}$. 

Denote by $\cC_{q}$ the $q$-dimensional (oriented) polytopes in $\cC$. We define $\SF_{p}(P)$ to be $\SF_p(U)$ for any basic open set $U$ that meets the relative interior of $P$. Define the groups of $(p,q)$-chains and integral $(p,q)$-chains by
\begin{equation*}
	\rmC_{p,q}(X) = \bigoplus_{P\in \cC_{q}} \SF_{p}(P)
	\quad \text{and} \quad 
	\rmC_{p,q}(X,\Z) = \bigoplus_{P\in \cC_{q}} \ssup\SF^{\Z}_{p}(P).
\end{equation*}
A $(p,q)$-chain $Z$ is expressed as
\begin{equation*}
	Z = \sum_{P\in \cC_{q}} [P, v], \quad v\in \SF_{p}(P).
\end{equation*}
 By convention, we set $-[P,v] = [P,-v]$. 

Suppose $Q$ is a codimension one face of $P$. The inclusion $Q\subset P$ induces a morphism  
$\ssub \pi_{P,Q} \colon \SF_{p}(P) \to \SF_{p}(Q)$. 
Define the boundary map
\begin{equation*}
	\partial_{q}\colon \rmC_{p,q}(X) \to \rmC_{p,q-1}(X), 
	\quad 
	\partial_p([P,v]) = 
	\sum_{Q \ssubface P}  
	\sgn(P,Q)[Q,\ssub\pi_{P,Q}(v)]
\end{equation*} 
This defines a chain complex $\rmC_{p,\bul}(X)$: 
\begin{equation*}
	\cdots \rightarrow \rmC_{p,q+1}(X) \xrightarrow{\partial_{q+1}} \rmC_{p,q}(X) \xrightarrow{\partial_{q}} \rmC_{p,q-1}(X) \rightarrow \cdots
\end{equation*}
We have a similar chain complex $\rmC_{p,\bul}(X,\Z)$ for integral chains. By  \cite[Prop.~7]{MikhalkinZharkov:Eigenwave}, the tropical homology group $\rmH_{p,q}(X,\R)$ is the $q$-th homology group of $\rmC_{p,\bul}(X)$, and similarly $\rmH_{p,q}(X,\Z)$ is the $q$-th homology group of $\rmC_{p,\bul}(X,\Z)$.

\subsection{Tropical subvarieties, tropical cycles and the cycle class map}
Our presentation is close to \cite{MikhalkinRau}, see \cite{GS-sheaf} for a sheaf-theoretic description in terms of Minkowski weights.
Let $X$ be a compact tropical variety of dimension $d$. A \emph{connected tropical subvariety} of $X$ is the topological closure $Y$ in $X$ of a connected tropical variety $(\mathring{Y},\omega_{\mathring{Y}})$ 
with an embedding $\mathring{Y}\hookrightarrow X_s$ for a sedentarity type $s\in \SedT(X)$. We set $\omega_Y = \omega_{\mathring{Y}}$, and view (roughly speaking) $(Y, \omega_Y)$ as a tropical variety with possible singularities along the boundary $Y\setminus \mathring{Y}$.

A \emph{tropical cycle} of dimension $k$ is a formal sum of connected tropical subvarieties of dimension $p$ in $X$. These form a group which we denote by $\cZ_{p}(X)$. 

In this definition, we distinguish between, e.g., $(Y,\omega_{1}) + (Y,\omega_{2})$ and $(Y,\omega_{1} + \omega_{2})$. However, when we quotient $\cZ_{p}(X)$ by rational equivalence in \S~\ref{sec:rationalFunctionsDivisors}, these two cycles become equal. See Remark \ref{rmk:weights-higher-gcd}.
 
Let $(Y, \omega)$ be a connected compact tropical variety of dimension $k$, and let $\cC$ be a polyhedral structure on $Y$. 
The \emph{fundamental class} of $Y$ is the canonical element in $\rmH_{k, k}(Y, \Z)$, defined by 
\begin{equation*}
	[Y] \coloneqq \sum_{P\in \cC_{k}} \omega(P) [P, \ss \fv_{P}] \in \rmH_{k,k}(Y, \Z)
\end{equation*}
Note that this definition makes sense even if the polyhedral structure is not rational as the polytopes $P$ appearing in the above sum are full dimensional, and hence are rational. 

The \emph{cycle class map} for $X$ 
\begin{equation*}
	\cl_{X} \colon \cZ_{p}(X) \to \rmH_{p,p}(X, \Z)
\end{equation*}
is defined as follows. Given a connected tropical subvariety $(Y, \omega_Y)$ of $X$ defined as the closure of $\iota\colon \mathring Y \hookrightarrow X_s$, we consider the $(p,p)$--chain
$\sum_{Q\in \cC_p} \omega(Q) [\comp{Q}, \ss \fv_{Q}]$, where $\cC$ is a polyhedral structure on $\mathring{Y}$, and $\comp Q$ is the closure of $Q$ in $X$. The boundary of this $(p,p)$--chain is zero. The corresponding element of $\rmH_{p,p}(X,\Z)$ is $\cl_{X}(Y)$. We extend by linearity this map to the full cycle space $\cZ_{p}(X)$.

\subsection{Pushforward of cycles}
Let $X$ and $X'$ be connected compact tropical varieties. 
A continuous function $\varphi \colon X\to X'$ is a \emph{tropical morphism} if $\varphi(\ssX_{\conezero})$ is contained in a single sedentarity stratum of $X_{s}'$ of $X'$ and such that 
there are compatible atlases of charts $\{(\ssW_i,\phi_i)\}_{i\in I}$ for $X$ and $\{(\ssW_j',\phi_j')\}_{j\in J}$ for $\comp{X_{s}'}$ so that, for each $i\in I$ there is a $j\in J$ with $\varphi(\ssW_i) \subset \ssW_{j}'$ the composition $\phi_j'\circ \varphi \circ \phi_{i}^{-1}$ is induced by a  $\Z$-linear map of lattices $\ssN_{i} \to \ssN_{j}'$.

Suppose $\varphi\colon X\to X'$ is a proper tropical morphism. The pushforward $\varphi_{*} \colon \cZ_p(X)\to \cZ_p(X')$ is defined in the following way. 

Let $Y\hookrightarrow X$ be a connected $p$-dimensional tropical subvariety of $X$.

First, replace $X$ with the closure of the sedentarity stratum of $X$ containing $\mathring Y$. 
Then replace $X'$ with closure of the sedentarity stratum of $X'$ containing $\varphi(\ssX_{\conezero})$. 
This means that $\varphi(\ssX_{\conezero}) \subset \ssX_{\conezero}'$. 

Let $\cC$ be a polyhedral structure of $\mathring Y$. 
By subdividing $\cC$ if necessary, we may assume that $\varphi(P)$ is contained in a rational polytope of $\ssX_{\conezero}'$ for each full dimensional polytope $P$ in $\cC$. 
Then $\mathring Y' = \varphi_{*}(Y)$ is defined as a set by taking the union of all $\varphi(P)$ for $P\in \cC_p$ such that $\varphi(P)$ is $p$-dimensional.

Next we define a weight function on $\mathring Y'$. 
Let $\mathring Y''\subset \ssub{\mathring Y}_{\reg}'$ be the subset of all points $y'$ such that $\varphi^{-1}(y')$ is a finite subset of $\ssub{\mathring Y}_{\reg}$. 
Define $\omega_{\mathring Y'} \colon \mathring Y'' \to \Z\setminus \{0\}$ 
to be
\begin{equation*}
	\omega_{\mathring Y'}(y') = \sum_{y\in \varphi^{-1}(y')} |N_{y'} / d\varphi(N_{y})| \omega_{Y}(y) 
\end{equation*}
where $d\varphi\colon \TT_{y}Y \to \TT_{y'}X'$ is the induced map on tangent spaces $\TT_{y}Y$ and $\TT_{y'}X'$, $N_{y} \subset \TT_{y}Y$ is the lattice of integer points, and $N_{y'}$ is the lattice of integer points of $d\varphi(\TT_{y}Y)$. 
By \cite[Prop.~6.2.1]{MikhalkinRau}, the function $\omega_{\mathring Y'}$ extends to a weight function $\omega_{Y'}$ on $\ssub{\mathring Y}_{\reg}'$ so that $(\mathring Y',\omega_{\mathring Y'})$ is a tropical variety. 
Let $Y'$ be the closure of $\mathring Y'$. 
This is a connected compact tropical subvariety of $X'$, and $\varphi_{*}(Y,\omega_{Y}) = (Y',\omega_{Y'})$. 
By extending linearly, this defines $\varphi_{*}\colon \cZ_p(X)\to \cZ_p(X')$.

\subsection{The monodromy operator}
\label{sec:monodromy}
 For a compact connected tropical variety $X$, we define in this section the monodromy operator
\begin{equation*}
	\rmN \colon  \rmH_{p,q}(X, \R) \to \rmH_{p+1,q-1}(X,\R).
\end{equation*} 

Suppose $P$ is a  polytope of $\compSigma$ and $Q\subfaceq P$, with $\sigma = \sed(P)$. 
Because the recession cone of $P$ is a cone of $\Sigma^{\sigma}$, there is a unique largest face $R$ of $P$ with $\sed(R) = \sigma$ such that  $\ssub\pi_{P,Q}(R) = Q$. 
A \emph{marking} of $P$ is a choice of point $\sso{Q}\in Q$ for each face $Q$ of $P$ such that, if $\sed(Q) \neq \sigma$, then $\sso{Q} = \ssub\pi_{P,Q}(\sso{R})$ where $R$ is the face defined above. 
A \emph{marked polytope} is a pair consisting of a  polytope $P$ together with a marking. 

Let $\cC$ be a polyhedral structure on $X$ and choose points $\{\sso{P}\}_{P\in \cC}$ such that for each $P\in \cC$, the pair $(\phi(P), \{\phi(\sso{Q}) \st Q \subfaceq P\})$ is a marked polytope, where $(W,\phi)$ is a compatible chart with $P\subset W$. 
When $Q\subface P$, we write $\ssub w_{P,Q} = \ssub \pi_{P,Q}(\phi(\sso{P})) - \phi(\sso{Q})$. At the level of chains 
\begin{equation}
\label{eq:monodromy-marked-polytope}
	\rmN([P, v]) = \sum_{Q\ssubface P} \sgn(P,Q)[Q, \ssub w_{P,
	Q} \wedge \ssub \pi_{P,Q}(v)].   
\end{equation}
where $\sgn(P,Q)$ is the sign function defined in \S~\ref{sec:tropical-homology}. 
By linear extension, this defines the monodromy operator as a map 
\begin{equation*}
	\rmN \colon \rmC_{p,q}(X) \to \rmC_{p+1,q-1}(X). 
\end{equation*}
Denote by $\rmB_{p,q}(X) = \image(\partial_{p} \colon \rmC_{p-1,q+1}(X) \to \rmC_{p,q}(X))$, i.e.,  the subgroup of $(p,q)$--boundaries of $X$. 

\begin{theorem}
\label{thm:monodromy-chain-map}
The map $\rmN$ defined above is a chain map, and induces a map 
\begin{equation}
\label{eq:monodromy-mod-Bpq}
	\rmN \colon  \frac{\rmC_{p,q}(X)}{\rmB_{p,q}(X)} \to \frac{\rmC_{p+1,q-1}(X)}{\rmB_{p+1,q-1}(X)}.
\end{equation} 
Furthermore, we have the following. 
\begin{enumerate}
	\item When restricted to a map $\rmN \colon \rmH_{p,q}(X) \to \rmH_{p-1,q+1}(X)$, this map is independent of the choice of markings and agrees with the eigenwave operator defined in \cite{MikhalkinZharkov:Eigenwave}.
	\item $\rmN\circ\cl_{X} \equiv 0$. 
	\item If $q=p+1$ and $\gamma \in \rmC_{p,p+1}(X)$ verifies $\partial\gamma$ is represented by an algebraic cycle in $\rmC_{p,p}(X) / \rmB_{p,p}(X)$, then $\rmN(\gamma) \in \rmC_{p+1,p}(X) / \rmB_{p+1,p}(X)$ does not depend on the choice of markings. 
	\item If $\varphi\colon X\to X'$ is a tropical morphism of tropical varieties, then $\rmN \circ \varphi_{*} = \varphi_{*} \circ \rmN$ at the level of tropical homology as in (1). Furthermore, for $\gamma$ as in (3), the boundary of $\varphi_{*}(\gamma)$ is represented by an algebraic cycle on $X'$ and we have $\rmN \circ \varphi_{*}(\gamma) = \varphi_{*} \circ \rmN(\gamma)$. 
\end{enumerate}
\end{theorem}
\begin{proof} The first, second and half of the forth assertions are known; the new parts are (3) and the second half of (4). We briefly discuss the proofs. 

The statement $N\circ\partial = \partial\circ N$ can be verified directly using the relation 
\begin{equation*}		
	\sgn(P,Q)\sgn(Q,R) = \sgn(P,Q')\sgn(Q',R)
\end{equation*}
for a pair of faces $Q$ and $Q'$ of codimension one in $P$ that share a face $R$ of codimension one.  
	
The independence from the choice of markings can be verified directly. Then, (1) follows from the definition of eigenwave operator, see~\cite[\S~2.D.]{JellRauShaw}.

To prove (2), note that a rational polytope $P$ with canonical multivector $\ss \fv_P$ and a face of codimension one $Q$ in $P$, we have $\ssub w_{P,Q}\wedge \ssub\pi_{P,Q}(\ss \fv_P)=0$. Using the definition of  $\cl_{X}$ and $\rmN$, this gives  $\rmN\circ\cl_{X} \equiv 0$.

We now show (3). Let $P$ be a polytope of dimension $p+1$ and let $v\in \ssup \SF_{p}^{\Z}(P)$ so that $[P, v]$ is in the support of $\gamma$. Choose a different marking $\sso{P}'$. Let $\rmN'$ be the corresponding monodromy map. Let $u = \phi(\sso{P}) - \phi(\sso{P}')$ and $\ssub u_{Q} = \ssub\pi_{P,Q}(u)$.   Then,
\begin{equation*}
	\rmN(\gamma) - \rmN'(\gamma) = \sum_{Q\ssubface P} \sgn(P,Q) [Q, \ssub u_{Q} \wedge v] = \partial \left( [P, u \wedge v] \right) \in \rmB_{p+1,p}(X).
\end{equation*}
Now, let $Q$ be a rational polytope of dimension $p$ and choose a different marking $\sso{Q}'$. Denote by $\rmN'$  the corresponding monodromy map. Let  $u = \phi(\sso{Q}) - \phi(\sso{Q}')$ and set $\gamma = \sum_{P} [P, \ssub z_P]$. Note that  
\begin{equation*}
	\partial \gamma = \sum_{P} \sum_{R\ssubface P} \sgn(P,R)[R, \ssub \pi_{P,R}(\ssub z_P)].
\end{equation*}
Because $\partial \gamma$ is the fundamental class of an algebraic cycle, we can alternatively write $\partial \gamma$ 
\begin{equation*}
	\partial \gamma = \sum_{R} \ssub a_{R} [R, \ss \fv_{R}]
\end{equation*}
where the sum is over rational polytopes $R$ of dimension $p$ and the $\ssub a_{R}$ are integers, all but finitely many of them are zero. Combining the above expressions yields
\begin{equation*}
	\rmN(\gamma) - \rmN'(\gamma) = \sum_{P\ssupface Q} \sgn(P,Q) [Q, u \wedge \ssub \pi_{P,Q}(z_{P})] =  \ssub a_{Q}[Q, u \wedge \ss \fv_{Q}] = 0. 
\end{equation*}

For assertion (4), we only prove the second part. If $\partial\gamma$ is represented by an algebraic cycle $\alpha$, then $\partial \varphi_{*}(\gamma)$ is represented by $\varphi_{*}(\alpha)$.  We express $\gamma = \sum_{P} [P, \ssub z_p]$. We claim that $\varphi_{*}\rmN([P, \ssub z_{P}]) - \rmN\varphi_{*}([P, \ssub z_{P}])$ is a boundary.  The only nontrivial verification is when $\varphi(P)$ is of dimension $p$ or $p+1$; we only prove the former. Note that in this case $\varphi_{*}([P, \ssub z_{P}]) = 0$.

We choose markings of the faces $Q$ of codimension $1$ in $P$ satisfying $\dim \varphi(Q) = p$ in such a way that $\varphi_{*}(w_{P,Q})$ is constant. 
Because of (3) and the fact that $\rmN$ commutes with $\partial$, it is sufficient to treat the claim for these choices of markings. 
Then, $\varphi_{*}\rmN([P, \ssub z_{P}]) - \rmN\varphi_{*}([P, \ssub z_{P}])  = \varphi_{*}\rmN([P, \ssub z_{P}])$ and
\begin{align*}
	\varphi_{*}\rmN([P, \ssub z_{P}])
	= \sum_{Q\ssubface P \\ \dim\varphi(Q) = p-1} \sgn(P,Q) [\varphi_{*}(Q), \varphi_{*}(\ssub w_{P,Q} \wedge \ssub z_{P})].
\end{align*}
Since $\varphi_{*}(\ssub w_{P,Q} \wedge \ssub z_{P})$ is constant, the claim that this is a boundary comes from the fact in ordinary homology that the image of the boundary of a chain is a boundary. 
\end{proof}

\subsection{K\"ahler tropical varieties and the weight-monodromy property}
\label{sec:KahlerWMP}
A tropical variety $X$ satisfies the \emph{weight monodromy property} if
\begin{equation}
\label{eq:WMP}
\tag{WMP}	\rmN^{p-q}\colon \rmH_{q,p}(X,\R) \to \rmH_{p,q}(X,\R) \qquad \textrm{is an isomorphism for all $p\geq q$}.
\end{equation}
 A broad class of tropical varieties that satisfy \eqref{eq:WMP} is given by projective tropical manifolds and more generally by K\"ahler tropical varieties, as was shown in the work~\cite{AminiPiquerez:Hodge}. 
 A tropical variety $X$ is \emph{K\"ahler} if there exists an atlas of charts $(\ssW_i, \phi_i)_{i\in I}$ with $\phi_i \colon \ssW_i \to \ssU_i\subset \compSigma_i$, such that each tropical fan $\ssSigma_i$ is K\"ahler, and there is an $\omega\in \rmH^{1,1}(X,\R)$ which restricts to a K\"ahler form in each chart.
 
 By a K\"ahler tropical fan, we mean a tropical fan whose star-fans have a Chow ring which satisfies the K\"ahler package and in addition each of its open subsets verifies the tropical Poincar\'e duality, for its tropical homology groups, see \cite{AminiPiquerez:HodgeFans, AminiPiquerez:Hodge}.  
 In particular, if $X$ is a projective tropical manifold, i.e., $X$ has an atlas of charts as above such that each $\ssSigma_i$ has the same support as the Bergman fan of a matroid, then $X$ is K\"ahler (and hence satisfies the \eqref{eq:WMP}). 
 This follows from the K\"ahler package for the Chow ring of  Bergman fans proved in~\cite{AHK}, and extended to all simplicial fans with support a Bergman fanfold in \cite{ADH, AminiPiquerez:HodgeFans}, combined with Poincar\'e duality for Bergman fanfolds proved in~\cite{JellShawSmacka}. 
 More generally, by~\cite{AminiPiquerez:HodgeFans, AminiPiquerez:HomologyFans}, any projective tropical variety locally modeled by quasilinear fans is K\"ahler. 
 
\subsection{Tropical abelian varieties}
\label{sec:tropicalAbelianVarieties}
 We recall the definition of tropical abelian varieties, and refer to~\cite{KawaguchiYamaki24} and the references therein for details.
 Let $X = \R^{g}/\Z^{g}$. Let $L$ be a lattice of rank $g$ in $\R^g$. 
 This $L$ puts an integral affine structure on $X$ which defines $X$ as a tropical variety. 
 Let $M$ be the dual lattice to $L$. 
 A \emph{polarization} on $X$ is the data of a symmetric positive definite bilinear form 
\begin{equation*}
	Q \colon \R^g \times \R^g \to \R
\end{equation*}
which verifies 
\begin{equation*}
	Q(L, \Z^g) \subseteq \Z.
\end{equation*}
The bilinear form $Q$ induces an embedding 
\begin{equation*}
	\lambda \colon \Z^g \hookrightarrow M.
\end{equation*}
A tropical variety $X$ is called \emph{tropical abelian variety} if it admits a polarization $Q$.

Note that the quotient $\rquot{M}{\lambda(\Z^g)}$ is a finite abelian group. Denote by $d_1, \dots, d_g$ the corresponding invariant  factors, i.e., $d_1|\cdots |d_g$ and $\rquot{M}{\lambda(\Z^g)} \simeq \bigoplus_{j}\rquot{\Z}{d_j\Z}$. 
The vector $(d_1, \dots, d_g)$ is called the \emph{type} of $Q$, and $Q$ is called a \emph{principal polarization} if it has type $(1, \dots, 1)$. 
In this case, the map $\lambda$ is an isomorphism, and we call $X$ principally polarized.  

\section{Tropical intermediate Jacobians and Abel--Jacobi map}
\label{sec:JH_and_AJ}

In this section, we provide the constructions of the tropical intermediate Jacobians and Abel--Jacobi maps and establish several of their fundamental properties.

\subsection{Main constructions} 
\label{sec:def-AJ}

Let $X$ be a tropical variety. Given a pair of integers $q \leq p$, the $(p,q)$--th define 
\begin{equation*}
\rmL_{p,q}(X) \coloneqq \image\left(\rmN^{p-q}\colon  \rmH_{q,p}(X,\Z) \to \rmH_{p,q}(X,\R) \right)	
\end{equation*}
The \emph{intermediate Jacobian} of $X$ is
	\begin{equation*}
		\JH_{p,q}(X) \coloneqq \frac{\rmH_{p,q}(X,\R)}{\rmL_{p,q}(X)}.
	\end{equation*}
\noindent When $X$ satisfies \eqref{eq:WMP}, e.g., if $X$ is K\"ahler, the map
\begin{equation*}
	\rmN^{p-q} \colon  \rmH_{q,p}(X,\R) \to \rmH_{p,q}(X,\R)
\end{equation*}
is an isomorphism, and so $\rmL_{p,q}(X)$ 
is a full-rank lattice of $\rmH_{p,q}(X,\R)$. 
In this case, $\JH_{p,q}(X)$ is a real torus of dimension $h_{p,q}\coloneqq \dim \rmH_{p,q}(X, \R)$. In the following, we only consider the situation where $p-q = 1$. 

\begin{proposition}
	If $\varphi \colon X\to X'$ is a tropical morphism of tropical varieties, then the pushforward $\varphi_{*} \colon \rmH_{p,q}(X) \to \rmH_{p,q}(X')$ induces a morphism  $\varphi_{*} \colon \JH_{p,q}(X) \to \JH_{p,q}(X')$.  
\end{proposition}
\begin{proof}
	This follows from Theorem \ref{thm:monodromy-chain-map}(4). 
\end{proof}

\begin{remark}
\label{rem:MZIntermediateJacobian}
For tropical manifolds, in the case the sum of $p$ and $q$ equals the dimension, Mikhalkin and Zharkov define in \cite{MikhalkinZharkov:Eigenwave} a tropical intermediate Jacobian. We compare their definition to the one we provide above. 
Let $X$ be a projective tropical manifold of dimension $d$, and let  $p \geq q$ with $p+q=d$. 
The intermediate Jacobian defined in \loccit  is the torus 
\begin{equation*}
	\ssup{\JH}_{p,q}^{\mathrm{MZ}}(X) = \frac{\rmH_{q,p}(X,\R)}{\rmH_{q,p}(X,\Z)}.
\end{equation*}
Since $X$ is a projective tropical manifold, it satisfies \eqref{eq:WMP} as discussed in \S~\ref{sec:KahlerWMP}. 
The map $\rmN^{p-q}$ provides an isomorphism from $\ssup{\JH}_{p,q}^{\mathrm{MZ}}(X)$ to $\JH_{p,q}(X)$. 
	
Since we are in the case $p+q=d$, there is a pairing on $\ssup{\JH}_{p,q}^{\mathrm{MZ}}(X)$ given by the monodromy: the pairing between elements $\gamma,\gamma' \in \rmH_{q,p}(X,\R)$ is the intersection pairing between the two complementary degree cycles $\rmN^{p-q}(\gamma)$ and $\gamma'$. 
This is a symmetric and bilinear form which is nondegenerate, as conjectured by Mikhalkin and Zharkov (Conjecture 2 in their paper), and proved in \cite{AminiPiquerez:Hodge}. 
But in general this pairing is not positive by the signature calculation given in \ibid
\end{remark}
 
Denote by $\cZ_p^{\circ}(X)$ the set of \emph{homologically-trivial} tropical cycles, i.e., the kernel of the cycle class map $\cl_X \colon \cZ_p(X) \to \rmH_{p,p}(X)$. 
The \emph{tropical Abel--Jacobi map} is the homomorphism 
\begin{equation*}
	\AJ\colon  \cZ_p^{\circ}(X) \to \JH_{p+1,p}(X)
\end{equation*}
defined in the following way. 
Given a $\gamma \in \rmC_{p,p,+1}(X)$, define $\Psi_{\gamma}\colon \rmZ^{p+1,p}(X,\R) \to \R$ by
\begin{equation*}
	\Psi_{\gamma}(\omega) = \langle \rmN(\gamma), \omega \rangle
\end{equation*}
where the pairing $\bk{-}{-}$ is the natural duality pairing between $\rmC_{p+1,p}$ and $\rmC^{p+1,p}$.

\begin{proposition}
\label{prop:Psi-additative}
Given $\gamma,\gamma' \in \rmC_{p,p+1}(X, \R)$, we have
\begin{equation*}
	\Psi_{\gamma + \gamma'} = \Psi_{\gamma} + \Psi_{\gamma'}.
\end{equation*}
For any $\theta \in \rmC_{p,p}(X, \R)$, we have $\Psi_{\partial \theta} = 0$. 
In particular, $\Psi_{\gamma} = \Psi_{\gamma + \partial \theta}$.
\end{proposition}
\begin{proof}
The first additivity statement is clear. 
The second statement on $\partial \theta$ follows from the identity
\begin{equation*}
	\Psi_{\partial \theta}(\omega) = 
	\langle \rmN ( \partial \theta), \omega \rangle = 
	\langle \rmN (\theta), d\omega \rangle = 0. \qedhere 
	\end{equation*}
\end{proof}

\begin{proposition}
\label{prop:Psi-gamma-cohomology}
If $\partial \gamma$ is in the kernel of $\rmN$, e.g., if $\partial \gamma$ is the fundamental chain of a tropical $p$-cycle, then we have an induced homomorphism
\begin{equation*}
	\Psi_{\gamma}\colon \rmH^{p+1,p}(X,\R) \to \R.
\end{equation*}
\end{proposition}

\begin{proof}
	If $N\partial \gamma = 0$ and $\omega = d\beta$, then
	\begin{equation*}
		\Psi_{\gamma}(d\beta) = \langle  \rmN (\gamma), d\beta \rangle = \langle \partial \rmN (\gamma), \omega \rangle = \langle \rmN (\partial \gamma), \omega \rangle  = 0 
	\end{equation*}
	where the second-to-last equation follows from Theorem \ref{thm:monodromy-chain-map}.
\end{proof}

\noindent By Proposition \ref{prop:Psi-gamma-cohomology} and the universal coefficients theorem, we may view $\Psi_{\gamma}$ as an element of $\rmH_{p+1,p}(X, \R)$ for $\gamma$ such that $\partial \gamma$ is in the kernel of $\rmN$.  

\begin{proposition}
\label{prop:AJ-boundaries-zero}
For a homologically-trivial cycle $\alpha \in \rmZ_{p, p}(X, \Z)$ and two chains  $\gamma, \gamma' \in \rmC_{p, p+1}(X, \Z)$ such that $\partial \gamma = \partial \gamma' = \alpha$, the difference $\Psi_{\gamma} - \Psi_{\gamma'}$ lies in $\rmL_{p+1,p}(X)$. 
\end{proposition}

\begin{proof}
By Proposition \ref{prop:Psi-additative}, $\Psi_{\gamma} - \Psi_{\gamma'} = \Psi_{\gamma - \gamma'}$. 
Since $\gamma-\gamma' \in \rmZ_{p,p+1}(X, \Z)$, we have $\Psi_{\gamma - \gamma'}$ is represented by $\rmN(\gamma-\gamma')$ in $\rmH_{p+1,p}(X, \R)$, which lies in $\rmL_{p+1,p}(X)$. 
\end{proof}

\noindent Therefore, we define the \emph{tropical Abel--Jacobi map} by 
\begin{equation*}
	\AJ \colon  \cZ_p^{\circ}(X) \to \JH_{p+1,p}(X), \quad \AJ(\alpha) = \Psi_{\gamma}
\end{equation*}
where $\gamma\in \rmC_{p,p+1}(X, \Z)$ is any chain such that $\alpha = \partial \gamma$.
This tropical Abel--Jacobi map is functorial in the following sense. 
\begin{proposition}
\label{prop:AJ-functoriality-Z}
	If $\varphi \colon X \to X'$ is a tropical morphism, then
	\begin{equation*}
	\begin{tikzcd}
	\cZ_{p}^{\circ}(X) \arrow[r, "\AJ"]  \arrow[d, "\varphi_{*}"] & \JH_{p+1,p}(X) \arrow[d, "\varphi_{*}"]  \\
	\cZ_{p}^{\circ}(X') \arrow[r, "\AJ"]   & \JH_{p+1,p}(X')\mathrlap{.}
	\end{tikzcd}
	\end{equation*}
\end{proposition}

\begin{proof}
Because the pushforward $\varphi_{*}\colon \cZ_{p} \to \cZ_{p}(X')$ commutes with the cycle class map, we have an induced morphism $\varphi_{*}\colon \cZ_{p}^{\circ} \to \cZ_{p}^{\circ}(X')$.
Suppose $\alpha \in \cZ_{p}^{\circ}(X)$ and let $\alpha' = \varphi_{*}(\alpha)$. 
Let $\gamma\in \rmC_{p,p+1}(X)$ be a chain such that $\partial \gamma = \alpha$. 
So $\gamma' = \varphi_{*}(\gamma)$ satisfies $\partial (\gamma') = \alpha'$. 
By Theorem \ref{thm:monodromy-chain-map}(4), we have that $\rmN(\gamma') =\varphi_{*}(\rmN(\gamma)) + \partial \theta$ for some chain $\theta$, and so $\Psi_{\gamma'} = \varphi_{*}(\Psi_{\gamma})$ as maps $\rmZ^{p+1,p}(X') \to \R$.
\end{proof}

\subsection{Rational functions and their divisors}
\label{sec:rationalFunctionsDivisors}

Let $Y$ be a connected tropical variety of dimension $p+1$. 
A function $f\colon \ssY_{\conezero}\to \R$ is  \emph{rational} if there is a rational polyhedral structure $\cC$ on $\ssY_{\conezero}$ such that $f$ is integral affine linear on each cell of $\cC$, and moreover for each sedentarity type $\delta$ of codimension one, the slope along rays of $\cC$ corresponding to this sedentarity type are all the same; we denote this by $\mathrm{sl}_{\delta}(f)$. 
This definition is justified by the way rational functions on algebraic varieties tropicalize. 
(Note that we do not assume that the closures of the rational polytopes in $\cC$ form rational polytopes of $Y$.)

The divisor associated to a rational function was initially formulated in \cite{Mikhalkin-2006}, see also \cite{AllermannRau} and  \cite[\S~4.4]{MikhalkinRau}. 
Given a polytope $P \in \cC$ contained in $\ssY_{\conezero}$, denote by $\ssf{P}$ the linear part of the restriction of $f$ to $P$. For $Q\in \cC$ with $\dim{Q}=p$, define
\begin{equation*}
	\ss{\mathrm{ord}}_{Q}(f) = 
	-\sum_{P \ssupface Q} \omega_{Y}(P)\ssf{P}(\nvect_{P/Q}) + \ssf{Q}\left(\sum_{P \ssupface Q} \omega_{Y}(P) \nvect_{P/Q} \right) 
\end{equation*}
Let $\mathring Z_f$ be the support of $\ss{\mathrm{ord}}_{Q}(f)$, 
\begin{equation*}
	\mathring Z_{f} = \bigcup_{\substack{Q\in \cC \\ \ss{\mathrm{ord}}_{Q}(f) \neq 0}} Q,
\end{equation*}
and denote by $Z_f$ the closure of $\mathring Z_{f}$ in $Y$. The function $\mathrm{ord}_{\bul}(f)$ defines a weight function $\omega_{Z_{f}}$ on the set of regular points of $Z_{f}$. The connected components of $\mathring Z_{f}$ endowed with the restriction of $\omega_{Z_{f}}$ are connected tropical subvarieties of $Y$. We view $Z_f$ as an element of $\cZ_p(Y)$ defined as the sum of these connected tropical subvarieties. If $\ss{\mathrm{ord}}_{Q}(f)$ is everywhere $0$, then $Z_{f} = 0$. 
The \emph{divisor} of $f$ is the element of $\cZ_{p}(Y)$ given by
\begin{equation*}
	\mathrm{div}(f) = Z_f + \sum_{\delta \in \SedT{Y} \\ \dim Y_{\delta} = p} \mathrm{sl}_{\delta}(f) \comp{Y}_{\delta}. 
\end{equation*}

\noindent The following proposition follows from \cite[Thm.~4.15]{JellRauShaw}, but we include a proof for the sake of completeness. 
%\noindent The following proposition is likely well-known, but we give a proof for the sake of completeness. 
\begin{proposition}
\label{prop:div-hom-trivial}
	The fundamental class of $\mathrm{div}(f)$ is a homologically trivial cycle. 
\end{proposition}

Given a linear function $f$, the contraction by $f$ is defined on elementary $p$-forms by
\begin{equation*}
	\iota_f (v_0\wedge \cdots \wedge v_p) = 
	\sum_{i=0}^{p} (-1)^{i+1} f(v_i)\ v_0\wedge \cdots \wedge \hat{v}_i \wedge \cdots \wedge v_p.
\end{equation*}

\begin{proof}[Proof of Proposition \ref{prop:div-hom-trivial}]
The fundamental class of $\div(f)$ is
\begin{equation*}
	[\mathrm{div}(f)] = \sum_{|Q|=p} [Q, \ss{\mathrm{ord}}_{Q}(f) \ss \fv_{Q}].
\end{equation*}
Given a $(p+1)$--dimensional polytope $P \in \cC$, let $\ss \fv_{P} \in \wedge^{p+1} \ssN_{P}$ be the canonical multivector of $P$.  
Define
\begin{equation*}
	\gamma = \sum_{P} [P,\  \omega_{Y}(P) \, \iota_{\ssf{P}} \ss \fv_{P}].
\end{equation*}
We claim that $\partial \gamma = \mathrm{div}(f)$. Write
\begin{equation*}
	\partial \gamma = \sum_{Q} \sgn(P,Q)[Q, \ssub {(\partial \gamma)}_{Q}] 
	= \sum_{Q} \sgn(P,Q) \omega_{Y}(P) [Q, \iota_{\ssf{P}} \ss \fv_{P}]. 
\end{equation*}

Given a $p$-dimensional polytope $Q \in \cC$ with $Q \subset \ssY_{\conezero}$, fix generators $v_1,\ldots, v_p \in \ssN_{Q}$  so that $\ss \fv_{Q} = v_1\wedge \cdots \wedge v_{p}$ and  $\ss \fv_{P} = \sgn(P, Q) \nvect_{P/Q} \wedge \ss \fv_{Q}$ for each cell $P \ssupface Q$. 
By the Leibniz rule, we have
\begin{align*}
	\iota_{\ssf{Q}}\left[\left(\sum_{P \ssupface Q} \omega_{Y}(P)\nvect_{P/Q} \right)\wedge  \ss \fv_{Q} \right] 
	= -\ssf{Q}\left(\sum_{P \ssupface Q} \omega_{Y}(P)\nvect_{P/Q} \right) \ss \fv_{Q} 
	- \left(\sum_{P \ssupface Q} \omega_{Y}(P)\nvect_{P/Q} \right) \wedge \iota_{\ssf{Q}} \ss \fv_{Q} 
\end{align*}
The left hand side equals to zero by the balancing condition. 
Using this, we compute
\begin{align*}
	\ssub{(\partial \gamma)}_{Q} &=  \sum_{P\ssupface Q} \omega_{Y}(P) \iota_{\ssf{P}} (\nvect_{P/Q} \wedge \ss \fv_{Q}) \\
	&= -\sum_{P \ssupface Q} \left(\omega_{Y}(P)\ssf{P}(\nvect_{P/Q}) \ss \fv_{Q}\right)
	- \left(\sum_{P \ssupface Q} \omega_{Y}(P)\nvect_{P/Q} \right) \wedge \iota_{\ssf{Q}} \ss \fv_{Q} \\
	&= -\sum_{P \ssupface Q} (\omega_{Y}(P)\ssf{P}(\nvect_{P/Q}) \ss \fv_{Q})
	+ \ssf{Q}\left(\sum_{P \ssupface Q} \omega_{Y}(P)\nvect_{P/Q} \right) \ss \fv_{Q} \\
	&=\ss{\mathrm{ord}}_{Q}(f) \ss \fv_{Q} 
\end{align*}

Now suppose $\sed(Q) = \delta \neq \conezero$. 
Let $P$ be the $(p+1)$--dimensional polytope whose boundary contains $Q$. 
Then $\omega_{Y_{\delta}}(Q) = \omega_{Y}(P)$ and $\iota_{\ssf{P}}(\ss \fv_{P}) = \mathrm{sl}_{\delta}(f) \ss \fv_{Q}$ where $\ss \fv_{Q}$ is the canonical multivector of $Q \subset Y_{\delta}$, as required. 
\end{proof}
	
Let $X$ be a $d$-dimensional connected tropical variety. 
Denote by $\Prin_p(X)$ the subgroup of $\cZ_p(X)$ generated by $\varphi_{*}(\mathrm{div}(f))$ where $Y$ is a connected tropical variety of dimension $p+1$, $\varphi\colon Y\to X$ is a tropical morphism, and $f$ is a rational function on $Y$. 
By the functoriality of the cycle class map and the previous proposition, we have  $\Prin_p(X)\subseteq \cZ_p^{\circ}(X)$. 
Define the \emph{Chow groups}
\begin{equation*}
	\rmA_{p}(X) \coloneqq \cZ_p(X) / \Prin_p(X) \quad \text{and} \quad
	\ChowTriv{p}{X} \coloneqq \cZ_p^{\circ}(X) / \Prin_p(X)
\end{equation*}
of cycles and homologically-trivial cycles on $X$ modulo rational equivalence, respectively. 

\begin{remark}
	\label{rmk:weights-higher-gcd}
Consider connected tropical subvarieties $(Y,\omega_1)$ and $(Y,\omega_2)$ with the same support $Y$. 
Then we have the relation $(Y,\omega_1) + (Y,\omega_2) = (Y, \omega_1 + \omega_2)$ in $\rmA_{p}(X)$. 
Indeed, we can write the difference of these two cycles as the pushforward of the divisor of a rational function on a tropical variety of the form $C \times Y$ under the projection map to the second factor, where $C$ is a genus-0 tropical curve. 
\end{remark}

\begin{proposition}
\label{prop:pushforwardChow}
If $\varphi\colon X\to X'$ is a tropical morphism, then the push-forward $\varphi_{*}\colon \cZ_{p}(X)\to \cZ_{p}(X')$ induces morphisms
\begin{equation*}
	\varphi_{*} \colon \rmA_{p}(X) \to \rmA_{p}(X') \quad \text{and}\quad \varphi_{*} \colon \ChowTriv{p}{X}  \to \ChowTriv{p}{X'}.
\end{equation*}
\end{proposition}

\begin{proof}
Clearly we have $\varphi_{*}(\Prin_{p}(X)) \subseteq \Prin_{p}(X')$, which yields the first map. 
The second one is a consequence of the functoriality of the cycle class map. 
\end{proof}

\begin{theorem}
\label{thm:AbelJacobiPrin}
	The Abel--Jacobi map vanishes on $\Prin_p(X)$. 
\end{theorem}
\noindent As a result, we view the Abel--Jacobi map as a morphism 
\begin{equation*}
	\AJ\colon \ChowTriv{p}{X} \to \JH_{p+1, p}(X). 
\end{equation*}

\begin{proof}[Proof of Theorem \ref{thm:AbelJacobiPrin}]
Let $Y$ be a $(p+1)$--dimensional tropical variety and $\varphi\colon Y\to X$ be a tropical morphism. 
We show that  $\AJ(\div(f)) = 0$ in $\JH_{p+1,p}(Y)$ for any rational function $f$ on $Y$. 
From this we deduce that $\AJ(\varphi_{*}(\div(f))) = 0$ in $\JH_{p+1,p}(X)$ by functoriality of the Abel--Jacobi map as in Proposition \ref{prop:AJ-functoriality-Z}, and the theorem follows.

Using the notation in the paragraph preceding Proposition \ref{prop:div-hom-trivial}, we have $\AJ(\mathrm{div}(f)) = \Psi_{\gamma}$ where $\gamma$ is the $(p,p+1)$--chain from the proof of that proposition: 
\begin{equation*}
	\gamma = \sum_{P} [P,\  \omega_{Y}(P) \, \iota_{\ssf{P}} \ss \fv_{P}].
\end{equation*}
Suppose $Q \subset \ssY_{\conezero}$. We write 
\begin{equation*}
	\rmN(\gamma) = \sum_{Q} [Q, \ssub u_{Q}].
\end{equation*}
By the Leibniz rule for $\iota_{\ssf{P}}$ and the fact that $\ssub w_{P,Q}$ is parallel to $\ss \fv_{P}$, we have 
\begin{equation*}
	\ssub u_{Q} = \sum_{P\ssupface Q} \omega_{Y}(P) \ssub w_{P,Q} \wedge \ssub \pi_{P,Q}(\iota_{\ssf{P}} \ss \fv_{P}) =
	\begin{cases}
	-\sum_{P\ssupface Q} \omega_{Y}(P) \ssf{P}(\sso{P} - \sso{Q}) \ss \fv_{P} & \text{if } \sed(Q) = \conezero, \\
	0 & \text{otherwise.}
	\end{cases}
\end{equation*}	
(See \S~\ref{sec:tropical-homology} and \S~\ref{sec:monodromy} for the definition of $\ssub \pi_{P,Q}$ and $\ssub w_{P,Q}$, respectively. Furthermore, we simplify our expressions by omitting the mention of $\phi$ appearing in the atlas of charts.)
We show that $\rmN(\gamma) = \partial \xi$ where  
\begin{equation*}
	\xi = -\sum_{P}\omega_{Y}(P)[P, \ssf{P}(\sso{P})\ss \fv_{P}].	
\end{equation*}
Since $\ss \fv_{P}$ is a $(p+1)$--multivector, we have that $(\partial \xi)_{Q} = 0$ when $\sed(Q) \neq \conezero$. 
For $Q$ whose sedentarity is $\conezero$, we have
\begin{align*}
	\ssub{(\partial\xi)}_{Q} &= -\sum_{P\ssupface Q} \omega_{Y}(P) \ssf{P}(\sso{P})\ss \fv_{P} = -\sum_{P\ssupface Q} \omega_{Y}(P) \ssf{P}(\sso{P} - \sso{Q} +\sso{Q})\ss \fv_{P} \\
	&= -\sum_{P\ssupface Q} \omega_{Y}(P) \ssf{P}(\sso{P} - \sso{Q})\ss \fv_{P} - \sum_{P\ssupface Q} \omega_{Y}(P) \ssf{P}(\sso{Q})\ss \fv_{P}.
\end{align*}
Using $\ssf{P}(\sso{Q}) = \ssf{Q}(\sso{Q})$ and the balancing condition on $Y$, we have
\begin{equation*}
	\sum_{P\ssupface Q} \omega_{Y}(P) \ssf{P}(\sso{Q})\ss \fv_{P} = \ssf{Q}(\sso{Q})\sum_{P\ssupface Q} \omega_{Y}(P) \ss \fv_{P} = \ssf{Q}(\sso{Q})\ssub{(\partial[Y])}_{Q} = 0. 
\end{equation*}
Therefore, $\ssub{(\partial\xi)}_{Q} = \ssub u_{Q}$ for all $Q$, and so $\Psi_{\gamma}\equiv 0$, as required. 
\end{proof}

\begin{proof}[Proof of Theorem \ref{thm:IntroAJ}]
	This theorem follows from Proposition \ref{prop:pushforwardChow} and Theorem \ref{thm:AbelJacobiPrin}.
\end{proof}

\subsection{An obstruction to algebraic equivalence}

Let $X$ be a tropical variety. 
Recall from the introduction that we set
\begin{equation*}
	\rmK_{p,q}(X) = \image(\rmN^{p-q-1}\colon \rmH_{q+1,p-1}(X,\Z) \to \rmH_{p,q}(X,\R))
\end{equation*}
and
\begin{equation*}
	\rmQ_{p,q}(X) = \frac{\rmH_{p,q}(X,\R)}{\rmK_{p,q}(X)}.
\end{equation*}
The monodromy map $\rmN \colon \rmH_{p,q}(X,\R)\to \rmH_{p+1,q-1}(X,\R)$ induces a map
\begin{equation*}
	\rmN \colon \JH_{p,q}(X) \to \rmQ_{p+1,q-1}(X).
\end{equation*}
We define an obstruction in $\rmQ_{p+2,p-1}(X)$ to triviality under algebraic equivalence of tropical cycles. 
Suppose $Z_1$ and $Z_2$ are connected tropical subvarieties of $X$. 
We write $Z_1 \sim Z_2$ if there is a tropical curve $B$, two points $b_1,b_2\in B$, and a subvariety $W\subset X\times B$ such that
\begin{equation*}
	(\pi_1)_{*}[W \cap (X\times b_1) - W \cap (X\times b_2)] = Z_1-Z_2
\end{equation*}
where $\pi_1\colon X\times B \to X$ is the first projection and $\cap$ is stable intersection, see \cite[\S~3.6]{MaclaganSturmfels}. 
This induces a reflexive and symmetric on $\cZ_{p}(X)$. The transitive closure of this relation is called \emph{algebraic equivalence}, and cycles in the same equivalence class are called \emph{algebraically equivalent}. 
By \cite[Lem.~6]{Zharkov}, algebraically equivalent cycles are homologically equivalent. 

\begin{theorem}
\label{thm:obstructionAlgebraicEquivalence}
	If a $p$-cycle $\alpha$ is algebraically equivalent to 0, then 
	\begin{equation*}
		\rmN(\AJ(Z)) = 0 \quad \text{in} \quad \rmQ_{p+2,p-1}(X).
	\end{equation*}
\end{theorem}

\begin{proof} Let $Z$ be an element of $\cZ_{p}^{\circ}(X)$ representing $\alpha$ in the Chow ring. 
Because the basic equivalences $Z_1\sim Z_2$ generates algebraic equivalence, we may assume, without loss of generality, that $Z = Z_1 - Z_2$ for tropical subvarieties $Z_1, Z_2$, that $B$ is a tropical curve, and $W\subset X\times B$ is a $(p+1)$--dimensional tropical variety that contains $Z_1$ and $Z_2$ as fibers. 
By functorality of $\rmN$ and $\AJ$ as in Theorem \ref{thm:monodromy-chain-map} and Proposition \ref{prop:AJ-functoriality-Z}, we have the following commutative diagram:
\begin{equation*}
	\begin{tikzcd}
	\cZ_{p}^{\circ}(W) \arrow[r, "\AJ"]  \arrow[d] & \JH_{p+1,p}(W) \arrow[r, "\rmN"] \arrow[d] & \rmQ_{p+2,p-1}(W) \arrow[d] \\
	\cZ_{p}^{\circ}(X \times B) \arrow[r, "\AJ"]  \arrow[d] & \JH_{p+1,p}(X \times B) \arrow[r, "\rmN"] \arrow[d] & \rmQ_{p+2,p-1}(X\times B) \arrow[d] \\
	\cZ_{p}^{\circ}(X) \arrow[r, "\AJ"]   & \JH_{p+1,p}(X) \arrow[r, "\rmN"]  & \rmQ_{p+2,p-1}(X)\mathrlap{.}
	\end{tikzcd}
\end{equation*}
Viewing $Z$ as a cycle in either $W$ or $X$, we see that $Z\in \cZ_{p}^{\circ}(X)$ is the image of $Z \in \cZ_{p}^{\circ}(X)$ under the composition of the vertical arrows on the left of the above diagram. 
Since $\dim(W) = p+1$, we have $\rmQ_{p+2,p-1}(W) = 0$, and hence $\rmN(\AJ(Z)) = 0$.  
\end{proof}

\subsection{Tropical varieties with torsion Abel--Jacobi image}
\label{sec:torsion-AJ}

Motivated by questions on torsion of the Ceresa class (see \S~\ref{sec:torsion-Ceresa}), we derive a general criterion that ensures the image of a homologically trivial cycle under the Abel--Jacobi map is torsion.  

\begin{theorem}
\label{thm:torsion-AJ}
	 Suppose $X$ is a tropical variety that satisfies \eqref{eq:WMP}. 
	 Assume that the monodromy $\rmN$ is rational. 
	 Then for any homologically trivial cycle in $X$, its Abel--Jacobi image is torsion. 
\end{theorem}

\begin{proof}
To prove the claim, note that the monodromy operator for $\Gamma$ has rational coefficients. 
Suppose $\alpha$ is a homologically trivial $p$-cycle of $X$. 
By assumption, the monodromy map  $\rmN \colon \rmH_{p,p+1}(X,\Z) \to \rmH_{p+1,p}(X,\R)$ takes values in $\rmH_{p+1,p}(X,\Q)$. 
This implies the inclusion $\rmL_{p+1,p}(X) \subseteq \rmH_{p+1,p}(X, \Q)$. 
Since $X$ satisfies \eqref{eq:WMP}, these two lattices have the same rank, so we can find a positive integer $n$ such that $n \rmH_{p+1,p}(X, \Z) \subseteq \rmL_{p+1,p}(X)$. 
The $(p, p+1)$--chain $\gamma$ which bounds $\alpha$ has integral coefficients, therefore $\Psi_{\gamma}$ (as defined in \S~\ref{sec:def-AJ}) belongs to $\rmH_{p+1,p}(X,\Q)$. 
There exists a positive integer $m$ such that $m\Psi_{\gamma}$ belongs to $\rmH_{p+1,p}(X,\Z)$. 
We conclude $mn\Psi_{\gamma} \in \rmL_{p+1,p}(X)$. This proves that $\AJ(\alpha)$ is torsion.
\end{proof}

\subsection{Tropical Albanese}
\label{sec:Tropical-Albanese}

Let $X$ be a $d$-dimensional K\"ahler tropical variety.  
The tropical \emph{Albanese} variety of $X$ is 
\begin{equation*}
	\Alb(X) = \JH_{1,0}(X,\mathbb{R}).
\end{equation*}
The Albanese $\Alb(X)$ is a compact torus with the integral affine structure given by the lattice $\rmH_{1,0}(X, \Z)$. 

Suppose that the K\"ahler class $\omega$ on $X$ belongs to $\rmH^{1,1}(X, \Q)$ and the monodromy takes values in rational tropical (co)homology. 
In this case, $\Alb(X)$ is endowed with a polarization which makes it a tropical abelian variety, see \S~\ref{sec:tropicalAbelianVarieties}. 
This is given as follows. 

By the weight-monodromy property \eqref{eq:WMP}, and assumption on rationality of $\rmN$, the map 
\begin{equation*}
\rmN \colon \rmH^{1,0}(X, \Q) \to \rmH^{0,1}(X, \Q)
\end{equation*}
is an isomorphism (see \cite{AminiPiquerez:Hodge} for a cohomological description of the monodromy map $\rmN$). By Hodge--Riemann property~ \loccit, the pairing 
\begin{align*}
\rmH^{1,0}(X, \Q) \times \rmH^{1,0}(X, \Q) &\to \Q\\
\alpha, \beta &\mapsto -\deg(\alpha \rmN(\beta) \omega^{d-1})	
\end{align*}
is positive definite. 
Using the duality between $\rmH_{1,0}(X, \Q)$ and $\rmH^{1,0}(X,\Q)$, this induces a polarization on $X$. 
By fixing a base-point $\flat \in X$, the Abel--Jacobi map specializes to 
\begin{equation*}
	\AJ\colon  X \to \ChowTriv{0}{X} \to \Alb(X)
\end{equation*}
where the first map is $x\mapsto [x]-[\flat]$.  The tropical Albanese variety satisfies the following universal property. 
\begin{proposition}
\label{prop:albaneseUniversalProperty}
	If $A$ is a tropical abelian variety and $X \to A$ is a morphism, then there is a unique morphism $\Alb(X) \to A$ so that the following diagram commutes
	\begin{equation*}
		\begin{tikzcd}
			X \arrow[rd] \arrow[r, "\AJ"] & \Alb(X)  \arrow[d] \\
			   & A \mathrlap{.}
		\end{tikzcd}
	\end{equation*}
\end{proposition}

\begin{proof}
	By functoriality of the Abel--Jacobi map in Theorem \ref{thm:IntroAJ}, we get a map $\Alb(X) \to \Alb(A)$. By Proposition \ref{prop:monodromyAbelianVariety}, we have that $\Alb(A) \cong A$.  
\end{proof}

\section{Intermediate Jacobians of  tropical compact tori}
\label{sec:monodromyAbelianVarieties}

In this section, let $X = \R^{g}/\Z^{g}$ with an integral affine structure given by a rank $g$ lattice $L\subset \R^{g}$. 
Since $X$ is a real compact torus, its tangent bundle is trivial, so the cosheaves $\SF_{p}$ and $\ssup\SF_{p}^{\Z}$ are isomorphic to the constant cosheaves $\wedge^{p}\R^{g}$ and $\wedge^{p}L$, respectively. 
Moreover, $\rmH_1(X,\Z) $ is canonically isomorphic to $\Z^{g}$.  
We therefore have isomorphisms
\begin{equation}
\label{eq:tropicalHomologyTropAbelianVariety}
\begin{array}{l}
	\rmH_{p,q}(X,\Z) = \rmH_{q}(X,\ssup\SF_{p}^{\Z}) \cong  \wedge^{q} \Z^g \otimes_{\Z} \wedge^{p} L, \\[2mm]
\rmH_{p,q}(X,\R) = \rmH_{q}(X,\SF_{p}) \cong  \wedge^{q} \R^g \otimes \wedge^{p} \R^{g}.
\end{array}
	\end{equation} 
Define the linear map
\begin{align*}
	&\varphi \colon  \wedge^{q} \Z^g \otimes_{\Z} \wedge^{p} L \to \  \wedge^{q-1} \Z^g \otimes_{\Z} \Z^{g} \wedge (\wedge^{p} L) \  \subset \wedge^{q-1} \R^{g} \otimes \wedge^{p+1} \R^{g},  \\
	& u_1 \wedge \cdots \wedge u_{q} \otimes v \mapsto \sum_{k=1}^{q} (-1)^{k} u_{1} \wedge \cdots \wedge \widehat{u_{k}} \wedge \cdots \wedge u_q \otimes u_k \wedge v.
\end{align*}

\begin{proposition}
\label{prop:monodromyAbelianVariety}
The following diagram commutes
\begin{equation*}
	\begin{tikzcd}
	\rmH_{p,q}(X,\Z) \arrow[r, "\rmN"]  \arrow[d] & \rmH_{p+1,q-1}(X,\R)\arrow[d]  \\
	\wedge^{q} \Z^g \otimes_{\Z} \wedge^{p} L \arrow[r, "\varphi"]  & \wedge^{q-1} \R^g \otimes \wedge^{p+1} \R^{g}\mathrlap{.}
	\end{tikzcd}
\end{equation*}
In this case, the tropical intermediate Jacobian of the tropical compact torus $X$ is given by
\begin{equation*}
	\JH_{p,q}(X) \cong \frac{\wedge^{q}\R^{g} \otimes \wedge^{p}\R^{g}}{\image(\varphi^{q-p}\colon \wedge^{p}\Z^{g} \otimes \wedge^{q}L \to \wedge^{q}\R^{g} \otimes \wedge^{p}\R^{g})}\,.
\end{equation*}
\end{proposition}

\begin{figure}
	\includegraphics[height=2.5cm]{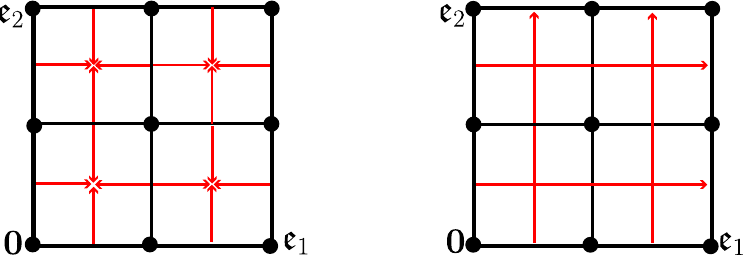}
	\caption{An illustration of the proof of Proposition \ref{prop:monodromyAbelianVariety}}
	\label{fig:monodromyExample}
\end{figure}

\begin{proof}
Denote by $\e_1,\ldots, \e_g$ the standard basis of $\Z^{g}$. A fundamental domain in the universal cover $\pi \colon \R^{g}\to X$ of the tropical compact torus $X$ is $[0,1]^{g}$.
We give $X$ a cubical complex structure with $2^{g}$ maximal cells by taking the barycentric subdivision of each copy of $[0,1]$ and forming the induced subdivision on the product $[0,1]^{g}$, see Figure \ref{fig:monodromyExample}. 
Let $I_0 = [0,\frac{1}{2}]$ and $I_{1} = [\frac{1}{2},1]$. 
The cells of this cubical complex are in bijection with the set of functions 
\begin{equation*}
	 f\colon \{1,\ldots,g\} \to \{0, \tfrac{1}{2}, 1, I_0, I_1 \} 
\end{equation*}
where the cell corresponding to the function $f$ is
\begin{equation*}
	P_f = \pi(\tilde{P}_{f}) \subset X, \quad \text{ where } \tilde{P}_f = \prod_{k=1}^{g} f(k) \subset [0,1]^{g}.
\end{equation*}
Because $\pi\rest{\tilde{P}_{f}}$ is an isomorphism to its image, we identify $\tilde{P}_{f}$ with $P_{f}$. 
Denote by $\sso{f}$ the barycenter of $P_{f}$. 
We write $h\ssubface f$ if  $P_{h}$ is a codimension one face of $P_{f}$. 
When $h\ssubface f$, there is a unique $d\in \{1,\ldots,g\}$ such that $f(d)$ is an interval and $h(d)$ is a vertex of $f(d)$. 
So $\op{f} - \op{h} $ equals $\frac{1}{4} \e_{d}$ or $-\frac{1}{4} \e_{d}\, $  if $h(d)$ is larger or smaller than the barycenter of $f(d)$, respectively. 
We denote this $d$ by $d_{f,h}$. 
The monodromy operator on chains is given by
\begin{equation*}
	\rmN([P_{f}, v]) = \sum_{h\ssubface f} \sgn(f,h)[P_{h}, (\op{f} - \op{h}) \wedge v]
\end{equation*}
where $\sgn(f,h) = \sgn(P_{f}, P_{h})$ is the sign function in \S~\ref{sec:tropical-homology}. 
We consider the case $q=g$, the proof is similar for $q<g$.  
Let $\cube_{g}$ denote the set of functions corresponding to the $g$-dimensional cells of $X$; these are exactly the functions $f\colon \{1,\ldots,g\} \to \{I_0,I_1\}$. 
So 
\begin{equation*}
	\rmH_{p,g}(X, \Z) = \left\{[\sum_{f \in \cube_{g}} P_{f}, v  ] \st  v \in \wedge ^{p}L \right\}.
\end{equation*}
Under the identification in Equation \eqref{eq:tropicalHomologyTropAbelianVariety}, the element of $\rmH_{p,g}(X,\Z)$ given by
\begin{equation*}
	[\sum_{f \in \cube_{g}} P_{f}, v] \quad \text{ corresponds to } \quad  \e_1 \wedge \cdots \wedge \e_g \otimes_{\Z} v.
\end{equation*}
Then, we compute
\begin{align*}
	\rmN([\sum_{f \in \cube_{g}} P_{f}, v]) 
	&= \sum_{f \in \cube_{g}}\sum_{h\ssubface f} \sgn(f,h) [P_{h}, (\op{f} - \op{h}) \wedge v] \\
	&= \sum_{k=1}^{g} \left( \sum_{f\in \cube_g}  \sum_{h:d_{f,h} = k}  \sgn(f,h) [P_{h}, (\op{f} - \op{h}) \wedge v]   \right).
\end{align*}
We claim that
\begin{equation*}
	\sum_{f\in \cube_g} \sum_{h:d_{f,h} = k}   \sgn(f,h) [P_{h}, (\op{f} - \op{h}) \wedge v]  \mapsto (-1)^{k} \e_1 \wedge \cdots \wedge \widehat{\e_k} \wedge \cdots \wedge \e_g \otimes \e_i \wedge v
\end{equation*}
via the isomorphism in Equation \eqref{eq:tropicalHomologyTropAbelianVariety}. We prove this for $k=1$. Then, the sum on the left is
\begin{equation}
\label{eq:dfh_eq_1}
	 \sum_{f\in \cube_g \\ f(1) = I_0} \sum_{h:d_{f,h} = 1}  \sgn(f,h) [P_{h}, (\op{f} - \op{h}) \wedge v] + \sum_{f\in \cube_g \\ f(1) = I_1} \sum_{h:d_{f,h} = 1}  \sgn(f,h) [P_{h}, (\op{f} - \op{h}) \wedge v].
\end{equation}
Consider the sum on the left, i.e., assume that $f(1)=I_0$. 
Then $h(1)$ is either $0$ or  $\frac{1}{2}$. 
In the former case, we have $\sgn(f,h) = -1$ and $\op{f} - \op{h} = \frac{1}{4}\e_1$ and in the latter case we have $\sgn(f,h) = 1$ and $\op{f} - \op{h} = -\frac{1}{4}\e_1$.
So the sum on the left equals
\begin{equation*}
	\sum_{h\in\cube_{g-1} \\ h(1)=0}  -[P_{h}, \tfrac{1}{4}\e_1 \wedge v] +\sum_{h\in\cube_{g-1} \\ h(1)=1/2}  [P_{h}, -\tfrac{1}{4}\e_1 \wedge v] 
	= -\tfrac{1}{4}\sum_{h\in\cube_{g-1} \\ h(1)=0}  [P_{h}, \e_1 \wedge v]  -\tfrac{1}{4}\sum_{h\in\cube_{g-1} \\ h(1)=1/2}  [P_{h}, \e_1 \wedge v]. 
\end{equation*}
The two sums on the right are $(p+1,g-1)$--cycles that are homologous to each other.  
Under the identification in Equation \eqref{eq:tropicalHomologyTropAbelianVariety}, these cycles correspond to
\begin{equation*}
	-\tfrac{1}{4}\e_2\wedge \e_3 \wedge \cdots \wedge \e_g \otimes_{\Z} \e_1 \wedge v \in \wedge^{g-1} \Z^{g} \otimes_{\Z} \Z^{g} \wedge (\wedge^{p} L).
\end{equation*}  
In a similar way, the second double summation in \eqref{eq:dfh_eq_1} also equals to two $(p+1,g-1)$--cycles, each corresponding to the above element in $\wedge^{g-1} \Z^{g} \otimes_{\Z} \Z^{g} \wedge (\wedge^{p} L)$. 
This proves the claim for $k=1$, and the proof of the claim for $k\geq 2$ is similar. 
In total, we get that $\rmN([\sum_{f \in \cube_{g}} P_{f}, v ])$ corresponds to
\begin{equation*}
	 \sum_{k=1}^{g} (-1)^{k}\e_1\wedge \cdots \wedge \widehat{\e_k} \wedge \cdots \wedge \e_g \otimes_{\Z} \e_k \wedge v = \varphi(\e_1\wedge \cdots \wedge \e_g \otimes_{\Z} v),
\end{equation*} 
as required. 
\end{proof}

\section{Curves and their Jacobians}
\label{sec:curvesJacobians}
Let $\Gamma$ be a connected genus-$g$ tropical curve. 
Fix a model $(G,\ell)$ with $G = (V, E)$ and an orientation on the edges of $G$. 
We suppose that $G$ has no separating edges. 
Denote by $\Div_0(\Gamma)$ the group of degree-0 divisors on $\Gamma$. 
Note that $\Div_0(\Gamma)$ is exactly the group $\cZ_0^{\circ}(\Gamma)$ of homologically trivial $0$-cycles on $\Gamma$.  
\begin{theorem}
\label{thm:tropicalJHCurves}
The tropical intermediate Jacobian $\JH_{1,0}(\Gamma)$ of $\Gamma$ is canonically identified with the Jacobian $\Jac(\Gamma)$ of $\Gamma$: 
\begin{equation*}
	\Jac(\Gamma) \cong \JH_{1,0}(\Gamma).
\end{equation*}
Furthermore, we have the following commutative diagram:
\begin{equation*}
	\begin{tikzcd}
	\Div_0(\Gamma) \arrow[r, "\AJ"]  \arrow[d, "="] & \Jac(\Gamma) \arrow[d, "\cong"]  \\
	\Div_0(\Gamma) \arrow[r, "\AJ"]  & \JH_{1,0}(\Gamma)\mathrlap{.}
		\end{tikzcd}
\end{equation*}
\end{theorem}
\noindent In the above commutative diagram, the top $\AJ$ is the tropical Abel--Jacobi map as defined in \cite{MikhalkinZharkov:TropicalCurves, BakerFaber, BakerNorine}, whereas the bottom $\AJ$ is the Abel--Jacobi map defined in \S~\ref{sec:JH_and_AJ}.

We begin by recalling the definition of the Jacobian $\Jac(\Gamma)$. 
With $\Gamma$ understood, we write $\ssH_{\Z} = \rmH_1(\Gamma,\Z)$ and $\ssH = \ssH_{\Z}\otimes \R$. 
Denote by $\ssub\Omega_{\Z} = \ssub\Omega_{\Z}(\Gamma)$ the space of integral harmonic 1-forms on $\Gamma$, and $\Omega = \ssub\Omega_{\Z} \otimes \R$. Given any path $\gamma$ in $\Gamma$, define
\begin{equation}
\label{eq:evgamma}
	\ev{\gamma} \in \Omega^{\vee} \quad \omega \mapsto \ev{\gamma}(\omega) \coloneqq \int_{\gamma} \omega. 
\end{equation}
Under the embedding
\begin{equation*}
	\ev{}\colon  \ssH_{\Z} \hookrightarrow \Omega^{\vee} \quad \gamma \mapsto \ev{\gamma},
\end{equation*}
$\ssH_{\Z}$ is realized as a full-rank lattice inside $\Omega^{\vee}$. 
	The \emph{Jacobian} of $\Gamma$ is the compact real torus 
\begin{equation*}
	\Jac(\Gamma) \coloneqq \frac{\Omega^{\vee}}{\ssH_{\Z}}.
\end{equation*}
The tropical curve $\Gamma$ has an intrinsic tropical structure making it into a K\"ahler tropical variety of dimension 1. Each point $\flat \in \Gamma$ of valence $d$ has a neighborhood identified with a neighborhood of $0$ in the 1-skeleton of the fan of the projective space $\P^{d-1}$. The Abel--Jacobi map gives an embedding of $\Gamma$ into $\Jac(\Gamma)$ as a tropical subvariety (we use that $G$ has no separating edges). 
 
The tropical structure of $\Gamma$ may be described in the following way. Fix an orientation for the edges of $\Gamma$. 
The tangent space at each point of $\Jac(\Gamma)$ is naturally identified with $\Omega^{\vee}$, and its integral affine structure is given by $\ssub\Omega_{\Z}\subset \Omega$. 
Given an edge $e$, denote by $\src{e}$ its \emph{source} vertex and $\dst{e}$ its \emph{target} vertex. 
We have:
\begin{equation*}
	\ssup\SF_{0}^{\Z}(e) = \Z, \quad \ssup\SF_{1}^{\Z}(e) = \Z \cdot 1_{e} = \Z \cdot \frac{\ev{e}}{\ell(e)} \subset \ssub\Omega_{\Z}^{\vee}
\end{equation*}
where $1_{e}$ is the unit tangent vector to $e$ with respect to the orientation of $e$. 
Given a vertex $v\in V$, we have
\begin{equation*}
	\ssup\SF_{0}^{\Z}(v) = \Z, \quad \ssup\SF_{1}^{\Z}(v) = \rquot{(\bigoplus_{e\ssupface v} \ssup\SF_{1}^{\Z}(e)) }{\langle \sum_{e\ssupface v} \sgn(e, v)\ 1_e \rangle} \subset \ssub\Omega_{\Z}^{\vee}.
\end{equation*}
Here, $\sgn(e, v) = 1$ if $v=\src{e}$ and  $\sgn(e, v) = -1$ if $v=\dst{e}$. 
Using this description, we have:
	\begin{equation}
	\label{eq:tropicalHomologyCurves}	
		\rmH_{0,1}(\Gamma,\Z) = \ssH_{\Z},\quad \rmH_{1,0}(\Gamma,\Z) = \ssub\Omega_{\Z}^{\vee}, \quad \rmH_{0,1}(\Gamma,\R) = H, \quad \rmH_{1,0}(\Gamma,\R) = \Omega^{\vee}.
	\end{equation}
The monodromy operator is defined at level of chains $\rmN\colon \rmC_{0,1}(\Gamma, \Z) \to \rmC_{1,0}(\Gamma, \R)$ in the following way. 
Given an edge $e$, let $o_e$ be the midpoint of $e$, and $o_v = v$ for each vertex $v$. 
Then,
\begin{equation}
\label{eq:monodromyAndIntegration}
	\rmN([e, 1]) 
	= \sum_{v\ssubface e} \sgn(e,v)\, [v, (o_e - o_v)] 
	= \tfrac{1}{2}([\src{e}, \ev{e}] + [\dst{e}, \ev{e}]) 
	\equiv \ev{e} \mod \rmB_{1,0}(\Gamma, \R).
\end{equation}

\begin{proposition}
\label{prop:monodromyAndIntegrtion}
Under the identifications in \eqref{eq:tropicalHomologyCurves}, we have the identification $\rmN = \ev{}$.
\end{proposition}
\begin{proof}
Under the identification $H\cong \rmH_{0,1}$, the cycle $\sum n_e e$ in $H$ corresponds to $\sum n_e [e,1_e]$ in $\rmH_{0,1}$; we denote both by $\gamma$. 
We have	
\begin{equation*}
	\rmN(\gamma) 
	= \sum_{e} n_e \rmN([e,1_e]) 
	= \sum_{e} n_e \ev{e} 
	= \ev{\sum_e n_e e} 
	= \ev{\gamma}. \qedhere
\end{equation*}
\end{proof}

\begin{proof}[Proof of Theorem \ref{thm:tropicalJHCurves}] 
The canonical identification $\Jac(\Gamma) \cong \JH_{1,0}(\Gamma)$ is a direct consequence of Proposition \ref{prop:monodromyAndIntegrtion}. 
Consider the $0$-cycle $\rflat-\flat$ for $\flat, \rflat\in \Gamma$. 
Let $\eta$ be an oriented path in $\Gamma$ with $\partial \eta = \rflat-\flat$. 
The top Abel--Jacobi map sends $\rflat-\flat$ to the linear form  $\langle \ev{\eta}, -\rangle$ on $\Omega$. 
The bottom Abel--Jacobi map sends $\rflat-\flat$ to the linear form $\langle \rmN(\eta), - \rangle$ on $\rmH^{1,0}(\Gamma)$, which can be identified with $\Omega$ by duality and the identification $\rmH_{1,0} \cong \Omega^{\vee}$.  
These two forms agree under the identification in Equation \eqref{eq:monodromyAndIntegration}. 
\end{proof}

Fixing a basis of $\ssH_{\Z} = \rmH_{1}(\Gamma,\Z)$ yields isomorphisms $\Omega^{\vee} \cong \R^{g}$ and $\ssH_{\Z}\cong \Z^{g}$. Define 
\begin{equation}
\label{eq:polarizationOnJacobian}
	\ssQ_{\Gamma}
	\colon  \Omega^{\vee} \times \Omega^{\vee} \to \R; \quad \quad \ssQ_{\Gamma}\left(\sum_{e\in E} x_e \ev{e} , \sum_{e\in E} y_e \ev{e} \right) = \sum_{e\in E} \ell(e) x_e y_e
\end{equation}
which is a symmetric, positive-definite, bilinear form on $\Omega^{\vee}$. So $\Jac(\Gamma)$ is identified with the tropical abelian variety $\R^{g} / \Z^{g}$ whose polarization is given by the pairing in Equation \eqref{eq:polarizationOnJacobian}. 
The lattice $L$ dual to $H \cong \Z^{g}$ is identified with $\ssub\Omega_{\Z}^{\vee}$. 
The intermediate Jacobian of $\Jac(\Gamma)$ is therefore given by Proposition \ref{prop:monodromyAbelianVariety}. 

\section{The tropical Ceresa class}

Let $\Gamma$ be a connected genus $g\geq 2$ tropical curve. 
Corresponding to a point $\flat \in \Gamma$ are two maps:
\begin{equation*}
	\Gamma \to \Div^{0}(\Gamma) \quad x\mapsto [x] - [\flat] \quad \text{ and } \quad [\flat] - [x].
\end{equation*}
Composing with the tropical Abel--Jacobi map $\Div^{0}(\Gamma) \to \Jac(\Gamma)$ produces two maps $\Gamma \to \Jac(\Gamma)$. 
Denote the image of $\Gamma$ under these two maps by $\Gamma_{\flat}$ and $\Gamma_{\flat}^{-}$, respectively. 
As these maps contract all separating edges, we may assume $\Gamma$ has no separating edges. 
The \emph{tropical Ceresa cycle} $\nu_{\flat}(\Gamma) \in \cZ_{1}(\Jac(\Gamma))$ is the $1$-cycle 
\begin{equation*}
	\nu_{\flat}(\Gamma) = [\Gamma_{\flat}] - [\Gamma_{\flat}^{-}].
\end{equation*}

We prove below that $\nu_{\flat}(\Gamma)$ is nullhomologous and we produce an explicit bounding $(1,2)$--chain. 
Another $(1,2)$--chain may be found in \cite[Eq.~17]{Ritter}. 
We define the \emph{tropical Ceresa class} of $\Gamma$ based at $\flat$ to be 
\begin{equation*}
	\cc_{\flat}(\Gamma) = \AJ(\nu_{\flat}(\Gamma)) \in \JH_{2,1}(\Jac(\Gamma)). 
\end{equation*}

\subsection{An explicit bounding chain for the Ceresa cycle}
\label{sec:Ceresa-class-explicit}

\begin{figure}
	\includegraphics[width=0.8\textwidth]{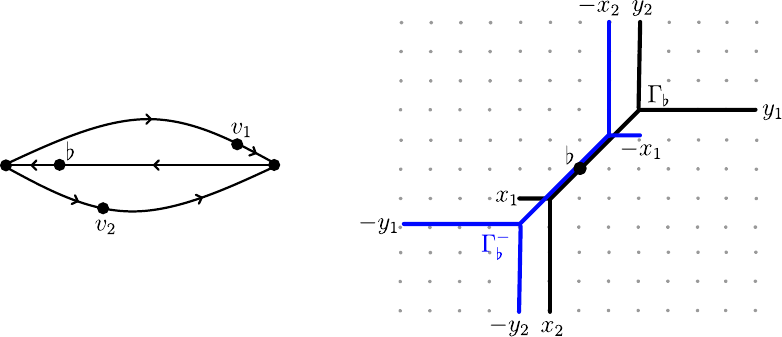}
	\caption{A genus 2 tropical curve $\Gamma$ and the trees $S$ and $S'$ in the universal cover of $\Jac(\Gamma)$}
	\label{fig:CeresaExample}
\end{figure}

Fix a model $(G,\ell)$ of $\Gamma$ with $G = (V,E)$. 
Fix a point $\flat \in V$ and fix an orientation of each edge. 
Let $T = (V,F)$ be a spanning tree of $G$. 
Choose points $D = \{v_{\ve} \st \ve \in F^{c}\}$ such that $v_{\ve}$ lies in the relative interior of the edge $\ve$. 
Denote by $\pi\colon \R^g \to \Jac(\Gamma)$ the universal covering space for $\Jac(\Gamma)$. 
We may assume that $\pi(\mathbf{0}) = \flat$ where $\mathbf{0}\in \R^{g}$ is the origin.  
Let $S$ be the closure of a connected component of $\pi^{-1}(\Gamma \setminus D)$. 
The map $\pi \colon S \to \Gamma$ is one-to-one except over the points in $D$ where the fiber over each of these points has size 2; set $\pi^{-1}(v_{\ve}) = \{x_{\ve}, y_{\ve}\}$.  
See Figure \ref{fig:CeresaExample} for a genus $2$ example. 

 The model $(G,\ell)$ induces a model on $S$. 
 The underlying graph is $G_{S} = (V_S, E_S)$, where $V_{S} = V \cup \{x_{\ve},\ y_{\ve} \st \ve \in \ssF^c\}$ and $E_{S}$ consists of the edges in $F$,  and $2g$ edges obtained by splitting $\ve$ at $v_{\ve}$. 
 The orientations of the edges in $E$ induce orientations on the edges in $E_{S}$. Assume that $x_{\ve}$ has its adjacent edge oriented \emph{away} from it and $y_{\ve}$ has its adjacent edge oriented \emph{towards} it.  

For an oriented edge $e$, denote by $\src{e}$ and $\dst{e}$ the source and target vertex of $e$, respectively. 
Let $\unt{e}\in \rmH_{1,0}(\Jac(\Gamma),\R)$ the unit vector in the direction $\dst{e} - \src{e}$. 
Under the isomorphism $\rmH_{1,0}(\Jac(\Gamma)) \cong \Omega^{\vee}$ in Formula \eqref{eq:evgamma}, $\unt{e}$ corresponds to $\ev{e}/\ell(e)$.
Given points $a_1,\ldots, a_k \in \R^g$, denote by $\pibkk{a_1,\ldots, a_k}$ the image under $\pi$ of the (singular) simplex formed by these points.  
The cycle class of $\Gamma_{\flat}$ is represented by
\begin{equation*}
	Y(\Gamma_{\flat}) \coloneqq \sum_{e\in E_{S}} \Bigl[ \pibkk{\src{e}, \dst{e}}, \unt{e} \Bigr]  \in \rmC_{1,1}(\Jac(\Gamma),\Z). 
\end{equation*}
Define the $(1,1)$ and $(1,2)$--chains 
\begin{equation*}
	Z(\Gamma_{\flat}) 
	= \sum_{\ve\in \ssF^c} \Bigl[ \pibkk{x_{\ve}, y_{\ve}} ,\unt{\ve} \Bigr]
	\quad \text{and} \quad 
	\eta(\Gamma_{\flat}) 
	= \sum_{e\in E_S} 
	\Bigl[\pibkk{ \mathbf{0}, \src{e}, \dst{e} }, \unt{e} \Bigr] - 
	\sum_{\ve\in \ssF^c} 
	\Bigl[ \pibkk{ \mathbf{0}, x_{\ve}, y_{\ve}}, \unt{\ve} \Bigr].
\end{equation*}

\begin{proposition}
\label{prop:YtoZ}
The cycle $Y(\Gamma_{\flat})$ is homologous to $Z(\Gamma_{\flat})$ and
\begin{equation*}
	Y(\Gamma_{\flat}) - Z(\Gamma_{\flat}) = \partial \eta(\Gamma_{\flat}).
\end{equation*}
\end{proposition}

\noindent The following Lemma is a consequence of \cite[Lem.~6.3]{MikhalkinZharkov:TropicalCurves}.

\begin{lemma}
	\label{lem:balancing}
	Given for any $v\in V$, we have
	\begin{equation*}
	\sum_{\src{e} = v} \unt{e} - \sum_{\dst{e}= v} \unt{e} = 0.
	\end{equation*}
\end{lemma}

\begin{proof}[Proof of Proposition \ref{prop:YtoZ}] First we compute
\begin{align*}
	&\partial \Bigl[\pibkk{ \mathbf{0}, \src{e}, \dst{e} }, \unt{e} \Bigr] = \Bigl[\pibkk{\src{e}, \dst{e}} - \pibkk{\mathbf{0}, \dst{e}} + \pibkk{\mathbf{0},\src{e}}, \unt{e} \Bigr], \\
	&\partial \Bigl[\pibkk{\mathbf{0}, x_{\ve}, y_{\ve}}, \unt{\ve} \Bigr] = 
	\Bigl[\pibkk{x_{\ve}, y_{\ve}} - \pibkk{\mathbf{0}, y_{\ve}} + \pibkk{\mathbf{0}, x_{\ve}}, \unt{\ve}\Bigr].
	\end{align*}
Let $L = \{x_\ve, y_{\ve} \st \ve \in F^{c}\} \subset V_{S}$. 
By Lemma \ref{lem:balancing},
\begin{align*}
	\partial\left(\sum_{e\in E_S} \Bigl[\pibkk{ \mathbf{0}, \src{e}, \dst{e} }, \unt{e} \Bigr] \right) 
	=& \ Y(\Gamma_{\flat}) + 
\sum_{v\in V_S\setminus L} \left[ \pibkk{\mathbf{0}, v}, \left(\sum_{\src{e} = v} \unt{e} - \sum_{\dst{e}= v} \unt{e} \right)\right]\\ 
&+ \sum_{\ve \in F^{c}} \Bigl[
\pibkk{\mathbf{0}, x_{\ve}} - \pibkk{\mathbf{0}, y_{\ve}}, \unt{\ve} \Bigr] \\
=& \ Y(\Gamma_{\flat}) + \sum_{\ve \in F^{c}} \Bigl[
\pibkk{ \mathbf{0}, x_{\ve}} - \pibkk{\mathbf{0}, y_{\ve}}, \unt{\ve} \Bigr].
\end{align*}
Finally, we have
\begin{align*}
	\partial\left(\sum_{\ve \in F^{c}}\Bigl[\pibkk{ \mathbf{0}, x_{\ve}, y_{\ve}}, \unt{\ve} \Bigr]  \right) 
	&= \sum_{\ve \in F^c } \Bigl[ \pibkk{x_{\ve}, y_{\ve}}, \unt{\ve} \Bigr] + \sum_{\ve\in F^{c}} \Bigl[ -\pibkk{\mathbf{0}, y_{\ve}} + \pibkk{\mathbf{0}, x_{\ve}}, \unt{\ve} \Bigr] \\
	&= Z(\Gamma_{\flat}) + \sum_{\ve\in F^{c}} \Bigl[-\pibkk{\mathbf{0}, y_{\ve}} + \pibkk{\mathbf{0}, x_{\ve}}, \unt{\ve} \Bigr].
\end{align*}
from which the proposition follows. 	
\end{proof}

The cycle class of $\Gamma_{\flat}^{-}$ is represented by the chain 
\begin{equation*}
	Y(\Gamma_{\flat}^{-}) 
	= \sum_{e\in E_{S}} \Bigl[ \pibkk{-\src{e}, -\dst{e}}, -\unt{e} \Bigr] \in \rmC_{1,1}(\Jac(\Gamma),\Z). 
\end{equation*}
Similar to Proposition \ref{prop:YtoZ}, $Y(\Gamma_{\flat}^{-})$ is homologous to $Z(\Gamma_{\flat}^{-})$ where
\begin{equation*}
	Z(\Gamma_{\flat}^{-}) 
	= \sum_{\ve\in F^{c}} \Bigl[ \pibkk{-x_{\ve}, -y_{\ve}}, -\unt{\ve} \Bigr].
\end{equation*}
In fact, $Y(\Gamma_{\flat}^{-}) - Z(\Gamma_{\flat}^{-}) = \partial (\eta(\Gamma_{\flat}^{-}))$ where 
\begin{gather*}
	\eta(\Gamma_{\flat}^{-}) = 
	\sum_{e\in E_{S}} \Bigl[\pibkk{ \mathbf{0}, -\src{e}, -\dst{e} }, -\unt{e}\Bigr] \; 
	- \; \sum_{\ve \in F^{c}} \Bigl[ \pibkk{ \mathbf{0}, -x_{\ve}, -y_{\ve}}, -\unt{\ve} \Bigr]. 
\end{gather*}
So 
\begin{equation*}
	\partial(\eta(\Gamma_{\flat}) - \eta(\Gamma_{\flat}^{-})) = Y(\Gamma_{\flat}) - Y(\Gamma_{\flat}^{-}) -
	(Z(\Gamma_{\flat}) - Z(\Gamma_{\flat}^{-})).
\end{equation*}
We must find a bounding chain for the nullhomologous cycle $Z(\Gamma_{\flat}) - Z(\Gamma_{\flat}^{-})$. 
Denote by $\rho(a,b,c,d)$ the image under $\pi$ of the rectangle with vertices $a,b,c,d \in \R^{g}$ in counterclockwise order. 
We have that
 \begin{equation*}
	\lambda_{\flat}(\Gamma) = \sum_{\ve \in F^{c}} [\rho(x_{\ve},y_{\ve},-x_{\ve},-y_{\ve}), \unt{\ve}] \quad \text{satisfies} \quad
	\partial(\lambda_{\flat}(\Gamma)) = 
	Z(\Gamma_{\flat}) - Z(\Gamma_{\flat}^{-}).
\end{equation*}
This discussion yields the following theorem. 
\begin{theorem}
\label{thm:connecting-chain-ceresa}
A chain $\xi_{\flat}(\Gamma) \in \rmC_{1,2}(\Jac(\Gamma),\Z)$ such that $\partial \xi_{\flat}(\Gamma) = Y(\Gamma_{\flat}) - Y(\Gamma_{\flat}^{-})$ is given by
\begin{equation*}
	\xi_{\flat}(\Gamma) = \eta(\Gamma_{\flat}) - \eta(\Gamma_{\flat}^{-}) + \lambda_{\flat}(\Gamma). 
\end{equation*}
\end{theorem}

\begin{figure}	
\includegraphics[width=\textwidth]{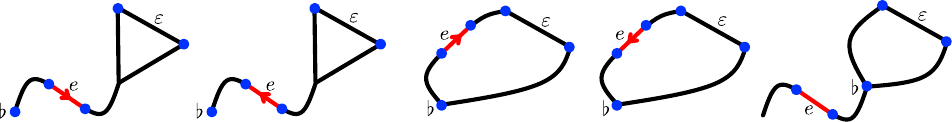}
	\caption{The values of $\sgn_{T}^{\flat}(e,\ve)$ are, left to right: $2$, $-2$, $1$, $-1$, $0$.}
	\label{fig:sgn_flat}
\end{figure}

\noindent Using the bounding chain $\xi_{\flat}(\Gamma)$ above, we derive a useful and explicit description of the Ceresa class $\cc_{\flat}(\Gamma)$. 
For $e\in F$ and $\ve \in \ssF^c$, we define the sign function  $\sgn_{T}^{\flat}(e,\ve)$ in terms of the paths in $T$ from $\flat$ to $\src{\ve}$ or $\dst{\ve}$ passing through $e$.
Specifically, 
\begin{equation*}
	\sgn_{T}^{\flat}(e,\ve) = 
	\begin{cases}
		2 & \text{ if } e \text{ is on the path from } \flat \text{ to both } \src{\ve} \text{ and } \dst{\ve}, \text{ oriented away from } \flat, \\
		1 & \text{ if } e \text{ is on the path from } \flat \text{ to one of } \src{\ve} \text{ or } \dst{\ve}, \text{ oriented away from } \flat, \\
		0 & \text{ if } e \text{ is not on the path from } \flat \text{ to either } \src{\ve} \text{ or } \dst{\ve},\\
		-1 & \text{ if } e \text{ is on the path from } \flat \text{ to one of } \src{\ve} \text{ or } \dst{\ve}, \text{ oriented towards } \flat, \\
		-2 & \text{ if } e \text{ is on the path from } \flat \text{ to both } \src{\ve} \text{ and } \dst{\ve}, \text{ oriented towards } \flat. \\
	\end{cases}
\end{equation*}
See Figure \ref{fig:sgn_flat}. 
The cycles $\pibkk{x_{\ve},y_{\ve}}$ and $\pibkk{-y_{\ve},-x_{\ve}}$ define the same class in $\rmH_1(\Jac(\Gamma), \Z)$; we denote this class by $\cycb{\ve}$. 
Under the isomorphism $\rmH_{1}(\Jac(\Gamma), \Z) \cong \rmH_{1}(\Gamma,\Z)$, the cycle $\cycb{\ve}$ corresponds to the homology class of the fundamental cycle contained in $T\cup\ve$. 

\begin{theorem}
\label{thm:ceresa-explicit}
Given a tropical curve $\Gamma$ with model $(G = (V,E),\ell)$, a point $\flat\in V$, and spanning tree $T = (V, F)$, the tropical Ceresa class $\cc_{\flat}(\Gamma)$ is given by
\begin{equation*}
	\cc_{\flat}(\Gamma) = \sum_{e\in F \\ \ve \in \ssF^c} \sgn_{T}^{\flat}(e,\ve) \, \ell(e) \, \cycb{\ve} \otimes (\unt{\ve} \wedge \unt{e}). 
\end{equation*}
\end{theorem}
\noindent We use the description of $\JH_{2,1}(\Jac(\Gamma))$ given in Proposition \ref{prop:monodromyAbelianVariety}.
\begin{proof}[Proof of Theorem~\ref{thm:ceresa-explicit}]
By Proposition \ref{prop:AJ-boundaries-zero}, $\cc_{\flat}(\Gamma)$ only depends on the bounding chain of a representative of $\cl(\nu_{\flat}(\Gamma))$ up to $\rmB_{2,1}(\Jac(\Gamma))$. 
A the chain $\xi_{\flat}(\Gamma)$ bounds $\cl(\nu_{\flat}(\Gamma))$ by Theorem \ref{thm:connecting-chain-ceresa}.
To compute $\rmN(\xi_{\flat}(\Gamma))$, we must fix points $\op{P}$ for each codimension zero and one cell $P$. 
With
\begin{gather*}
\sigma_{e} = \pibkk{\mathbf{0}, \src{e}, \dst{e}}, \quad
\sigma_{e}' = \pibkk{\mathbf{0}, -\src{e}, -\dst{e}}, \quad
\tau_{\ve} = \pibkk{\mathbf{0}, x_{\ve}, y_{\ve}},  \quad
\tau_{\ve}'= \pibkk{\mathbf{0}, -x_{\ve}, -y_{\ve}}, \\
\rho_{\ve} = \rho(x_{\ve}, y_{\ve}, -x_{\ve}, -y_{\ve}),
\end{gather*}
set
\begin{gather*}
	\op{\sigma_{e}} = \dst{e}, \quad
	o_{\sigma_{e}'} = -\dst{e}, \quad
	\op{\tau_{\ve}} = y_{\ve}, \quad 
	o_{\tau_{\ve}'} = -y_{\ve}, \quad 
	o_{\rho_{\ve}} = y_{\ve}, \\ 
	o_{\pibkk{\mathbf{0}, v}} = v, \quad 
	o_{\pibkk{\pm x_{\ve}, \pm y_{\ve}}} = \pm y_{\ve}, \quad 
	o_{\pibkk{y_{\ve}, -x_{\ve}}} = y_{\ve}, \quad  
	o_{\pibkk{-y_{\ve}, x_{\ve}}} = x_{\ve}.
\end{gather*}
Using this, we compute 
\begin{align*}
&\begin{array}{ll}
	\rmN \left([\sigma_{e}, \unt{e}]\right) = 0, 
	& \rmN\left([\sigma_{e}', -\unt{e}]\right) = 0, \\
	\rmN \left([\tau_{\ve}, \unt{\ve}] \right) = \Bigl[\pibkk{\mathbf{0}, x_{\ve}}, \cycb{\ve} \wedge \unt{\ve} \Bigr], 
	& \rmN \left([\tau_{\ve}', -\unt{\ve}]\right) = \Bigl[\pibkk{\mathbf{0}, -x_{\ve}}, \cycb{\ve}\wedge \unt{\ve} \Bigr], \\
\end{array} \\
	&\rmN ([\rho_{\ve}, \unt{\ve}]) = 
	\Bigl[\pibkk{-x_{\ve},-y_{\ve}}, (2y_{\ve}-\mathbf{0}) \wedge \unt{\ve} \Bigr] + \Bigl[\pibkk{-y_{\ve}, x_{\ve}}, \cycb{\ve} \wedge \unt{\ve} \Bigr].
\end{align*}
Then,
\begin{align*}
	\rmN(\eta(\Gamma_{\flat}) - \eta(\Gamma_{\flat}^{-})) &= \sum_{\ve\in F^{c}} \Bigl[
	-\pibkk{\mathbf{0}, x_{\ve}} + \pibkk{\mathbf{0}, -x_{\ve}}, \cycb{\ve} \wedge \unt{\ve} \Bigr] \\
	&\equiv  \sum_{\ve\in F^{c}} \Bigl[
	\pibkk{x_{\ve}, -x_{\ve}}, \cycb{\ve} \wedge \unt{\ve} \Bigr] \mod \rmB_{2,1}(\Jac(\Gamma)).
\end{align*}
	
\noindent Altogether, we have
\begin{align*}
	\rmN(\xi_{\flat}(\Gamma)) 
	\equiv \sum_{\ve\in F^{c}} & \left(\Bigl[\pibkk{x_{\ve}, -x_{\ve}}, \cycb{\ve} \wedge \unt{\ve}\Bigr] + \Bigl[\pibkk{-x_{\ve}, -y_{\ve}}, (2y_{\ve} - \mathbf{0}) \wedge \unt{\ve}\Bigr] \right. \\
	& + \left. \Bigl[\pibkk{-y_{\ve}, x_{\ve}}, \cycb{\ve} \wedge \unt{\ve}\Bigr] \right)
	\mod \rmB_{2,1}(\Jac(\Gamma)).
\end{align*}
As $\pibkk{x_{\ve}, -x_{\ve}} + \bkk{-y_{\ve}, x_{\ve}} \equiv  \pibkk{-y_{\ve}, -x_{\ve}}$ modulo $\rmB_{0,1}(\Jac(\Gamma))$, we may combine the three terms in the above sum to get 
\begin{align*}
	\rmN(\xi_{\flat}(\Gamma)) \equiv \sum_{\ve\in F^{c}} \Bigl[\pibkk{-y_{\ve}, -x_{\ve}}, (\cycb{\ve} - 2y_{\ve}) \wedge \unt{\ve} \Bigr] 
	\equiv \sum_{\ve\in F^{c}} \Bigl[\cycb{\ve}, (-y_{\ve} - x_{\ve}) \wedge \unt{\ve} \Bigr].
\end{align*}
Because
\begin{equation*}
	-y_{\ve} - x_{\ve} = (-y_{\ve} - \mathbf{0}) - (x_{\ve} - \mathbf{0}) = -(y_{\ve} - \mathbf{0}) - (x_{\ve} - \mathbf{0}),
\end{equation*}
we have that
\begin{equation*}
	\cc_{\flat}(\Gamma) = \sum_{\ve\in F^{c}}  [\cycb{\ve},  \unt{\ve} \wedge (x_{\ve} - \mathbf{0})  + \unt{\ve} \wedge (y_{\ve} - \mathbf{0})].
\end{equation*}
Let $A_{\ve}$, respectively $B_{\ve}$, be the edges in $F$ contained in the unique path in $S$ from $\mathbf{0}$ to $x_{\ve}$, respectively from $\mathbf{0}$ to $y_{\ve}$. 
For $e\in A_{\ve}$, let $\alpha_{e} \in \{\pm 1\}$ be $1$ if $e$ is oriented away from $\flat$ and $-1$ if $e$ is oriented towards $\flat$. Define $\beta_{e} \in \{\pm 1\}$ for $e\in B_{\ve}$ similarly. 
Then,
\begin{equation*}
	(x_{\ve} - \flat) = k\unt{\ve}
 	+ \sum_{e\in A_{\ve}} \alpha_{e}  \ \ell(e) \unt{e}, \quad
 	(y_{\ve} - \flat) = (\ell(\ve) - k)\unt{\ve}
 	+ \sum_{e\in B_{\ve}} \beta_{e} \ \ell(e) \unt{e} 
\end{equation*}
for some $0\leq k \leq \ell(\ve)$. 
The theorem now follows from the description of $\JH_{2,1}(\Jac(\Gamma))$ as in Proposition \ref{prop:monodromyAbelianVariety}.   
\end{proof}

\subsection{Dependence on the basepoint}
Given two points $\flat, \rflat \in \Gamma$, we consider the difference between $\cc_{\flat}(\Gamma)$ and $\cc_{\rflatind}(\Gamma)$. 
Fix a basis of cycles $\cycb{1}, \ldots, \cycb{g}$ of $\ssH_{\Z} \subset \Omega^{\vee}$, and let $\unt{1}, \ldots, \unt{g} \in \ssub\Omega_{\Z}^{\vee}$ be a dual basis with respect to $Q$ (so $Q(\cycb{i}, \unt{j}) = \delta_{ij}$). 
Define the homology class
\begin{equation}
\label{eq:trop-omega}
	\omega  = [\cycb{1}, \unt{1}] + \cdots + [\cycb{g}, \unt{g}] \in \rmH_{1,1}(\Jac(\Gamma)). 
\end{equation}

\begin{proposition}
The homology class $\omega \in \rmH_{1,1}(\Jac(\Gamma))$ is independent of basis. 
Furthermore,
\begin{equation*}
	\omega = \cl(\Gamma_{\flat}).
\end{equation*}
\end{proposition}

\begin{proof}
Suppose $\cycb{1}', \ldots, \cycb{g}', \unt{1}', \ldots, \unt{g}'$ is another pair of dual bases. 
There is linear automorphism $C\colon \ssH_{\Z} \to \ssH_{\Z}$ such that $C\cycb{i} = \cycb{i}'$. 
Then $\unt{i}'=C^{-T}\unt{i}$ where  $C^{T}\colon \Omega^{\vee} \to \Omega^{\vee}$ is the unique linear automorphism such that $Q(Cx,y) = Q(x,C^{T}y)$.  
Say $C = [c_{ij}]$ and $C^{-T} = [d_{ij}]$. Then
\begin{align*}
	\sum_{i=1}^{g} [\cycb{i}', \unt{i}'] &= \sum_{i=1}^{g} [C\cycb{i}, C^{-T}\unt{i}] = \sum_{i=1}^{g} \sum_{k,\ell}[c_{ik}\cycb{k}, d_{i\ell}\unt{\ell}] = \sum_{k, \ell} [\cycb{k}, \unt{\ell}] \sum_{i=1}^{g} c_{ik} d_{i\ell} \\ 
	&= \sum_{k,\ell} [\cycb{k}, \unt{\ell}] (C^{T}C^{-T})_{k\ell}  = \sum_{k,\ell} [\cycb{k}, \unt{\ell} ] \delta_{k\ell} = \omega.
\end{align*} 
The last statement follows from this and Proposition \ref{prop:YtoZ}. 
\end{proof}

\begin{theorem}
\label{thm:dependenceOnBasePoint}
	We have
	\begin{equation*}
	\cc_{\flat}(\Gamma) - \cc_{\rflatind}(\Gamma)
		 = -2\AJ(\flat - \rflat) \wedge \omega.
	\end{equation*}
\end{theorem}

\begin{proof}
This theorem may be deduced from Theorem \ref{thm:ceresa-explicit}, but we provide a more direct proof here. 
First, consider the nullhomologous cycle $\Gamma_{\flat} - \Gamma_{\rflatind}$. Representatives of $\cl(\Gamma_{\flat})$ and $\cl(\Gamma_{\rflatind})$ are given by
\begin{equation*}
	Y(\Gamma_{\flat}) = \sum_{e\in E} \Bigl[\pibkk{\src{e}, \dst{e}}, \unt{e} \Bigr] \quad \text{and} \quad Y(\Gamma_{\rflatind}) = \sum_{e\in E} \Bigl[\pibkk{\src{e}+\rflat, \dst{e}+\rflat}, \unt{e} \Bigr],
\end{equation*}
respectively. The $(2,1)$--chain
\begin{equation*}
\eta = \sum_{e\in E} [\eta_{e}, \unt{e}] \quad  \text{ where }\quad 
	\eta_e = \rho(\src{e}, \dst{e}, \dst{e} + \rflat, \src{e} + \rflat )	
\end{equation*}
satisfies $\partial \eta = Y(\Gamma_{\flat}) - Y(\Gamma_{\rflatind})$. For each vertex $v\in V$ and edge $e\in E$ set 
\begin{equation*}
	o_{\pibkk{v, v+\rflatind}} = v + \rflat, \quad
	o_{\bkk{\src{e}, \dst{e}}} = \dst{e}, \quad
	o_{\bkk{\src{e}+\rflatind, \dst{e}+\rflatind}} = \dst{e}+\rflat, \quad
	o_{\eta_e} = \dst{e} + \rflat.		
\end{equation*}
Then,
\begin{equation*}
	\rmN(\eta_e) 
	= [\pibkk{\src{e}, \dst{e}},(\rflat - \flat) \wedge \unt{e}] - [\pibkk{\src{e}, \src{e}+\rflat}, \ell(e)\unt{e} \wedge \unt{e}] 
	= [\pibkk{\src{e}, \dst{e}}, (\rflat - \flat) \wedge \unt{e}],
\end{equation*}
and so
\begin{equation*}
	\rmN(\eta) = \sum_{e\in E} [\pibkk{\src{e}, \dst{e}}, (\rflat-\flat) \wedge \unt{e}] =  \AJ(\rflat-\flat) \wedge Y(\Gamma_{\flat}).
\end{equation*}
Since $Y(\Gamma_{\flat})$ and $\omega$ are homologous, we have
\begin{equation*}
	\AJ(\Gamma_{\flat} - \Gamma_{\rflatind}) =  -\AJ(\flat-\rflat) \wedge \omega.
\end{equation*}
Similarly,
\begin{equation*}
	\AJ(\Gamma_{\flat}^{-} - \Gamma_{\rflatind}^{-}) =  \AJ(\flat-\rflat) \wedge \omega,
\end{equation*}
and therefore
\begin{equation*}
	\cc_{\flat}(\Gamma) - \cc_{\rflatind}(\Gamma)
	= \AJ(\Gamma_{\flat} - \Gamma_{\rflatind}) - \AJ(\Gamma_{\flat}^{-} - \Gamma_{\rflatind}^{-})
	= -2 \AJ(\flat-\rflat) \wedge \omega
\end{equation*}
from which the theorem follows. 
\end{proof}

\subsection{The unpointed Ceresa class}

With Theorem \ref{thm:dependenceOnBasePoint} in mind, we make the following definition. 

\begin{definition}
Define 
\begin{equation*}
	\overline{\JH}_{p+1,p}(\Jac(\Gamma)) = \frac{\rmH_{p+1,p}(\Jac(\Gamma),\R)}{\omega \wedge \rmH_{p,p-1}(\Jac(\Gamma),\R) + \rmL_{p+1,p}(\Jac(\Gamma))}.
\end{equation*}
The \emph{unpointed} Ceresa class $\overline{\cc}(\Gamma)$ is the image of $\cc_{\flat}(\Gamma)$ in $\overline{\JH}_{2,1}(X)$.
\end{definition}
\noindent By Theorem \ref{thm:dependenceOnBasePoint}, the class $\cc_{\flat}(\Gamma)$ is independent of $\flat$. 

Next, we derive a new formula for $\overline{\cc}(\Gamma)$.  
Suppose $(G,\ell)$ is a model of $\Gamma$ with $G = (V,E)$, fix an orientation of each edge, and let $T = (V,F)$ be a spanning tree of $G$. Given $\ve\in \ssF^c$, there is a unique cycle contained in $T\cup \ve$; call this cycle $\gamma_{T}(\ve)$. 
Given $e\in F$ and $\ve \in \ssF^c$, define 
\begin{equation*}
	\sgnbar_{T}(e,\ve) = 
	\begin{cases}
		1 & \text{ if } e \notin \gamma_{T}(\ve) \text{ and } e \text{ points towards } \ve,\\
		-1 & \text{ if } e \notin \gamma_{T}(\ve) \text{ and } e \text{ points away from } \ve, \\
		0 & \text{ if } e \in \gamma_{T}(\ve).
	\end{cases}
\end{equation*}

\begin{theorem}
\label{thm:unpointedCeresaFormula}
Given a tropical curve $\Gamma$ with model $(G = (V,E), \ell)$ and spanning tree $T = (V,F)$, its unpointed tropical Ceresa class is 
\begin{equation*}
	\overline{\cc}(\Gamma) = \sum_{e\in F \\ \ve \in \ssF^c}  \sgnbar_{T}(e,\ve) \, \ell(e) \, \cycb{\ve} \otimes (\unt{\ve}\wedge \unt{e}) \quad \text{ in }\quad \overline{\JH}_{2,1}(\Jac(\Gamma)).
\end{equation*}
\end{theorem}

\begin{proof}
Fix a basepoint $\flat\in V$. Given $e\in F$ and $\ve\in \ssF^c$, a case-by-case analysis yields 
\begin{equation*}
	\sgn_{T}^{\flat}(e,\ve) - \sgnbar_{T}(e,\ve) = 
	\begin{cases}
		1 & \text{ if }	e \text{ points away from } \flat \text{ in } T, \\
		-1 & \text{ if } e \text{ points towards  } \flat \text{ in } T.
	\end{cases}
\end{equation*}
In particular, this expression is independent of $\ve \in \ssF^c$; call it $\alpha_{e}$. 
Using Theorem \ref{thm:ceresa-explicit},  we have
	\begin{equation*}
		\cc_{\flat}(\Gamma) - \sum_{e\in F \\ \ve \in \ssF^c}  \sgnbar_{T}(e,\ve) \, \ell(e) \, \cycb{\ve}\otimes (\unt{\ve}\wedge \unt{e}) = \omega \wedge \sum_{e\in F} \ell(e) \alpha_e \,  \unt{e},
	\end{equation*} 
	as required. 
\end{proof}

\subsection{The Ceresa--Zharkov class and algebraic equivalence}
\label{sec:CZ-class}
Given a tropical curve $\Gamma$ and a point $\flat \in \Gamma$, define the class $\mathbf{w}(\Gamma)$ in $\rmQ_{3,0}(\Jac(\Gamma))$ by
\begin{equation*}
	\mathbf{w}(\Gamma) = \rmN(\cc_{\flat}(\Gamma)) \in \rmQ_{3,0}(\Jac(\Gamma)). 
\end{equation*}
By Theorem \ref{thm:obstructionAlgebraicEquivalence}, if $\mathbf{w}(\Gamma) \neq 0$ in $\rmQ_{3,0}(\Jac(\Gamma))$, then the tropical Ceresa cycle $\nu_{\flat}(\Gamma)$ is not algebraically equivalent to 0. 
We expect that this recovers \cite[Prop.~1.1]{Ritter}. 
Using the results from Appendix \ref{sec:MoritaClass}, we see that this class generalizes the \emph{Ceresa--Zharkov class} of a tropical curve as defined in \cite{CoreyLi} to curves with arbitrary edge-lengths. 

\medskip

Theorems \ref{thm:ceresa-explicit} and \ref{thm:unpointedCeresaFormula} give efficient ways of calculating the Ceresa class of a tropical curve. 
We illustrate this with two examples, the complete graph on $4$ vertices $K_4$ and the trivalent loop of 3 loops $TL_3$.  

\subsection{The complete graph on 4 vertices $K_4$}
\label{sec:K4}
\begin{figure}
	\centering	
	\begin{minipage}{0.45\textwidth}
	\centering
		\includegraphics[width=0.5\textwidth]{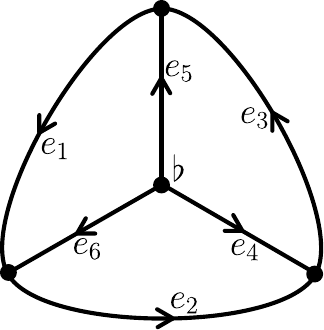}
	\end{minipage}
	\begin{minipage}{0.45\textwidth}
	\centering
		\includegraphics[width=0.5\textwidth]{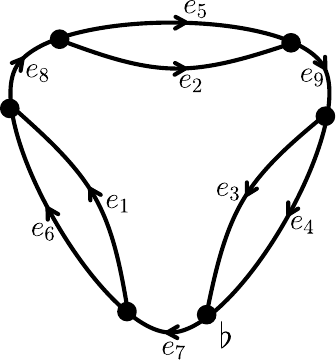}
	\end{minipage}
		\caption{A genus 3 tropical curve whose underlying graph is $K_4$, and a genus 4 tropical curve whose underlying graph is $TL_3$}
		\label{fig:K4}
	\end{figure}

Suppose $\Gamma$ is the tropical curve whose underlying graph is $K_{4}$, the left-hand graph in Figure \ref{fig:K4}. 
For each edge $e_i$, denote by $\ell_i = \ell(e_i)$ and $\unt{i} = \unt{e_i}$.  
Let $T$ be the spanning tree formed by the edges $e_4, e_5, e_6$. The vectors $\unt{1}, \unt{2}, \unt{3}$ form a basis of $\rmH_{1,0}(\Gamma,\Z)$, and we have
\begin{equation*}
	\unt{4} = \unt{3} - \unt{2}, \quad \unt{5} = \unt{1} - \unt{3}, \quad \unt{6} = \unt{2} - \unt{1}.
\end{equation*}
For $i=1,2,3$, let $\cycb{i} = \cycb{e_i}$. 
The matrix $\ssQ_{\Gamma}$ with respect to the basis $\cycb{1}, \cycb{2}, \cycb{3}$ (i.e., the matrix whose $(i,j)$--th entry is $\ssQ_{\Gamma}(\cycb{i}, \cycb{j})$) is given by 
\begin{equation*}
	\ssQ_{\Gamma} = 
	\begin{bmatrix}
	\ell_1 + \ell_5 + \ell_6	 & -\ell_6 & -\ell_5 \\
	-\ell_6 & \ell_2 + \ell_4 + \ell_6 & -\ell_4 \\
	-\ell_5 & -\ell_4 & \ell_3 + \ell_4 + \ell_5
	\end{bmatrix}.
\end{equation*} 
The columns of this matrix describe the coefficients of $\cycb{i}$ with respect to the basis $\unt{1}, \unt{2}, \unt{3}$. 
Explicitly:
\begin{align*}
	\cycb{1} &= \ell_1\unt{1} + \ell_5\unt{5} -\ell_6\unt{6} 
	= (\ell_1 + \ell_5 + \ell_6) \unt{1} - \ell_6 \unt{2} - \ell_{5} \unt{3}, \\
	 \cycb{2} &= \ell_2\unt{2} - \ell_4\unt{4} +\ell_6\unt{6} 
	=  -\ell_6 \unt{1} + (\ell_2 + \ell_4 + \ell_6) \unt{2} - \ell_4\unt{3}, \\
	\cycb{3} &= \ell_3\unt{3} + \ell_4\unt{4} -\ell_5\unt{5} 
	= -\ell_{5} \unt{1} -\ell_{4} \unt{2} + (\ell_3 + \ell_4 + \ell_5) \unt{3}.
\end{align*}
This description gives a practical way to compute the monodromy map on the Jacobian via Proposition \ref{prop:monodromyAbelianVariety}. 
Indeed, we have that $\rmL_{2,1}(\Jac(\Gamma))$ is generated by the forms
\begin{equation*}
	\rmN(\aab{i}{j}{k}) = \cycb{i}\otimes \cycb{j} \wedge \unt{k} - \cycb{j}\otimes \cycb{i} \wedge \unt{k} 
\end{equation*}
for $1\leq i<j\leq 3$ and $1\leq k\leq 3$, which can be expressed in the basis 
\begin{equation*}
	\abb{i}{j}{k} \quad \text{for} \quad 1\leq i\leq 3, \; 1\leq j<k\leq 3
\end{equation*}
of $\R^{3}\otimes \wedge^2\R^{3}$ using the descriptions of $\cycb{i}$ and $\unt{j}$ above. 
Using Theorem \ref{thm:ceresa-explicit} we compute the pointed Ceresa class to be
	\begin{align}
	\label{eq:K4Pointed}
	\begin{array}{rll}
		\cc_{\flat}(\Gamma) &=& 
		  \ell_4( \cycb{2} \otimes (\unt{2} \wedge \unt{4}) + \cycb{3}\otimes (\unt{3} \wedge \unt{4}) ) \\
		& & + \; \ell_5( \cycb{1}\otimes (\unt{1} \wedge \unt{5}) + \cycb{3}\otimes (\unt{3} \wedge \unt{5}) ) \\
		& & + \; \ell_6( \cycb{1} \otimes (\unt{1} \wedge \unt{6}) +  \cycb{2} \otimes (\unt{2} \wedge \unt{6})).
	\end{array}
	\end{align} 
The formula in Theorem \ref{thm:unpointedCeresaFormula} yields the same representative for $\overline{\cc}(\Gamma)$. 
We compute in $\overline{\JH}_{2,1}(\Jac(\Gamma))$
\begin{align*}
	\overline{\cc}(\Gamma) &= 
	\cc_{\flat}(\Gamma) - \omega \wedge (\ell_4\unt{4} + \ell_5\unt{5} + \ell_6\unt{6}) + \rmN([\cycb{1} \wedge \cycb{2}, \unt{1}])  \\
	&= \ell_2 \cycb{1} \otimes (\unt{1} \wedge \unt{2})  
	- \ell_5 (\cycb{2} \otimes (\unt{1} \wedge \unt{2}) + \cycb{2} \otimes (\unt{2} \wedge \unt{3}) - \cycb{2}\otimes (\unt{1} \wedge \unt{3})),
\end{align*}
which recovers the computation in \cite[Ex.~7.2]{CoreyEllenbergLi}. 
Finally, we compute $\rmQ_{3,0}(\Jac(\Gamma))$ and the class $\mathbf{w}(\Gamma)$ in $\rmQ_{3,0}(\Jac(\Gamma))$. 
The subgroup $\rmK_{3,0}(\Jac(\Gamma))$ of $\wedge^{3}\R^3 \cong \R$ is generated by 
\begin{align*}
	\rmN^2(\aab{i}{j}{k}) = 2 \ \cycb{i} \wedge \cycb{j} \wedge \unt{k}
\end{align*}
for $1\leq i<j\leq 3$ and $1\leq k\leq 3$, and	
\begin{align*}
	\mathbf{w}(\Gamma)  
	&= -(\ell_5\ell_6 + \ell_4\ell_6 + \ell_4\ell_5)\ \unt{1} \wedge \unt{2} \wedge \unt{3} 
	= \ell_1 \ell_4\, \unt{1} \wedge \unt{2} \wedge \unt{3} \quad \text{ in } \; \rmQ_{3,0}(\Jac(\Gamma)).
\end{align*}
This recovers the computation in \cite{CoreyLi, Ritter, Zharkov}. 

\subsection{The trivalent loop of 3 loops $TL_3$}

Suppose $\Gamma$ is the tropical curve whose underlying graph is $TL_{3}$, the right-hand graph in Figure \ref{fig:K4}. 
As in the $K_4$ example,  for each edge $e_i$, denote by $\ell_i = \ell(e_i)$ and $\unt{i} = \unt{e_i}$.  
Let $T$ be the spanning tree formed by the edges $e_5$, $e_6$, $e_7$, $e_{8}$, $e_{9}$. 
Then $\unt{1},\unt{2},\unt{3},\unt{4}$ forms a basis of $\rmH_{1,0}(\Gamma,\R)$ and we have
\begin{equation*}
	\unt{5} = - \unt{2} + \unt{3} + \unt{4}, \quad 
	\unt{6} 	= -\unt{1} + \unt{3} + \unt{4}, \quad 
	\unt{7} = \unt{8} = \unt{9} = \unt{3} + \unt{4}.
\end{equation*}
For $i=1,2,3,4$, let $\cycb{i} = \cycb{e_i}$, which forms a basis of $\rmH_{0,1}(\Gamma)$. 
The matrix $\ssQ_{\Gamma}$ with respect to this basis is given by 
\begin{equation*}
	\ssQ_{\Gamma} = 
	\begin{bmatrix}
	\ell_1 + \ell_6 & 0 & -\ell_6 & -\ell_6 \\
	0 & \ell_2 + \ell_5 & -\ell_5 & -\ell_5 \\
	-\ell_{6} & -\ell_{5} & \ell_{3} + \ell_{5} + \ell_{6} + \ell_{7} + \ell_{8} + \ell_{9} &  \ell_{5} + \ell_{6} + \ell_{7} + \ell_{8} + \ell_{9} \\
	-\ell_{6} & -\ell_{5} &  \ell_{5} + \ell_{6} + \ell_{7} + \ell_{8} + \ell_{9} & \ell_{4} + \ell_{5} + \ell_{6} + \ell_{7} + \ell_{8} + \ell_{9}
	\end{bmatrix}.
\end{equation*} 
Using Theorem \ref{thm:ceresa-explicit} we compute the pointed Ceresa class $\cc_{\flat}(\Gamma)$ in $\JH_{2,1}(\Jac(\Gamma))$ to be
\begin{align*}
	\cc_{\flat}(\Gamma) &= 
	\ell_{5} (\cycb{1} \otimes (\unt{1} \wedge \unt{5}) 
	+ \cycb{3}\otimes (\unt{3} \wedge \unt{5}) 
	+ \cycb{4}\otimes (\unt{4} \wedge \unt{5})) \\
	&+ \ell_6 (
	\cycb{1}\otimes (\unt{1} \wedge \unt{6}) 
	+2\cycb{2}\otimes (\unt{2} \wedge \unt{6}) 
	+\cycb{3}\otimes (\unt{3} \wedge \unt{6}) 
	+\cycb{4}\otimes (\unt{4} \wedge \unt{6})) \\
	&+ \ell_7 (
	2\cycb{1}\otimes (\unt{1} \wedge \unt{7}) 
	+2\cycb{2}\otimes (\unt{2} \wedge \unt{7}) 
	+\cycb{3}\otimes (\unt{3} \wedge \unt{7}) 
	+\cycb{4}\otimes (\unt{4} \wedge \unt{7})) \\
	&+ \ell_8 (
	2\cycb{2}\otimes (\unt{2} \wedge \unt{8}) 
	+\cycb{3}\otimes (\unt{3} \wedge \unt{8}) 
	+\cycb{4}\otimes (\unt{4} \wedge \unt{8})) \\
	&+ \ell_9 (
	\cycb{3}\otimes (\unt{3} \wedge \unt{9}) 
	+\cycb{4}\otimes (\unt{4} \wedge \unt{9})).
\end{align*}
The formula for $\overline{\cc}(\Gamma)$ in $\overline{\JH}_{2,1}(\Jac(\Gamma))$ from Theorem \ref{thm:unpointedCeresaFormula} yields
\begin{align*}
	\overline{\cc}(\Gamma) = 
	& -\ell_5 \cycb{1} \otimes (\unt{1} \wedge \unt{5}) + \ell_6  \cycb{2} \otimes (\unt{2} \wedge \unt{6}) -\ell_7 (
		\cycb{1} \otimes (\unt{1} \wedge \unt{7}) 
		+\cycb{2} \otimes (\unt{2} \wedge \unt{7}) ) \\
	& +\ell_8 (
		\cycb{1} \otimes (\unt{1} \wedge \unt{8}) 
		- \cycb{2} \otimes (\unt{2} \wedge \unt{8})) +\ell_9 (
		\cycb{1} \otimes (\unt{1} \wedge \unt{9}) 
		+ \cycb{2} \otimes (\unt{2} \wedge \unt{9})).
\end{align*}
This recovers the calculation in \cite[Ex.~7.6]{CoreyEllenbergLi}.
The class $\mathbf{w}(\Gamma)$ in $\rmQ_{3,0}(\Jac(\Gamma))$ is represented by
\begin{align*}
	\mathbf{w}(\Gamma)  = \ &  \ell_5\ell_6\ \unt{1} \wedge \unt{5} \wedge  \unt{6} - \ell_5\ell_6\ \unt{2} \wedge \unt{5} \wedge  \unt{6} 
	- \ell_6\ell_7\  \unt{1} \wedge \unt{6} \wedge \unt{7} - \ell_5\ell_7\ \unt{2} \wedge \unt{5} \wedge \unt{7} \\
	& + \ell_6\ell_8\ \unt{1} \wedge \unt{6} \wedge \unt{8} - \ell_5\ell_8\ \unt{2} \wedge \unt{5} \wedge \unt{8} 
	+ \ell_6\ell_9\ \unt{1} \wedge \unt{6} \wedge \unt{7} + \ell_5\ell_9\ \unt{2} \wedge \unt{5} \wedge \unt{9} \\
	= \ & -2 \ell_5\ell_6\ (\unt{1} \wedge \unt{2} \wedge  \unt{3} + \unt{1} \wedge \unt{2} \wedge  \unt{4}). 
\end{align*}
When we set the edge lengths $\ell_{7}$, $\ell_{8}$, $\ell_{9}$ each equal to $0$, the underlying graph of $\Gamma$ is the loop of three loops $L_3$, and this  recovers the calculation from \cite[Prop.~5.9]{CoreyLi} and \cite[Ex.~4.2]{Ritter}.

\section{Further results and questions}

In this final section, we prove some complementary results and raise further questions. 

\subsection{Torsion Ceresa class}
\label{sec:torsion-Ceresa}
Beauville in \cite{Beauville21} produces a nonhyperelliptic curve whose Ceresa class is torsion (later with Schoen, they show that the Ceresa cycle of this curve is torsion modulo algebraic equivalence \cite{BeauvilleSchoen}). 
Separate examples may be found in \cite{BLLS,  LagaShnidman-picard, LagaShnidman-vanishing, LilienfeldtShnidman}. 
Determining whether the locus in the moduli space of curves $\cM_{g}$ of curves with torsion Ceresa class is dense has seen significant recent activity. 
It is known that the positive-dimensional components of this locus are not dense in $\cM_g$, see \cite{GaoZhang, KerrTayou}. 
For tropical curves, we have the following. 

\begin{theorem}
	The set of tropical curves with torsion tropical Ceresa class is dense in $M_{g}^{\trop}$.
\end{theorem}

\begin{proof} Let $G = (V,E)$ be a stable graph of genus $g$, and $\ell \colon E \to \Q_{>0}$ an edge length function taking only rational values. 
Let $\Gamma$ be the corresponding metric graph of genus $g$. 
The monodromy operator for $\Gamma$ has rational coefficients. 
By Theorem \ref{thm:torsion-AJ}, the Ceresa class of $\Gamma$ is torsion. 
Since these tropical curves form a dense subset of the tropical moduli space $M_{g}^{\trop}$, this proves the theorem. 
\end{proof}

\subsection{Ceresa cycle and hyperelliptic curves}

If $X$ is a hyperelliptic curve and $\flat\in X$ is a Weierstrass point, then the Ceresa cycle $[X_{\flat}] - [X_{\flat}^{-}] = 0$ as an algebraic cycle. 
In particular, this implies that $[X_{\flat}] - [X_{\flat}^{-}]$ is algebraically equivalent to $0$ for any $\flat$. 
On the other hand, if $X$ is a curve and $\flat\in X$ have the property that $[X_{\flat}] - [X_{\flat}^{-}]$ is trivial (rationally, algebraically, or under the Abel--Jacobi image), it is unknown whether $X$ must be a hyperelliptic curve. 
As mentioned in \S~\ref{sec:torsion-Ceresa}, there are instances of nonhyperelliptic curves whose Ceresa class is torsion (respectively, whose Ceresa cycle is torsion modulo algebraic equivalence).

Similar to the algebraic setting, if $\Gamma$ is a hyperelliptic tropical curve (in the sense of \cite{Baker, Chan-Hyperelliptic}) and $\flat \in \Gamma$ is a Weierstrass point, then  the tropical Ceresa cycle $\nu_{\flat}(\Gamma)$ equals $0$ as a cycle. 
Motivated by the discussion in the previous paragraph and the explicit description of the tropical Ceresa class in Theorems \ref{thm:ceresa-explicit} and \ref{thm:unpointedCeresaFormula}, we ask the following question. 

\begin{question}
	Is there a nonhyperelliptic tropical curve $\Gamma$ such that either the pointed $\cc_{\flat}(\Gamma)$ or the unpointed tropical Ceresa class $\overline{\cc}(\Gamma)$ is 0?
\end{question}

In the hyperelliptic case, we have the following.
\begin{proposition}
	Suppose $\Gamma$ is a hyperelliptic tropical curve, and let $\flat\in \Gamma$. If $\cc_{\flat}(\Gamma) = 0$, then $\flat$ is a Weierstrass point.  
\end{proposition}

\begin{proof} 
We provide a sketch of the proof. The Weierstrass points of a hyperelleptic tropical curve are those points $\flat$ such that $2\flat$ is linearly equivalent to 0. 
Suppose $\rflat$ is a Weierstrass point of $\Gamma$. 
By monodromy invariance of $\omega$ as defined in Equation \eqref{eq:trop-omega}, we have $2\AJ(\flat-\rflat) = 0$, so $2\flat$ is linearly equivalent to $2\rflat$, which is linearly equivalent to 0. 
So $\flat$ is a Weierstrass point, as required. 
\end{proof}

\subsection{Tropical Albanese and Roitman's theorem}

Let $X$ be a complex smooth projective variety of dimension $d$. 
The Albanese of $X$ is the intermediate Jacobian $\Alb(X) = \JH^{2d-1}(X)$.
Roitman's theorem asserts that the Abel--Jacobi map
\begin{equation*}
	\AJ \colon \rmA_{0}^{\circ}(X) \to \Alb(X)
\end{equation*}
induces an isomorphism on torsion points.
Furthermore, if $\rmA_{0}(X)$ is representable, e.g., if $\rmA_{0}(X)$ is finite dimensional, then the Abel--Jacobi map above is an isomorphism. See \cite[Ch.~10]{Voisin:HodgeII} for details. 

Now suppose $X$ is a smooth and compact K\"ahler tropical variety whose monodromy is rational in the sense of \S~\ref{sec:Tropical-Albanese}. 
As we demonstrated in that section, the Albanese of $X$ is a tropical abelian variety. 
Does the analog of Roitman's theorem for such tropical varieties hold? More precisely, we ask the following question.

\begin{question}
	Does the tropical Abel--Jacobi map $\AJ \colon \ChowTriv{0}{X} \to \Alb(X)$ induce an isomorphism on torsion points? 
\end{question}

\appendix

\section{Relation to the Morita class}
\label{sec:MoritaClass}

In \cite{CoreyEllenbergLi}, the authors define and study a related tropical Ceresa class that is closely related to the Morita class in \cite{Morita:Extension-Johnson}. The benefit of this class is that it computes the $\ell$-adic Ceresa class of an algebraic curve $C$ defined over $K=\CCt$. In this appendix, we discuss the relationship between this class and our unpointed tropical Ceresa class $\overline{\cc}(\Gamma)$. 

\subsection{The Johnson homomorphism}
We recall some concepts on mapping class groups, a general reference is \cite{FarbMargalit}. 
Denote by $\ssSigma_g$ the genus $g$ topological surface and let $H = \rmH_{1}(\ssSigma_{g}, \Z)$. 
Throughout we assume that $g\geq 3$. 
The \emph{algebraic intersection pairing} is a skew-symmetric bilinear pairing $H\times H \to \Z$ that computes the signed intersections $\langle a, b \rangle$ of representatives of the homology classes $a$ and $b$. 
This defines a canonical symplectic form $\omega \in \wedge^{2}H$. 

The \emph{mapping class group} of $\ssSigma_g$ is the group $\Mod(\ssSigma_{g})$ of isotopy classes of  orientation-preserving diffeomorphisms from $\ssSigma_g$ to itself. Given a simple closed curve $\gamma$ on $\ssSigma_{g}$, denote by $\rmT_{\gamma} \in \Mod(\ssSigma_{g})$ the (left-handed) Dehn twist about $\gamma$. 
The Torelli group $\ss \cI_{g}$ of $\ssSigma_g$ is the normal subgroup of $\Mod(\ssSigma_{g})$ consisting of those mapping classes $f\colon \ssSigma_g \to \ssSigma_g$ such that $f_{*}\colon H\to H$ is the identity. 
For $g\geq 3$, the Torelli group is generated by bounding pair maps \cite{Birman, Powell}. 
A \emph{bounding pair} is a pair of nonseparating, nonintersecting simple closed curves $\ell$ and $\gamma$ that have the same homology class, and a \emph{bounding pair map} is a mapping class of the form $\rmT_{\ell}\rmT_{\gamma}^{-1}$ for a bounding pair $(\ell, \gamma)$.   

Up to 2-torsion, the abelianization of the Torelli group is isomorphic to $(\wedge^{3}H) / H$ where  $H\subset \wedge^3 H$ via the embedding $h\mapsto \omega \wedge h$ \cite{JohnsonIII}. The isomorphism is induced by the Johnson homomorphism $J\colon \ss \cI_{g} \to (\wedge^{3} H) / H$, which we describe on bounding pair maps. Suppose $(\ell,\gamma)$ is a bounding pair. Cutting $\ssSigma_g$ along $\ell$ and $\gamma$ separates $\ssSigma_g$ into two connected surfaces; denote by $S$ the one on the left of $\ell$. The homology $\rmH_1(S,\Z)$ contains a maximal symplectic subspace; denote the corresponding symplectic form by $\omega_S$.  Then 
\begin{equation}
\label{eq:Johnson}
	J(\rmT_{\ell}\rmT_{\gamma}^{-1}) = \omega_{S} \wedge [\ell]
\end{equation}
where $[\ell]$ denotes the homology class of $\ell$.

Morita in \cite{Morita:Extension-Johnson} defines an extension of $2J$ to a cocycle $m\colon \Mod(\ssSigma_g) \to (\wedge^3H)/H$. We describe a method of obtaining such an extension that agrees with Morita's construction at the level of group cohomology. Let $\tau \in \Mod(\ssSigma_g)$ be a hyperelliptic involution of $\ssSigma_g$, i.e., a mapping class such that $\tau_{*}\colon H\to H$ equals $-I$ and $\tau^{2} = 1\in \Mod(\ssSigma_g)$. Define
\begin{equation*}
	m_{\tau}\colon \Mod(\ssSigma_g) \to (\wedge^3H)/H \quad f\mapsto J([f,\tau])
\end{equation*}
where $[f, \tau] = f \tau f^{-1} \tau^{-1}$ is the commutator. This is a group cocycle, its restriction to $\ss \cI_{g}$ is $2J$, its class in $\rmH^{1}(\Mod(\ssSigma_g), (\wedge^3H)/H)$ is independent of $\tau$.

\subsection{The Morita class of a tropical curve}
Next, we recall the construction in \cite{CoreyEllenbergLi}.  Suppose $\Lambda$ is a set of distinct isotopy classes of simple closed curves that have trivial pairwise geometric intersection. Furthermore, assume that each connected component of $\ssSigma_g \setminus \bigcup_{\gamma \in \Lambda} \gamma$ has genus equal to 0. The \emph{dual graph} of $\Lambda$ is a graph with
\begin{enumerate}
	\item a vertex $v_S$ for the closure $S$ of each connected component $\ssSigma_g \setminus \bigcup_{\gamma \in \Lambda} \gamma$, and
	\item an edge $e_{\gamma}$ connecting $v_{S}$ and $v_{S'}$ for each loop $\ell$ in the boundary of  $S$ and $S'$ (possibly  $S = S'$).
\end{enumerate}
Let $\Gamma$ be a tropical curve with a model $(G=(V,E),\ell)$ and  $\Lambda(G) = \{\gamma_e \, \colon  \, e\in E\}$  a fixed arrangement of pairwise nonintersecting isotopy classes of simple closed curves whose dual graph is $G$.  Define the multitwist $\rmT_{\Gamma} \in \Mod(\ssSigma_g)$ by
\begin{equation*}
	\rmT_{\Gamma} = \prod_{e \in E} \rmT_{\gamma_e}^{\ell(e)}.
\end{equation*}
Two Dehn twists $\rmT_{\gamma}$ and $\rmT_{\gamma'}$ commute if $\gamma$ and $\gamma'$ have trivial geometric intersection, so the product above is unambiguous. The \emph{Morita class}\footnote{This class in \cite{CoreyEllenbergLi} is called the tropical Ceresa class. We use the term Morita class to distinguish this class from the tropical Ceresa class defined in this paper. }
 $m(\Gamma)$ of $\Gamma$ is the  cohomology class
\begin{equation*}
	m(\Gamma) \in \rmH^1(\langle \rmT_{\Gamma} \rangle, (\wedge^3H) / H) 
	\cong \rmH^1(\Z, (\wedge^3H) / H) 
\end{equation*}
given by the restriction of $m_{\tau}$ to the infinite cyclic subgroup $\langle \rmT_{\Gamma} \rangle$ of $\Mod(\ssSigma_g)$. A different choice of arrangement $\Lambda(G)$ produces a multitwist conjugate to $\rmT_{\Gamma}$, and hence the cohomology class $m(\Gamma)$ is unaffected. 

There is a nice expression for the image of the Morita class $m(\Gamma)$ in a particular subquotient $\overline{B}(\ssdelta_{\Gamma})$, which we now describe. Denote by $\beta_e\in H$ the homology class of $\gamma_e$ and let $Y = \Span\{\beta_e \, \colon  \,  e\in E \} \subset H$. Because each connected component of $\ssSigma_g \setminus \bigcup_{e \in E} \gamma_e$ has genus equal to 0, the subgroup $Y$ defines a (saturated) Lagrangian subgroup of the symplectic space $H$.  Given a spanning tree $T=(V,F)$ of $G$, there are $g$ edges in $F^{c}$, and the corresponding homology classes $\beta_{\ve}$ form a basis of $Y$. This extends to a symplectic basis $\{\alpha_\ve,\, \beta_\ve \st  \ve \in F^{c}\}$ where $\langle \alpha_e, \beta_{\ve} \rangle = \delta_{e\ve}$ for $e,\ve\in F^{c}$. 

Denote by $\ssdelta_{\Gamma}$ the symplectic representation of $\rmT_{\Gamma}$. By the formula for symplectic representation of Dehn twists:
\begin{equation*}
	(\rmT_{\gamma})_{*}([\zeta]) = [\zeta] + \langle [\zeta], [\gamma] \rangle [\gamma]
\end{equation*}
we have that the matrix of $\ssdelta_{\Gamma}$ with respect to the basis $\{\alpha_e,\, \beta_e \st e\in F^{c}\}$ is
\begin{equation*}
	\ssdelta_{\Gamma} = 
	\begin{bmatrix}
		\ss I_{g} & 0 \\
		\ssQ_{\Gamma} & \ss I_{g}
	\end{bmatrix}
\end{equation*}
where $\ssQ_{\Gamma}$ is the polarization matrix of $\Gamma$ as in \S~\ref{sec:curvesJacobians}. In particular $(\ssdelta_{\Gamma}-I)(H) \subset Y$.  The induced map $\ssdelta_{\Gamma}\colon  \wedge^3 H \to \wedge^3 H$ satisfies
{\footnotesize
\begin{align*}
	&(\ssdelta_{\Gamma} - I)(h_1\wedge h_2 \wedge h_3) 
	= (\ssdelta_{\Gamma} - I) (h_1) \wedge h_2 \wedge h_3 +  h_1 \wedge (\ssdelta_{\Gamma} - I) (h_2) \wedge h_3 +  h_1 \wedge h_2 \wedge (\ssdelta_{\Gamma} - I)(h_3) \\
	&+(\ssdelta_{\Gamma} - I) (h_1) \wedge (\ssdelta_{\Gamma} - I)(h_2) \wedge h_3 +  (\ssdelta_{\Gamma} - I)(h_1) \wedge  h_2 \wedge (\ssdelta_{\Gamma} - I)(h_3) +  h_1 \wedge (\ssdelta_{\Gamma} - I)(h_2) \wedge (\ssdelta_{\Gamma} - I)(h_3) \\ 
	&+(\ssdelta_{\Gamma} - I) (h_1) \wedge (\ssdelta_{\Gamma} - I)(h_2) \wedge (\ssdelta_{\Gamma} - I)(h_3).
\end{align*}
}
The operator $\ssdelta_{\Gamma}-I$ on $(\wedge^{3}H)/H $ is nilpotent, and so we may consider its corresponding weight filtration $W_{\bul}$ (see, \eg, \cite[\S~7]{Hain08}):
\begin{gather*}
	W_{-6} = \wedge^3Y, \quad 
	W_{-4} = (H \wedge Y \wedge Y) / Y \quad 
	W_{-2} = (H \wedge H \wedge Y) / H \quad
	W_{0} = (\wedge^3H)/ H, \\
	W_{-5} = W_{-6}, \quad
	W_{-3} = W_{-4}, \quad 
	W_{-1} = W_{-2}.
\end{gather*}
The $k$-th graded piece of this filtration is $\Gr_{k}(W_{\bul}) = W_{k}/W_{k-1}$. By the above formula for $\ssdelta_{\Gamma}$ and the fact that $(\ssdelta_{\Gamma}-I)(H) \subset Y$, we have that $(\ssdelta_{\Gamma}-I)(W_k) \subset W_{k-2}$, and so we have induced maps
\begin{equation*}
	\ssdelta_{\Gamma}-I \colon  \Gr_{k}(W_{\bul}) \to \Gr_{k-2}(W_{\bul}).
\end{equation*}
When $k=1$, this map is an isomorphism after tensoring with $\Q$. Explicitly,
\begin{equation}	
\label{eq:deltaI}
	(\ssdelta_{\Gamma}-I)(h_1\wedge h_2 \wedge y) = \ssQ_{\Gamma}(h_1) \wedge h_2 \wedge y + h_1 \wedge \ssQ_{\Gamma}(h_2) \wedge y
\end{equation}
and
{\small
\begin{equation*}
	(\ssdelta_{\Gamma}-I)^{-1}(h\wedge y_1 \wedge y_2) = 
	\tfrac{1}{2}(h \wedge \ssQ_{\Gamma}^{-1}(y_1) \wedge y_2  + h \wedge y_1 \wedge \ssQ_{\Gamma}^{-1}(y_2) -  \ssQ_{\Gamma}(h) \wedge \ssQ_{\Gamma}^{-1}(y_1) \wedge \ssQ_{\Gamma}^{-1}(y_2)).
\end{equation*}
}
Define
\begin{align*}
	\overline{A}(\ssdelta_{\Gamma}) &= \image(\rmH^1(\langle \rmT_{\Gamma} \rangle, (\wedge^3H) / H) \to \rmH^1(\langle \rmT_{\Gamma} \rangle, W_{-4}) \cong 
	W_{-4} / (\ssdelta_{\Gamma}-I)W_{-2}, \\
	\overline{B}(\ssdelta_{\Gamma}) &= \coker(\Gr_{-2}(W_{\bul}) \to \Gr_{-4}(W_{\bul})) \cong W_{-4} /( (\ssdelta_{\Gamma}-I)W_{-2} + W_{-6}). 
\end{align*}
Because $(\ssdelta_{\Gamma}-I)\colon \Gr_{-2}(W_{\bul}) \to \Gr_{-4}(W_{\bul})$ is rationally surjective, the groups $\overline{A}(\ssdelta_{\Gamma})$ and $\overline{B}(\ssdelta_{\Gamma})$ are finite. 
\noindent A main result of \cite{CoreyEllenbergLi} is that $m(\Gamma) \in \overline{A}(\ssdelta_{\Gamma})$, and is  therefore torsion. Quotienting by $W_{-6}$ yields a surjective homomorphism $\overline{A}(\ssdelta_{\Gamma}) \to \overline{B}(\ssdelta_{\Gamma})$. The image of the Morita class $m(\Gamma)$ in $\overline{B}(\ssdelta_{\Gamma})$, which we denote by $n(\Gamma)$, is given by
\begin{equation}
\label{eq:formula-morita}
	n(\Gamma) = \sum_{e\in E} J(\rmT_{\gamma_e} \rmT_{\tau(\gamma_e)}^{-1}).
\end{equation}

We end this section by describing a hyperelliptic involution $\tau$ that will be useful in the proof of Theorem \ref{thm:Morita-to-Ceresa}. 

\begin{proposition}
\label{prop:hyperelliptic-basis}
	Given a spanning tree $T=(V,F)$ of $G$, there is a hyperelliptic involution $\tau \in \Mod(\ssSigma_g)$ such that 
	\begin{equation*}
		\tau(\gamma_{\ve}) = -\gamma_{\ve}
	\end{equation*}
	for $\ve\in F^{c}$. That is, $\tau$ takes $\gamma_{\ve}$ to itself, but reverses its orientation. 
\end{proposition}

\begin{proof}
	Let $S$ be the surface obtained by cutting along $\gamma_{\ve}$ for $\ve \in F^{c}$. The surface $S$ is connected, has genus $0$, and $2g$ boundary components.  Choose a simple closed curve $\vartheta$ which separates the two copies $\gamma_{\ve,1}$ and $\gamma_{\ve,2}$ of $\gamma_{\ve}$ in $S$. Cutting along $\vartheta$ results in two homeomorphic surfaces $S_1, S_2$ of genus 0 and with $g+1$ boundary components with the $\gamma_{\ve, i}$'s in $S_{i}$. Choose a homeomorphism $\phi\colon S_1 \to S_2$ which sends $\gamma_{\ve,1}$ to $\gamma_{\ve,2}$. The homeomorphisms $\phi$ and $\phi^{-1}$ define a hyperelliptic involution $\tau$ on $\ssSigma_g$ with the required property. 
\end{proof}

For example, the curves $\beta_1$, $\ldots$, $\beta_5$ in Figure \ref{fig:sgn-surface} satisfy $\tau(\beta_i) = - \beta_i$ where $\tau$ is the hyperelliptic involution given by rotating the surface $180^{\circ}$ about the horizontal axis. 

\subsection{Comparison between the Morita class and the tropical Ceresa class}
\label{app:Comparison}

In this section, we provide the relationship between the Morita class $n(\Gamma)$ and the tropical Ceresa class $\overline{\cc}(\Gamma)$. First, we relate $\overline{B}(\ssdelta_{\Gamma})$ to the primitive intermediate Jacobian $\overline{\JH}_{2,1}(\Jac(\Gamma))$. 

\begin{proposition}
\label{prop:comparisonTopological}
We have an embedding
\begin{equation*}
		\Phi_{\Gamma}\colon  \overline{B}(\ssdelta_{\Gamma}) \to \overline{\JH}_{2,1}(\Jac(\Gamma)), \quad \alpha_{\ve_1} \wedge \beta_{\ve_2} \wedge \beta_{\ve_3} \mapsto [\cycb{\ve_1}, \unt{\ve_2} \wedge \unt{\ve_3}].
	\end{equation*}
\end{proposition}
\noindent See \S~\ref{sec:Ceresa-class-explicit} for the definitions of $\cycb{\ve}$ and $\unt{\ve}$.

\begin{proof}[Proof of Proposition \ref{prop:comparisonTopological}]
Let $X = \Span\{\alpha_\ve \st \ve \in F^{c}\}$. The maps $\alpha_{\ve} \mapsto \cycb{\ve}$ and $\beta_{\ve} \mapsto \unt{\ve}$ define isomorphisms
\begin{equation}
\label{eq:XY-tropical-homology}
X\wedge X \otimes  Y \cong \rmH_{1,2}(\Jac(\Gamma),\Z) \quad \text{and} \quad  X\otimes Y\wedge Y \cong \rmH_{2,1}(\Jac(\Gamma),\Z).
\end{equation} 
See Equation \eqref{eq:tropicalHomologyTropAbelianVariety}.
Set
\begin{equation*}
	B(\ssdelta_{\Gamma}) = \frac{X \otimes Y \wedge Y}{(\ssdelta_{\Gamma}-I)(X\wedge X \otimes Y)}.
\end{equation*}
Note that $\overline{B}(\ssdelta_{\Gamma}) \cong B(\ssdelta_{\Gamma}) / \omega \wedge Y$. By the description of $\ssdelta_{\Gamma}-I$ in Equation \eqref{eq:deltaI} and $\rmN$ in Propositions \ref{prop:monodromyAbelianVariety} and \ref{prop:monodromyAndIntegrtion},  we have a commutative diagram
\begin{equation*}
	\begin{tikzcd}
	X\wedge X \otimes Y \arrow[r, "\ssdelta_{\Gamma}-I"]  \arrow[d] & 
	X\otimes Y\wedge Y \arrow[r] \arrow[d] & B(\ssdelta_{\Gamma}) \arrow[d] \arrow[r] & 0\\
	\rmH_{1,2}(\Jac(\Gamma),\Z) \arrow[r, "\rmN"]  & \rmH_{2,1}(\Jac(\Gamma),\R) \arrow[r]  & \JH_{2,1}(\Jac(\Gamma)) \arrow[r] & 0
	\end{tikzcd} 
\end{equation*}
where the rows are exact. The left and middle vertical arrows are the isomorphisms in Equation \eqref{eq:XY-tropical-homology}. The right map produces the desired embedding $\Phi_{\Gamma}$. 
\end{proof}

We conjecture a formula (Conjecture \ref{conj:Johnson}) for the Johnson homomorphism that could be of independent interest. Using this, under the comparison morphism $\Phi_{\Gamma}$ from Proposition \ref{prop:comparisonTopological}, we have the following theorem. 

\begin{theorem}
\label{thm:Morita-to-Ceresa}
	Let $\Gamma$ be a tropical curve with integral edge lengths. 
	If Conjecture \ref{conj:Johnson} holds, then we have
	\begin{equation*}
		\Phi_{\Gamma}(n(\Gamma)) = \overline{\cc}(\Gamma).
	\end{equation*}
\end{theorem}

We now describe our conjectural formula for the Johnson homomorphism. 
Suppose there are nonseparating simple curves $\gamma, \gamma'$ and $\beta_1,\ldots,\beta_{g}$ such that 
\begin{enumerate}
	\item $\gamma$ and $\gamma'$ are homologous and disjoint from each $\beta_i$, and
	\item the curves $\beta_1,\ldots, \beta_g$  are pairwise disjoint and their homology classes form a basis for a Lagrangian subspace of $\rmH_1(\ssSigma_g, \Z)$.
\end{enumerate}

Let $S$ be the surface obtained by cutting along $\beta_1, \ldots, \beta_g$. 
This is a genus $0$ surface with $2g$ boundary components. 
The boundary components come in pairs corresponding to the curves $\beta_1, \ldots,\beta_g$.
 The curves $\gamma$ and $\gamma'$ divide $S$ into four types of regions, depending on whether the region lies to the left or to the right of $\gamma$ and $\gamma'$, respectively. 
 (Without loss of generality, we may assume that $\gamma$ and $\gamma'$ intersect transversally.)
 We say these regions are of type $LL$, $LR$, $RL$ and $RR$ (e.g., $LR$ means the region is to the \emph{left} of $\gamma$ and to the \emph{right} of $\gamma'$).

The curve $\beta_i$ is \emph{split} by $\gamma$ (respectively, $\gamma'$) if exactly one of the two copies of $\beta_i$ in $S$ is to the left of $\gamma$ (respectively, $\gamma'$). 
Otherwise, $\beta_i$ is \emph{nonsplit}. 
Since $\gamma$ and $\gamma'$ are homologous, we have that $\beta_i$ is split by $\gamma$ if and only if it is split by $\gamma'$. 
Furthermore, the copy of a split $\beta_{i}$ on the left of $\gamma$ also lies to the left of $\gamma'$. 
This implies that the split pairs only appear in the regions of types $LL$ and $RR$. 
Define a sign function  $\sgn_{\gamma,\gamma'} \colon \{\beta_1,\ldots, \beta_g\} \to \{0,\pm 1\}$ by
\begin{equation*}
	\sgn_{\gamma,\gamma'}(\beta_i) =
	\begin{cases}
		-1 & \text{ if } \beta_i \text{ is of type } LR, \\
		1 & \text{ if } \beta_i \text{ is of type } RL, \\
		0 & \text { otherwise.}
	\end{cases}
\end{equation*}
See Figure \ref{fig:sgn-surface}. Finally, suppose $\alpha_1, \ldots, \alpha_g$ are curves such that the homology classes $[\alpha_1]$, $\ldots$, $[\alpha_g]$, $[\beta_1]$, $\ldots$, $[\beta_g]$ form a symplectic basis of $\rmH_{1}(\ssSigma_g, \Z)$. 

\begin{conjecture}
\label{conj:Johnson}
The image of the Johnson homomorphism on the mapping class $\rmT_{\gamma} \rmT_{\gamma'}^{-1}$ is
\begin{equation*}
	\sum_{i=1}^{g} \sgn_{\gamma,\gamma'}(\beta_i) \, [\alpha_{i}] \wedge [\beta_{i}] \wedge [\gamma].
\end{equation*}
\end{conjecture}

\begin{proof}[Proof of Theorem \ref{thm:Morita-to-Ceresa}]
	There are $g$ curves in $\{\gamma_{\ve} \st \ve \in F^{c}\}$, and their homology classes form a basis for a Lagrangian subspace of $H$. Enumerate these as $\beta_1,\ldots,\beta_g$. Let $\alpha_1,\ldots,\alpha_g$ be curves so that the homology classes of the $\alpha_i$, $\beta_i$'s form a symplectic basis of $H$. By Proposition \ref{prop:hyperelliptic-basis}, each pair $(\gamma_{e},\tau(\gamma_e))$, together with the homology basis above, satisfies the hypotheses of Conjecture \ref{conj:Johnson}. So the theorem follows from this conjecture, the formula for $n(\Gamma)$ in Equation \eqref{eq:formula-morita}, and the formula for $\overline{\cc}(\Gamma)$ in Theorem \ref{thm:unpointedCeresaFormula}. 
\end{proof}

We describe a situation where we can prove Conjecture \ref{conj:Johnson}. 
Suppose there is a sequence of curves $\gamma = \gamma_1,\gamma_2,\ldots,\gamma_k = \gamma'$ such that
	\begin{itemize}
		\item each $\gamma_i$ is disjoint from $\beta_1,\ldots,\beta_g$, 
		\item $[\gamma_i] = [\gamma]$, and 
		\item $\gamma_i$ and $\gamma_{i+1}$ are disjoint. 
	\end{itemize} 
	We claim that Conjecture \ref{conj:Johnson} holds in this case. 

Before proving this claim, we comment about these conditions. 
As demonstrated in \cite[Thm~1.9]{Putman-08}, one can always find a sequence $\gamma_1,\ldots,\gamma_k$ satisfying the second and third conditions, but as communicated to us by  by Andrew Putman, it may not be possible to ensure the first condition. 
Nevertheless, we believe that Conjecture \ref{conj:Johnson} still holds. 
In our proof of Theorem \ref{thm:Morita-to-Ceresa}, we rely on the validity of this conjecture in the case $(\gamma,\gamma') = (\gamma_{e}, \tau(\gamma_{e}))$ for $e\in E$ for some hyperelliptic involution $\tau$. 

Let us prove the claim we made above. We proceed by induction on $k$. When  $k=2$, the curves $\gamma$ and $\gamma'
	$ are disjoint. 
Then the formula above follows from the description of the Johnson homomorphism in Equation \eqref{eq:Johnson}. 
	
	For the inductive step, assume that $\vartheta$ is a simple closed curve homologous to $\gamma$ such that $\vartheta$ and $\gamma'$ are disjoint, and the formula holds for the pair $(\gamma, \vartheta)$. We have
	\begin{align*}
		J(\rmT_{\gamma}\rmT_{\gamma'}^{-1}) &= 
		J(\rmT_{\gamma}\rmT_{\vartheta}^{-1}) 
		+ J(\rmT_{\vartheta}\rmT_{\gamma'}^{-1}) \\
		&=  \sum_{i=1}^{g} \sgn_{\gamma,\vartheta}(\beta_{i}) \,[\alpha_{i}] \wedge [\beta_{i}] \wedge [\gamma] 
		+  \sum_{i=1}^{g} \sgn_{\vartheta,\gamma'}(\beta_i) \, [\alpha_{i}] \wedge [\beta_{i}] \wedge [\gamma].
	\end{align*}
	This equals
	\begin{align*}
		 \sum_{i=1}^{g} \sgn_{\gamma,\gamma'}(\beta_i) \, [\alpha_{i}] \wedge [\beta_{i}] \wedge [\gamma]. 
	\end{align*}
	because of the relation
	\begin{equation*}
		\sgn_{\gamma,\vartheta}(\beta_i) + \sgn_{\vartheta,\gamma'}(\beta_i)  = \sgn_{\gamma,\gamma'}(\beta_i) 
	\end{equation*}
	which readily follows from the definition of the sign function.

\bibliographystyle{alpha}
\bibliography{biblio}

\end{document}